\documentclass[final]{siamltex}
\usepackage{amsmath, amssymb, amsfonts}
\usepackage{graphicx,color}
\usepackage{diagbox}
\usepackage{algorithm}
\usepackage{algpseudocode}
\usepackage{enumitem}
\usepackage{multirow}
\usepackage{subfig}
\usepackage{epstopdf}
\usepackage{comment}
\usepackage{braket}
\usepackage{caption}
\usepackage[dvipsnames]{xcolor}

\usepackage{adjustbox}

\newcommand{\Tr}{\text{Tr}}

\newcommand{\mc}[1]{\mathcal{#1}}

\global\long\def\R{\mathbb{R}}

\global\long\def\ra{\rightarrow}

\global\long\def\Tr{\text{Tr}}

\numberwithin{equation}{section}
\numberwithin{figure}{section}

\providecommand{\corollaryname}{Corollary}
\providecommand{\lemmaname}{Lemma}
\providecommand{\propositionname}{Proposition}
\providecommand{\remarkname}{Remark}
\providecommand{\theoremname}{Theorem}

\allowdisplaybreaks

\title{Multiscale semidefinite programming approach to positioning problems with pairwise structure}

\author{Yian Chen\thanks{Department of Statistics, University of Chicago, Illinois, IL 60637, USA. Email: {\tt yianc@uchicago.edu}}
\and
Yuehaw Khoo\thanks{Department of Statistics, University of Chicago, Illinois, IL 60637, USA. Email: {\tt ykhoo@uchicago.edu}}
\and
Michael Lindsey\thanks{Department of Mathematics, Courant Institute of Mathematical Sciences, New York University, New York, NY 10012, USA.  Email: \texttt{michael.lindsey@cims.nyu.edu}}
}

\begin{document}

\maketitle

\begin{abstract}
 We consider the optimization of pairwise objective functions, i.e., objective functions of the form $H(\mathbf{x}) = H(x_1,\ldots,x_N) = \sum_{1\leq i<j \leq N} H_{ij}(x_i,x_j)$ for $x_i$ in some continuous state spaces $\mathcal{X}_i$. Global optimization in this setting is generally confounded by the possible existence of spurious local minima and the impossibility of global search due to the curse of dimensionality. In this paper, we approach such problems via convex relaxation of the marginal polytope considered in graphical modeling, proceeding in a multiscale fashion which exploits the smoothness of the cost function. We show theoretically that, compared with existing methods, such an approach is advantageous even in simple settings for sensor network localization (SNL). We successfully apply our method to SNL problems, particularly difficult instances with high noise. We also validate performance on the optimization of the Lennard-Jones potential, which is plagued by the existence of many near-optimal configurations. We demonstrate that in MMR allows us to effectively explore these configurations.
\end{abstract}

\section{Introduction}
Determining the optimal coordinates $\mathbf{x} = (x_1,\ldots,x_N) \in \mc{S} := \bigoplus_{i=1}^N \mc{X}_i$ for a pairwise potential $H(\mathbf{x}) = \sum_{1\leq i<j \leq N} H_{ij}(x_i,x_j)$ is a problem that finds widespread application in engineering and the physical sciences. Most black-box methods, e.g., first- and second-order optimization and simulated annealing~\cite{van1987simulated}, rely only on local information to make updates and hence are liable to get stuck in local minima. It is natural to wonder whether a general and practical optimizer can make use of global information, evading the curse of dimensionality via the pairwise structure of the objective. In spite of the NP-hardness of optimizing general pairwise objectives~\cite{murty1985some}, such an optimizer may avoid the pitfalls of local optimization and could be hoped to succeed on important practical problems. Although progress along these lines has been made in the context of discrete (and especially 0-1) optimization (for example the celebrated 
semidefinite program~\cite{goemans1994879} for the max-cut problem), one nonetheless hopes to expand the class of soluble optimization problems. In particular, for continuous optimization problems that do not admit tractable moment-based relaxations (unlike, e.g., polynomial optimization~\cite{lasserre2001global,parrilo2003semidefinite}), to our knowledge there is a scarcity of general-purpose methods that are fundamentally global. Progress on this front could enhance the understanding of the energy landscape of molecular-dynamical models, which are generally defined by pairwise potentials and which play central roles in quantitative chemistry and biology.

As a model problem of molecular-dynamical type, we study the Lennard-Jones cluster~\cite{wales1997global} in two dimensions. The energy landscape of this model features many near-optimal local minima. It is believed that the energy landscape is `funnel-shaped'~\cite{doye1999double} in that by hopping between nearby local minima and making improvements in the energy, one can successfully find the global optimum. This perspective underlies the basin-hopping algorithm~\cite{wales1997global}, which has been successful for this type of problem. In addition to finding the global optimum, it is of interest to sample from the near-global optima. We will demonstrate that the methodology proposed in this paper achieves competitive performance on these tasks.

We shall also study sensor network localization (SNL) problem, in which one seeks to recover the absolute positions of sensors in the plane, given (possibly corrupted) measurements of their pairwise distances from one another. Semidefinite relaxations~\cite{so2007theory,nie2009sum,biswas2006semidefinite} have been proposed for this type of problem, but they are moment-based relaxations of an objective function that is not ideally suited to the random corruption model that we consider. As such it is of interest to consider a more flexible relaxation framework that can accommodate more general objectives. Moreover, with the exception of~\cite{nie2009sum}, these relaxations cannot localize certain instances (detailed in Section~\ref{section:marginal relaxation}) where the solution is unique.

Our approach begins by reformulating a general optimization problem over $\mathcal{S}$ as a linear program over probability measures $\mu = \mu(\textbf{x})$ on $\mathcal{S}$ with objective function $\mu \mapsto \langle H, \mu \rangle = \sum_{\mathbf{x} \in \mathcal{S}} H(\mathbf{x}) \mu(\mathbf{x})$. The optimizer of this linear program is the Dirac mass at the optimizer $\mathbf{x}^\star \in \mathcal{S}$ of the original optimization program (assuming that it is unique). Unfortunately, the dimension of the optimization space for the linear program is exponentially large in $N$, so this problem is intractable to solve directly.

However, the perspective of measure optimization (as adapted by~\cite{lasserre2001global, parrilo2003semidefinite}) allows for the key observation that the objective function can be rewritten in terms of the marginal distributions $\mu_{ij} = \mu_{ij}(x_i,x_j)$ (defined concretely in~\eqref{eq:2mar} below) as $\sum_{i<j} \langle H_{ij}, \mu_{ij} \rangle$. This allows us to apply marginal polytope relaxations in the style of variational inference (see~\cite{wainwright2008graphical} and references therein) to approximately optimize over the collection of 2-marginals $\mu_{ij}$. We consider a tightening of the local marginal polytope relaxation~\cite{peng2012approximate} that includes additional semidefinite constraints satisfied by the 2-marginals. This semidefinite relaxation, which we call the \emph{2-marginal relaxation}, was considered in the context of multi-marginal optimal transport in~\cite{khoo2019convex}. Under certain circumstances, one can guarantee exact recovery \emph{a posteriori} of the global minimizer. Otherwise, it is possible to define a natural guess $\mathbf{x} \in \mathcal{S}$ for the global minimizer in terms of the solution of the 2-marginal relaxation.

Now the 2-marginal relaxation requires a discrete state space to be numerically soluble. Thus for continuous optimization problems, we are forced to discretize the local state spaces $\mathcal{X}_i$ before we can apply the 2-marginal relaxation. Since the SDP becomes intractable as the discretization is refined, we pursue a multiscale approach, in which we perform global optimization on a coarsely discretized space via the 2-marginal relaxation and use the solution to define a bounding box for the global optimizer of the original problem. Then we refine our discretization on this bounding box and repeat until we locate the global optimizer. We call the algorithmic realization of this approach \emph{multiscale marginal relaxation} (MMR). We explain how for problems with many optimal or near-optimal solutions, MMR can also succeed in exploring these configurations.

In the problems considered, it is common that the solutions are trivially non-unique. For example, an arbitrary rigid motion applied to all the variables $x_1,\ldots,x_N$ may not change the cost. We explain how to remove such an ambiguity, which is crucial to the success of our algorithm. We also explain how permutation-invariance in the cost (which also introduces degeneracy) can in fact be exploited for computational efficiency.

\subsection{Outline}
In Section~\ref{section:marginal relaxation}, we detail the marginal relaxation and motivate its use by theoretically showing its advantage for certain simple SNL problems.  In Section \ref{sec:symmetricGeneral}, we also explain how permutation-symmetry among particles can be exploited for computational efficiency. In Section~\ref{sec:MMR}, we develop a multiscale approach to make the marginal relaxation tractable for continuous problems. In Section~\ref{sec:numerical}, we numerically verify the effectiveness of the proposed method for SNL problems and global optimization of the Lennard-Jones potential.

\section{Marginal relaxation for pairwise objectives}\label{section:marginal relaxation}
We consider global optimization on a state space with product structure, written $\mathcal{S} := \bigoplus_{i=1}^N \mathcal{X}_i$. For simplicity we assume that the `local state spaces' $\mathcal{X}_i$ are finite sets. In the sequel we shall specifically consider  $\mathcal{X}_i$ to be discretizations of $\R^d$, but for now we maintain a general perspective. Elements of $\mathcal{S}$ will be written $\mathbf{x} = (x_1,\ldots, x_N)$.
Our objective function $H:\mathcal{S} \rightarrow \R$ is assumed to have pairwise structure, i.e., to be of the form
\begin{equation}
\label{eq:pairwise}
H(\mathbf{x}) = \sum_{1\leq i<j \leq N} H_{ij}(x_i,x_j),
\end{equation}
where $H_{ij}:\mathcal{X}_i \times \mathcal{X}_j \ra \R$. Generalization to objective functions with triplet-wise structure, etc., can be pursued; see, e.g.,~\cite{KhooEtAl2019,LinLindsey2019} for analogous developments. We want to solve the minimization problem
\begin{equation}
\label{eq:problem}
E_0 := \inf_{\mathbf{x} \in \mathcal{S}} H(\mathbf{x}).
\end{equation}

\subsection{Global optimization as linear programming}
The point of departure for our method is the reformulation of the optimization problem~\eqref{eq:problem} as an optimization problem over probability measures on $\mathcal{S}$:
\begin{equation}
\label{eq:lp}
E_0 = \inf_{\mu \in \mathcal{P}(\mathcal{S})} \langle H, \mu \rangle.
\end{equation}
Here $\mathcal{P}(\mathcal{S})$ denotes the set of probability measures on $\mathcal{S}$, which can be identified with probability mass functions $\mu : \mathcal{S} \ra \R$. Accordingly, the angle brackets indicate the $L^2(\mathcal{S})$ inner product, so $\langle H, \mu \rangle$ denotes the expectation of $H$ with respect to the probability measure $\mu$. It is easy to see that in the case where the optimizer $\mathbf{x}^\star$ of~\eqref{eq:problem} is unique, the unique optimizer of~\eqref{eq:lp} is the delta measure $\mu(\mathbf{x}) = \delta_{\mathbf{x},\mathbf{x}^\star}$. When the optimizer of~\eqref{eq:problem} is not unique, the optimizers of~\eqref{eq:lp} are convex combinations of delta measures localized at the optimizers of~\eqref{eq:problem}.

Notice that the reformulated problem~\eqref{eq:lp} is a convex optimization problem and in fact a linear program. However, the optimization state space has been greatly enlarged. Indeed its dimension $\prod_{i=1}^N \vert \mathcal{X}_i\vert$ grows exponentially in $N$. If one thinks of the $\mathcal{X}_i$ as the state spaces for $N$ interacting particles, this means that direct solution of~\eqref{eq:lp} is limited to problems with very few particles. Although direct solution of~\eqref{eq:lp} is usually intractable, this reformulation serves as a springboard for the specification of relaxed problems that are tractable to solve. As we shall see, these problems will in fact be semidefinite programs.

\subsection{The 2-marginal relaxation}
To this end, given a measure $\mu \in \mathcal{P}(\mathcal{S})$, let $\mu_{ij} \in \mathcal{P}(\mathcal{X}_i\times \mc{X}_j)$ denote the \emph{2-marginals} obtained by marginalizing all components except $i,j$. Equivalently:
\begin{equation}
\label{eq:2mar}
\mu_{ij}(x_i,x_j) = \sum_{\mathbf{x}'\in \mathcal{S}\,:\,x_i' = x_i, x_j'=x_j} \mu(\mathbf{x}') .
\end{equation}
Notice that $\langle H, \mu \rangle = \sum_{i<j} \langle H_{ij}, \mu_{ij} \rangle.$ Hence~\eqref{eq:lp} can in turn be reformulated as
\begin{equation}
\label{eq:lp2}
E_{0} =\inf_{\{\mu_{ij}\}_{i<j} \in \mathcal{P}_2(\mathcal{S}) } \ \sum_{i<j} \langle H_{ij},\mu_{ij}\rangle,
\end{equation}
where $\mathcal{P}_2(\mathcal{S})$ is defined to be the set of collections $\{\mu_{ij}\}_{i<j}$ of \emph{jointly representable} 2-marginals, i.e., those collections $\{\mu_{ij}\}_{i<j}$ which can be obtained as the 2-marginals of a single $\mu \in \mathcal{P}(\mathcal{S})$. Although the dimension of the optimization state space of~\eqref{eq:lp2} has been reduced to $\sum_{i<j} \vert \mathcal{X}_i\vert \cdot \vert \mathcal{X}_j\vert$, the complexity of enforcing the joint representability constraint exactly remains intractable. But now that the complexity has been shifted into the constraints, the door has been opened to relaxation.

Specifically, we will write down several necessary conditions satisfied by any jointly representable collection $\{\mu_{ij}\}_{i<j}$ and define a relaxed problem that replaces the joint representability constraint in~\eqref{eq:lp2} with our necessary conditions. Before proceeding, it is convenient to identify the 2-marginals $\mu_{ij}(x_i,x_j)$ and other functions $\mathcal{X}_i \times \mc{X}_j \ra \R$ with matrices of size $\vert \mathcal{X}_i\vert \times \vert \mathcal{X}_j\vert$ via slight 
abuse of notation. Hence we may indicate $\langle H_{ij} , \mu_{ij} \rangle$ alternatively as $\Tr[H_{ij}^\top \mu_{ij}]$. Moreover, in this notation we have that $\mu_{ij} \mathbf{1}_{\vert \mathcal{X}_j \vert \times 1} \in \R^{\vert \mathcal{X}_i \vert}$ is the vector identified with the suitably defined \emph{1-marginal} $\mu_i (x_i)$, which can be viewed as either a function $\mathcal{X}_i \ra \R$ or a vector of length $\vert \mathcal{X}_i \vert$. (Here and throughout we use $\mathbf{1}_{p \times q}$ throughout to denote the $p \times q$ matrix of all ones.) Likewise $\mu_{ij}^\top \mathbf{1}_{\vert \mathcal{X}_i \vert \times 1} \in \R^{\vert \mathcal{X}_j \vert}$ corresponds to $\mu_j$.

Now we enumerate our necessary representability conditions. To get started, observe that the $\mu_{ij}$ must actually be elements of $\mathcal{P}(\mathcal{X}_i \times \mathcal{X}_j)$, i.e., we have the constraints $\mu_{ij} \geq 0$, $\Tr \left[ \mu_{ij}  \mathbf{1}_{\vert \mathcal{X}_i \vert \times \vert \mathcal{X}_j \vert} \right] = 1$ for all $i<j$. Next we have the \emph{local consistency constraints}, which specify that the 2-marginals agree on overlapping 1-marginals. These constraints can be enforced by introducing an extra set of optimization variables $\mu_i$ for the 1-marginals and enforcing $\mu_i = \mu_{ij} \mathbf{1}_{\vert \mathcal{X}_j \vert \times 1}$ and $\mu_j = \mu_{ij}^\top \mathbf{1}_{\vert \mathcal{X}_i \vert \times 1}$ for all $i < j$. Notice that the constraints $\Tr \left[ \mu_{ij}  \mathbf{1}_{\vert \mathcal{X}_i \vert \times \vert \mathcal{X}_j \vert} \right] = 1$ for all $i<j$ can then be equivalently specified by enforcing only $\mu_i^\top \mathbf{1}_{\vert \mc{X}_i\vert \times 1} = 1$ for all $i$.

Next we derive our final set of constraints, the \emph{global semidefinite constraint}. Let $\mu\in \mathcal{P}(\mc{S})$, and consider functions $f_i : \mc{X}_i \ra \R$ for $i=1,\ldots,N$. Abusing notation slightly, we may identify these with functions $\mc{S} \ra \R$ via the identification $f_i (\mathbf{x}) = f_i (x_i)$. Then compute: 
\begin{align*}
0 & \leq\left\langle \left(\sum_{i}f_{i}\right)^{2},\mu\right\rangle \\
 & =\sum_{i}\left\langle f_{i}^{2},\mu_{i}\right\rangle +\sum_{i\ne j}\left\langle f_{i}f_{j},\mu_{ij}\right\rangle \\
 & =\sum_{i} \sum_{x_i \in \mathcal{X}_i} f_{i}(x_{i})^{2}\mu_{i}(x_{i})+\sum_{i\neq j} \sum_{(x_i,x_j) \in \mathcal{X}_i \times \mathcal{X}_j}  f_{i}(x_{i})f_{j}(x_{j})\mu_{ij}(x_{i},x_{j}).
\end{align*}
In summary, we have derived for each collection $(f_i)$ the following linear inequalities satisfied by the 1- and 2-marginals:
\[
\sum_{i} \sum_{x_i \in \mathcal{X}_i} f_{i}(x_{i})^{2}\mu_{i}(x_{i})+\sum_{i\neq j} \sum_{(x_i,x_j) \in \mathcal{X}_i \times \mathcal{X}_j}  f_{i}(x_{i})f_{j}(x_{j})\mu_{ij}(x_{i},x_{j})\geq0.
\]
Identifying the $\mu_{ij}$ with matrices and the $\mu_i$, $f_i$ with vectors, we have equivalently:
\[
\left(\begin{array}{c}
f_{1}\\
f_{2}\\
\vdots\\
f_{N}
\end{array}\right)^{\top}{\left(\begin{array}{cccc}
\mathrm{diag}(\mu_{1}) & \mu_{12} & \cdots & \mu_{1N}\\
\mu_{21} & \mathrm{diag}(\mu_{2}) & \cdots & \mu_{2N}\\
\vdots & \vdots & \ddots & \vdots\\
\mu_{N1} & \mu_{N2} & \cdots & \mathrm{diag}(\mu_{N})
\end{array}\right)}\left(\begin{array}{c}
f_{1}\\
f_{2}\\
\vdots\\
f_{N}
\end{array}\right)\geq0.
\]
Recall that this inequality holds for all choices $(f_i)$, so our infinite collection of linear inequality constraints can be reformulated as the linear matrix inequality $G(\{\mu_i,\mu_{ij}\}) \succeq 0$, where $G$ is the map defined by 
\[
G(\{\mu_i, \mu_{ij}\}) := 
{\left(\begin{array}{cccc}
\mathrm{diag}(\mu_{1}) & \mu_{12} & \cdots & \mu_{1N}\\
\mu_{21} & \mathrm{diag}(\mu_{2}) & \cdots & \mu_{2N}\\
\vdots & \vdots & \ddots & \vdots\\
\mu_{N1} & \mu_{N2} & \cdots & \mathrm{diag}(\mu_{N})
\end{array}\right)}.
\]
Notice that $\mu_{ji} = \mu^\top_{ij}$, so $G$ can be evaluated only in terms of $\mu_{ij}$ for $i<j$. The condition $G(\{\mu_i,\mu_{ij}\}) \succeq 0$ is our global semidefinite constraint. See~\cite{KhooEtAl2019} for an alternate derivation in terms of the extremal points of $\mc{P}(\mc{S})$.

Collecting our necessary conditions, we have derived the following semidefinite relaxation of~\eqref{eq:lp2}, which we refer to as the \emph{2-marginal relaxation} (cf.~\cite{KhooEtAl2019,LinLindsey2019}): 
\begin{align}
\label{eq:sdp}
\underset{\{\mu_{ij}\in\mathcal{P}(\mc{X}_{i}\times \mc{X}_{j})\}_{i<j},\ \{\mu_{i}\in\mathcal{P}(\mc{X}_{i})\}}{\mbox{minimize}}\quad\quad\quad & \sum_{i<j} \Tr [ H_{ij}^\top \mu_{ij}] \\
\mathrm{\mbox{subject to }\quad\quad\quad} & \mu_{ij}\mathbf{1}_{\vert \mc{X}_j\vert \times 1}=\mu_{i},\ \ \mu_{ij}^{\top}\mathbf{1}_{\vert \mc{X}_i\vert \times 1}=\mu_{j},\quad i\neq j \nonumber \\
 & G(\{\mu_{ij}\})\succeq0. \nonumber
\end{align}

Since the constraints of the 2-marginal relaxation~\eqref{eq:sdp} are satisfied by the 1- and 2-marginals of any $\mu \in \mc{P}(\mc{S})$, the optimal value $E_0^{\text{sdp}}$ of~\eqref{eq:sdp} is guaranteed to be a lower bound for $E_0$. Furthermore, we are guaranteed exact recovery if the optimal $\mu_i$ in~\eqref{eq:sdp} are delta-measures, i.e., $\mu_i (x_i) = \delta_{x_i^\star}(x_i)$ for some $x_i^\star$, $i=1,\ldots N$. Indeed, if this is the case, then letting $\mathbf{x}^\star = (x_1^\star, \ldots, x_N^\star)$ we see that $H(\mathbf{x}^\star) = E_0^{\text{sdp}}$, hence $E_0 \leq E_0^{\text{sdp}}$. But we know automatically that $E_0^{\text{sdp}} \leq E_0$, so in fact $E_0 = E_0^{\text{sdp}}$, and moreover $\mathbf{x}^\star$ is an optimizer for~\eqref{eq:problem}. If this condition for exact recovery does not hold, we can still define a natural guess for the optimizer via $x_i^\star := \mathrm{argmax}_{x_i \in \mathcal{X}_i} \mu_i(x_i)$, effectively rounding our SDP solution to a point in the original optimization space $\mathcal{S}$.

\subsection{Symmetric case}\label{sec:symmetricGeneral}
For certain problems of interest (such as the Lennard-Jones potential optimization that we consider in our numerical experiments), we have that $\mathcal{X}_i = \mathcal{X}$ for all $i$, and $H_{ij} = H$ for all $i,j$, and moreover $H$ (viewed as a matrix) is symmetric. In particular the cost is invariant to permutations among the particles. In this case the optimization problem~\eqref{eq:problem} is generically degenerate with $N!$ optimizers due to the permutation symmetry. For such problems one can consider a specialized SDP relaxation as follows. Note that the derivation recovers a relaxation equivalent to the one considered for symmetric multi-marginal optimal transport in~\cite{khoo2019convex}, modulo the fact that the 1-marginals are not specified a priori in our case.

We describe how to derive a global semidefinite constraint in this setting. First note that the linear program~\eqref{eq:lp} can be restricted to an optimization problem over permutation-invariant probability measures. Then for such measures the 2-marginals $\mu_{ij} = \gamma$ are all the same, and moreover they are (as matrices) symmetric. Meanwhile, the 1-marginals $\mu_i = \rho$ are also all the same. Then let $\mu\in \mathcal{P}(\mc{S})$ be permutation-invariant, and consider a general function $f : \mc{X} \ra \R$. We may then define functions $f_i : \mc{S} \ra \R$ via the identification $f_i (\mathbf{x}) = f (x_i)$. Then compute: 
\begin{align*}
0 & \leq\left\langle \left(\sum_{i}f_{i}\right)^{2},\mu\right\rangle \\
 & = N \left\langle f^{2},\rho \right\rangle + N(N-1) \left\langle f f^\top, \gamma \right\rangle \\
 & = N f^\top \left[ \mathrm{diag}(\rho) + (N-1) \gamma \right] f,
\end{align*}
where we have again identified functions of one and two variables with vectors and matrices, respectively. Since the inequality must hold for all $f$, we obtain the constraint 
\begin{equation}
\label{eq:lambda}
\Lambda:= \frac{1}{N}\mathrm{diag}(\rho) + \frac{N-1}{N} \gamma \succeq 0.
\end{equation}
Here the seemingly arbitrary normalization of $\Lambda$ is meant to connect the notation to~\cite{khoo2019convex}. Together with appropriate local consistency constraints, this semidefinite constraint yields the SDP:
\begin{align}
\label{eq:sdp_symm}
\underset{\gamma\in\mathcal{P}(\mc{X}\times \mc{X}),\, \rho\in\mathcal{P}(\mc{X})}{\mbox{minimize}}\quad\quad\quad & \frac{N(N-1)}{2} \Tr [ H \gamma] \\
\mathrm{\mbox{subject to }\quad\quad\quad} & \gamma \mathbf{1}_{\vert \mc{X}\vert \times 1}= \rho, \ \gamma^\top = \gamma, \nonumber \\
 & \mathrm{diag}(\rho) + (N-1) \gamma \succeq 0. \nonumber
\end{align}

In the case of the Lennard-Jones potential to be considered below, we have that $\mathrm{diag}(H) = +\infty$. In this case, we may constrain $\mathrm{diag}(\gamma) = 0$ and replace $\mathrm{diag}(H)$ by zeros in the above. Notice that the size of this SDP is independent of the number of particles $N$.

In the case, e.g., of the Lennard-Jones potential, where no two particles can occupy the same position, one hopes that the 1-marginal $\rho$ obtained by solving~\eqref{eq:sdp_symm} is supported on $N$ distinct elements $x_1^\star, \ldots, x_N^\star$ of $\mathcal{X}$. Then any permutation of these points is a candidate solution for the original problem~\eqref{eq:problem}.

\subsection{Computing sublevel sets} Returning to the linear program~\eqref{eq:lp}, we will describe a modification that permits the computation of sublevel sets, which motivates some practical developments later on. First we may rewrite~\eqref{eq:lp} as follows:
\begin{eqnarray*}
\underset{\mu:\mathcal{S}\ra\R}{\mbox{minimize}} \quad\quad\quad &  & \left\langle H,\mu\right\rangle \\
\mbox{subject to}\quad\quad\quad &  & \left\langle 1,\mu\right\rangle =1\\
 &  & 0\leq\mu\leq1.
\end{eqnarray*}
Then consider the following modification:
\begin{eqnarray}
\label{eq:lpmod}
\underset{\mu:\mathcal{S}\ra\R}{\mbox{minimize}} \quad\quad\quad &  & \left\langle H,\mu\right\rangle \\
\mbox{subject to}\quad\quad\quad &  & \left\langle 1,\mu\right\rangle = t \nonumber \\
 &  & 0\leq\mu\leq1, \nonumber
\end{eqnarray}
where $t\geq 1$. It is easy to verify that the solution to~\eqref{eq:lpmod} is obtained as follows. Let $K = \lceil t \rceil$, and let $\mathbf{x}^{(1)}, \ldots , \mathbf{x}^{(K)}$ be such that $H(\mathbf{x}^{(k)}) \leq H(\mathbf{x})$ for all $k=1,\ldots,K$ and all $\mathbf{x} \notin \{ \mathbf{x}^{(1)}, \ldots , \mathbf{x}^{(K)} \}$. Roughly speaking, the $\mathbf{x}^{(k)}$ attain the $K$ lowest values of $H$. Then the measure 
\begin{equation}
\label{eq:mixture}
\mu := \sum_{k=1}^{K-1} \delta_{\mathbf{x}^{(k)}} + [1-(t-\lceil t \rceil)] \delta_{\mathbf{x}^{(K)}}
\end{equation}
is an optimizer for~\eqref{eq:lpmod}, and moreover, it is the unique optimizer if there is a gap between the $K$-th and $(K+1)$-th lowest values of $H$, counted with multiplicity. Denoting the $i$-th lowest value (counted with multiplicity) of $H$ by $E_i$, then if $E \in (E_K, E_{K+1})$, we have that the support of the optimizer $\mu^\star$ is the $E$-sublevel set of $H$.

Observe that it is equivalent (by scaling $\mu$) to consider the linear program
\begin{eqnarray}
\label{eq:lpmod2}
\underset{\mu:\mathcal{S}\ra\R}{\mbox{minimize}} \quad\quad\quad &  & \left\langle H,\mu\right\rangle \\
\mbox{subject to}\quad\quad\quad &  & \left\langle 1,\mu\right\rangle = 1 \nonumber \\
 &  & 0\leq\mu\leq b, \nonumber
\end{eqnarray}
where $b = t^{-1}$. This problem~\eqref{eq:lpmod2} carries the interpretation of upper-bounding the amount of probability that can be assigned to any state in the original linear program~\eqref{eq:lp}.

\subsection{Motivating example} Before moving forward, we  demonstrate the advantage of the 2-marginal relaxation~\eqref{eq:sdp} over the popular SNLSDP method~\cite{snlsdp} for sensor network localization in a simple toy setting.

Assume the state space $\mathcal{X}_i = \mathcal{X}$ for each variable is a 1D  grid points such that $|\mathcal{X}|=M$ where $M>N$. Let $x^\star_1,\ldots,x^\star_N \in \mathcal{X}$ be the ground truth positions of $N$ sensors. The observed distance matrix $D\in \mathbb{R}^{N\times N}$ is defined by
\begin{equation}
D(i,j) = \|x^\star_i-x^\star_j\|_2+\epsilon_{ij},\quad i,j=1,\ldots,N, 
\end{equation}
where $\epsilon_{ij}$ is a noise term. We assume in this section that $\epsilon_{ij}=0$.
Assume we  only have distance measurements on a subset of edges $E$ among the $N$ sensors. Then one can attempt to solve the SNL problem via minimizing
\begin{equation} 
\label{eq:cost_exp}
H(\mathbf{x}) = \sum_{(i,j)\in E}(\| x_i - x_j \|_2 - D(i,j))^q
\end{equation} 
for some choice of $q>0$. 

Assume that $E$ is a cycle, consisting of $N$ edges. Further assume that $\mathbf{x}^\star$ is the only $\mathbf{x}$ that minimizes \eqref{eq:cost_exp} (up to trivial global translation and reflection symmetry). We want to show that applying \eqref{eq:sdp} to the minimization problem \eqref{eq:cost_exp}---with some modifications---indeed recovers the solution $\mathbf{x}^\star$. Although this example seems rather trivial, the ability to solve it already shows the advantage of \eqref{eq:sdp} over the popular SDP method (SNLSDP), which often fails to recover $\mathbf{x}^\star$ in such an instance. We note that in the noiseless case where each $\epsilon_{ij}$ is zero, having a cycle in a 1D SNL problem automatically guarantees that \emph{generically} the minimizer of \eqref{eq:cost_exp} (with 0 cost) is unique \cite{jackson2007rigidity}, a property known as generic global rigidity. However, since in our scenario $x_1^\star,\ldots,x_N^\star$ lie on grid points (hence may not always be in generic position)
we need the additional assumption of the uniqueness of $\mathbf{x}^\star$ in minimizing \eqref{eq:cost_exp}.

As mentioned previously, in the case of SNL in 1D, there is a global translation and reflection symmetry, which we break by adding two anchor nodes. Without lost of generality, we can simply let $\mathbf{x}^\star_1 = 0$ and $\mathbf{x}^\star_N = D(1,N)$. The cycle $E$ consists of the edges $(1,2),\ldots,(N-1,N),(N,1)$, so  problem \eqref{eq:sdp} becomes
\begin{align}
\label{eq:sdp localization cycle}
\underset{\{\mu_{ij}\in\mathcal{P}(X_{i}\times X_{j})\}_{i<j},\ \{\mu_{i}\in\mathcal{P}(X_{i})\}}{\mbox{minimize}}\quad\quad\quad & \sum^{N-1}_{i=1} \Tr [ H_{i,i+1}^\top \mu_{i, i+1}] \\
\mathrm{\mbox{subject to }\quad\quad\quad} & \mu_{ij}\mathbf{1}_{\vert \mc{X}_j\vert \times 1}=\mu_{i},\ \ \mu_{ij}^{\top}\mathbf{1}_{\vert \mc{X}_i\vert \times 1}=\mu_{j},\quad i\neq j \nonumber \\
 & G(\{\mu_{ij}\})\succeq0 \label{eq:PSDconstraint} \\ 
 & \mu_1 = \delta_{x_1^\star}, \ \mu_N = \delta_{x_N^\star},
\end{align}
where the last two equality constraints correspond to the anchor constraints.

\vspace{2mm}
\begin{proposition}\label{thm:shortest}
Assuming that $E$ is a cycle and that there is a unique minimizer $\mathbf{x}^\star$ to \eqref{eq:cost_exp}, this minimizer is uniquely recovered by the SDP \eqref{eq:sdp localization cycle}.
\end{proposition}
\vspace{2mm}
\begin{proof}
Our strategy is to formulate the problem as a weighted shortest path problem. Then we shall recall an equivalent linear programming formulation of this problem, which is in fact equivalent to a modification of the problem~\eqref{eq:sdp localization cycle} in which the semidefinite constraint~\eqref{eq:PSDconstraint} is omitted. Since the recovered solution for this modified problem is then feasible for \eqref{eq:sdp localization cycle}, the result follows.

We now proceed with our reformulation of our sensor network localization problem as a shortest path problem. Recall that we assume without loss of generality that the cycle $E$ consists of the edges $(1,2),\ldots,(N-1,N),(N,1)$ and that $\mathbf{x}^\star_1$ and $\mathbf{x}^\star_N$ are the fixed anchor sensors. We construct a weighted directed graph $\mathcal{G} = (\mathcal{V},\mathcal{E})$ defining our shortest path problems as follows. Let $\mathcal{V}= \mathcal{V}_{1} \sqcup \ldots \sqcup \mathcal{V}_{N}$ where `$\sqcup$' denotes the disjoint union, $\mathcal{V}_{1}=\{\mathbf{x}^\star_1 \}$, $\mathcal{V}_{N}=\{ \mathbf{x}^\star_N \}$, and $\mathcal{V}_i = \mathcal{X}$ for all $i=2,\ldots,N-1$. The `layers' $\mathcal{V}_{2},\ldots,\mathcal{V}_{N-1}$ represent the sets of possible positions for the unknown sensors $\mathbf{x}_2, \ldots, \mathbf{x}_{N-1}$. Next denote the set of all directed edges from the $i$-th layer to the $j$-th layer by $(\mathcal{V}_{i},\mathcal{V}_{j})$, and define the edge set $\mathcal{E}= (\mathcal{V}_{1},\mathcal{V}_{2}) \cup \ldots \cup (\mathcal{V}_{N-1},\mathcal{V}_{N})$. Finally we define a weight function on the edges $f:\mathcal{E}\rightarrow \mathbb{R}$. For each pair of adjacent layers $i,i+1$, the weights on the edges connecting them are defined by the 2-marginal cost matrix $H_{i,i+1}$. A visualization of the shortest path problem for a simple four-sensor example is shown in Fig. \ref{fig:shortest_path}.

\begin{figure}[H]
  \centering
  \includegraphics[width=0.6\textwidth]{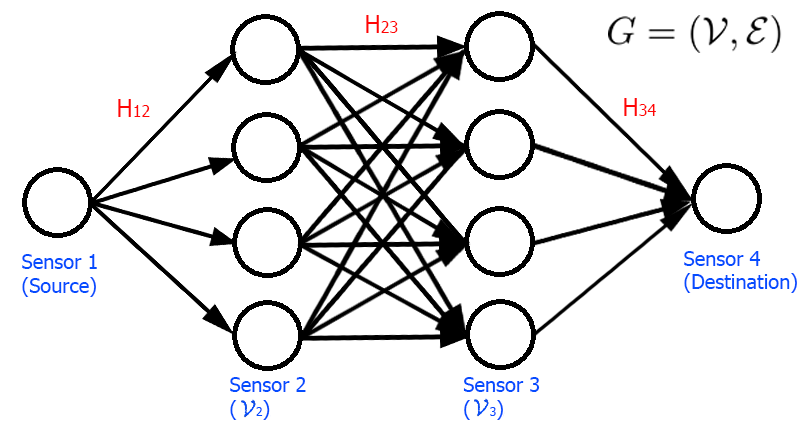}
  \caption{Shortest path problem reformulation for a four-sensor example. In this case, sensor 1 and sensor 4 are fixed anchors. Sensor 2 and sensor 3 have unknown locations. They each have a discrete state space with four positions. Therefore each layer has four nodes. The weights between each pair of layers corresponds to the corresponding blocks of the cost matrix in the 2-marginal relaxation, labeled in red text.}\label{fig:shortest_path}
\end{figure}

Evidently our SNL problem is then equivalent to the shortest path problem on the weighted graph $\mathcal{G}$ with source and destination vertices given by the elements of the first and last layers, respectively. Then we recall the well-known linear programming formulation \cite{shortest_path_lp} for this problem:

\begin{align}
\label{eq:shortest path}
\underset{\{\mu_{i,i+1}\in\mathcal{P}(X_{i}\times X_{i+1})\}_{i=1,\ldots,N-1}}{\mbox{minimize}}\quad & \sum^{N-1}_{i=1} \Tr [ H_{i,i+1}^\top \mu_{i, i+1}] \\
\mathrm{\mbox{subject to }\quad} & 
\mu_{i,i+1}\geq 0, \quad i=1,\ldots,N-1 \nonumber \\
& \mu_{i,i+1}^{\top}\mathbf{1}_{\vert \mc{X}_{i}\vert \times 1} = \mu_{i+1,i+2}\mathbf{1}_{\vert \mc{X}_{i+2}\vert \times 1},\quad i=1,\ldots,N-2\nonumber\\
 & \mu_{1,2}\mathbf{1}_{\vert \mc{X}_{2}\vert \times 1} = \delta_{x_1^\star}, \ \mu_{N-1,N}^{\top}\mathbf{1}_{\vert \mc{X}_{N-1}\vert \times 1} = \delta_{x_N^\star}.\nonumber
\end{align}
The last two constraints identify the source and destination vertices. By our uniqueness assumption the solution of the linear program~\eqref{eq:shortest path} is unique and in particular $\mu^\star_{i,i+1} = \delta_{x_i^\star}\delta_{x_{i+1}^\star}^\top$, where  $(x_1^\star,\ldots,x_N^\star)$ indicates the shortest path. The proof is finished by noticing that problem~\eqref{eq:sdp localization cycle} is equivalent to the shortest path linear program~\eqref{eq:shortest path} if one removes the positive semidefinite constraint in \eqref{eq:PSDconstraint}.
\end{proof}
\vspace{2 mm}

\begin{figure}[!htb]
\centering
\subfloat[SDLSDP solution]{{\includegraphics[trim=0.5cm 0.1cm 0.5cm 0.5cm, clip=true,width=0.4\textwidth]{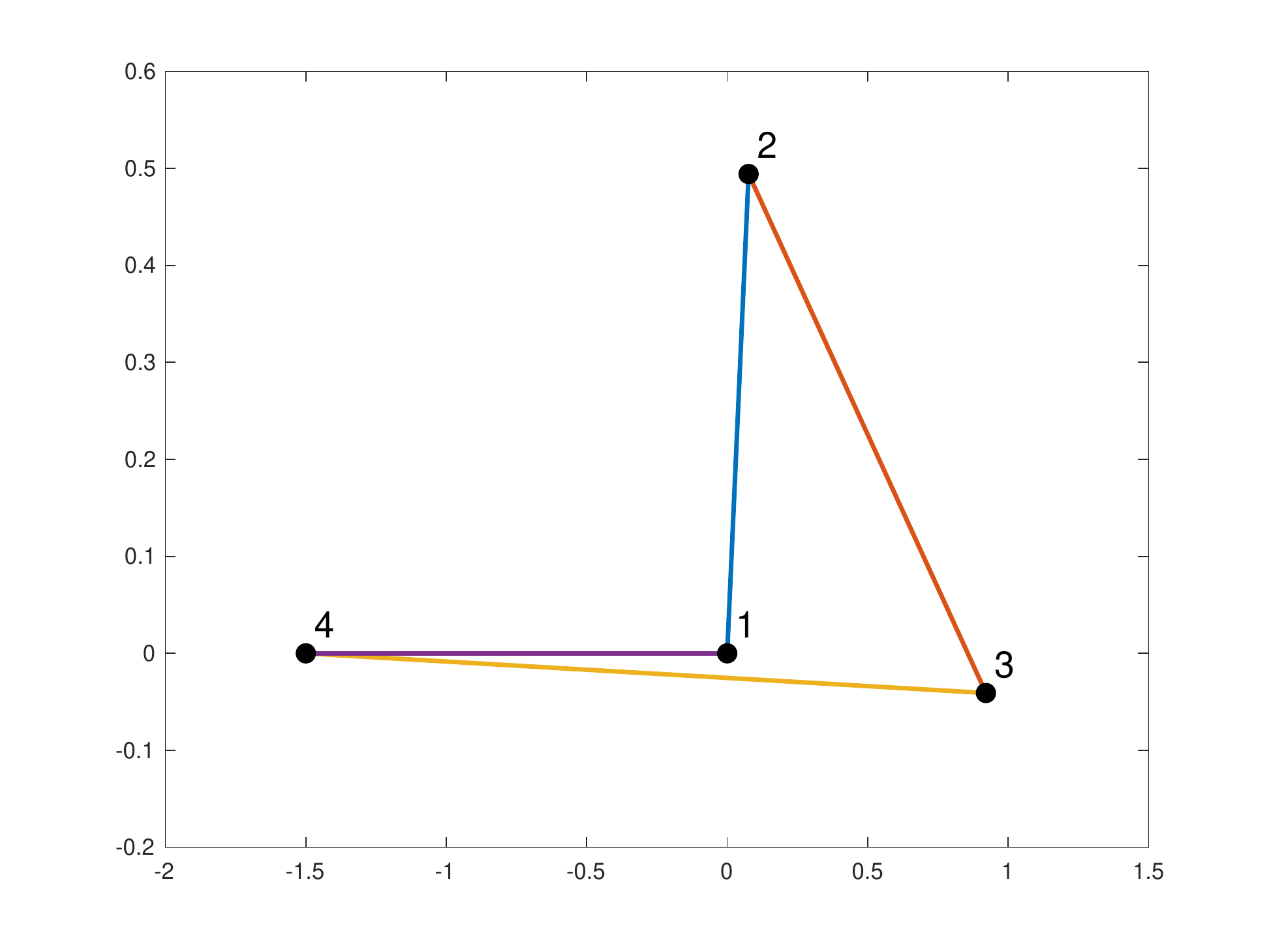}}}
\subfloat[2-marginal relaxation solution]{{\includegraphics[trim=0.5cm 0.1cm 0.5cm 0.5cm, clip=true,width=0.4\textwidth]{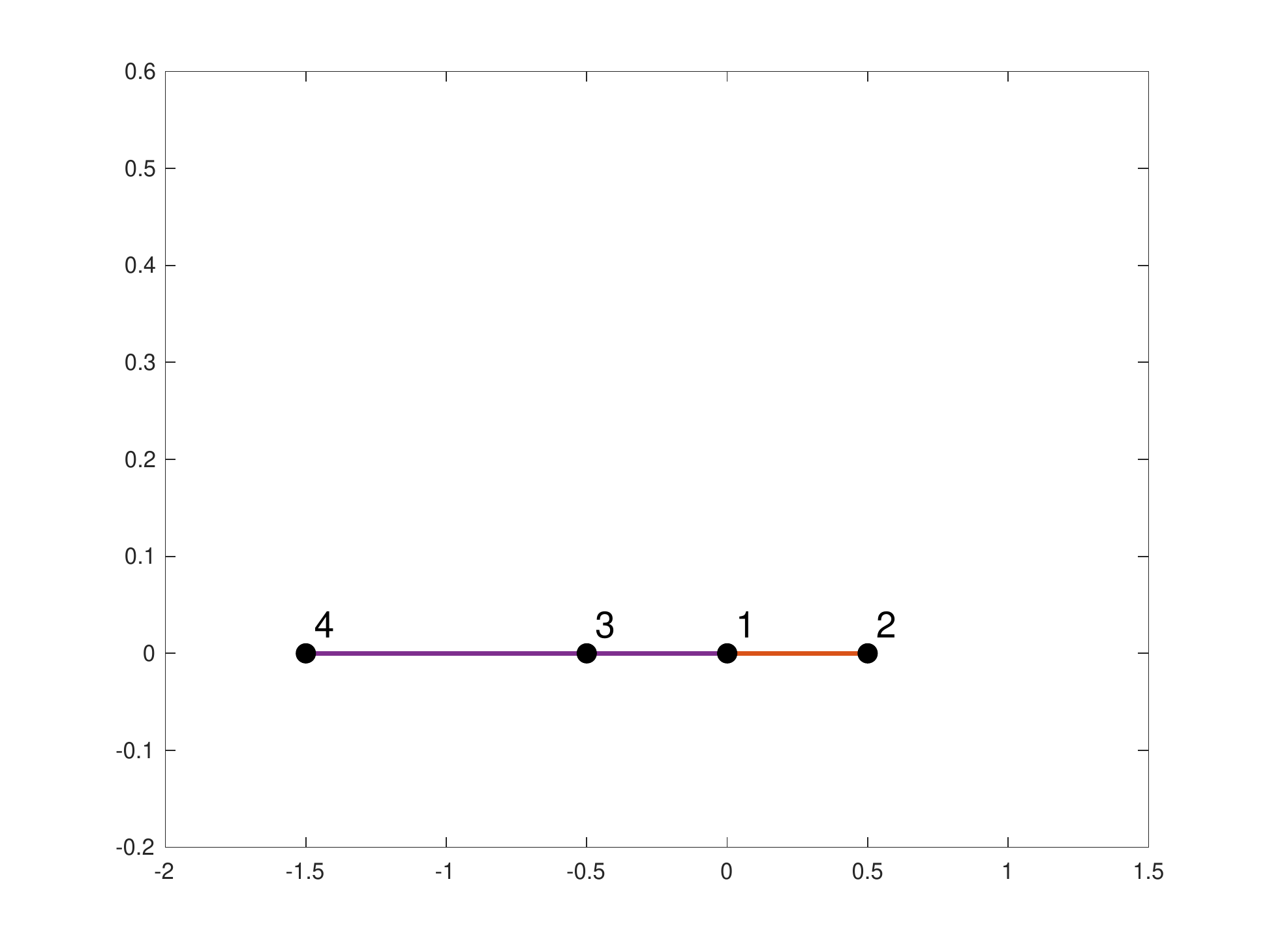}}}
  \caption{While 2-marginal relaxation uniquely recovers the ground truth solution in 1D, SNLSDP recovers the 2D embedding of the true solution}\label{fig:simple_circle_example}
\end{figure}

We have established that in the simple 1D cycle case, the 2-marginal relaxation uniquely recovers the ground truth solution. In turn we demonstrate numerically that SNLSDP may fail to recover this solution. Denote the unknown sensor position matrix by $X=[x_1,\ldots,x_N]$ and its Gram matrix $Y=X^{\top}X$. SNLSDP relaxes the sensor network localization problem to an SDP by relaxing the constraint $Y=X^{\top }X$ to $Y\succcurlyeq X^{\top }X$. Therefore, it is possible for feasible solutions $Y$ to have rank greater than 1, corresponding to recovered sensor positions living in a higher dimensional space. To construct a simple example, let us consider a cycle $E$ consisting of four sensors with $(x^\star_1,x^\star_2,x^\star_3,x^\star_4)=(0,0.5,-0.5,-1.5)$. We let $x^\star_1$ and $x^\star_4$ be the two anchor sensors and solve for $x^\star_2$ and $x^\star_3$. The available distance measurements are $D(1,2)=0.5$, $D(2,3)=1$, $D(3,4)=1$, $D(4,1)=1.5$. We are able to solve \eqref{eq:sdp localization cycle} to recover the ground truth. However, SNLSDP yields the 2D solution (obtained by factoring the recovered Gram matrix $Y$) $x_2=(0.0762,0.4942)$ and $x_3=(0.9210,-0.0409)$. Although the 2D system exactly matches the distance measurements and yields the objective function value of zero, SNLSDP fails to recover the 1D solution. Fig. \ref{fig:simple_circle_example} visualizes the solutions obtained via both SNLSDP and 2-marginal relaxation \eqref{eq:sdp localization cycle} in 2D space.

\section{A multiscale algorithm}\label{sec:MMR}

In this section, we detail a multiscale algorithm for solving the marginal relaxation \eqref{eq:sdp} in practice for continuous state spaces. We call this approach the multiscale marginal relaxation (MMR). Crucial to the success of the algorithm is the smoothness of the $H_{ij}$. For simplicity of presentation, we assume all particles live in the same space, i.e., $\mathcal{X}_1= \mathcal{X}_2=\cdots = \mathcal{X}_N = \mathcal{X}$. After describing the MMR algorithm, we will then propose an extra post-processing procedure to solve the problem of sampling near-optimal solutions when the output of MMR is ambiguous.

\subsection{MMR}

The main observation we have taken from the previous section is that since we expect \eqref{eq:sdp} to be a surrogate for \eqref{eq:lp}, we in turn expect the $\mu_i$ to be delta measures where $\supp{(\mu_i)}$ approximates the position of $x_i^\star$ in \eqref{eq:sdp}. If we have some understanding that $x_i^\star$ is unlikely to occur in some locations in $\mathcal{X}$, we can conveniently force $\mu_i$ to be zero in those locations. From the form of the positive semidefinite matrix $G$ in \eqref{eq:sdp}, where
\begin{equation}
\text{diag}(G) = \text{diag}([\mu_i]_{i=1}^N),
\end{equation}
zeroing entries of $\mu_i$'s effectively reduces the size of $G$. Therefore, a rough idea of where $x_i^\star$ could be allows us to access problems with larger $N$ and $M$. This motivates us to first solve a coarsened version of problem \eqref{eq:sdp}, where each $H_{ij}$ is coarsened according some coarse partition of $\mathcal{X}$. The support of $\mu_i$ on this coarse partition should provide a rough guess for the position of particle $i$. Now in a multiscale way we refine this guess. By identifying the supports of the $\mu_i$ for the coarse problem, we apply a finer partitioning \emph{only} within these supports and then solve a problem using our finer partition (of a restricted domain). 

The proposed algorithm is summarized in Alg.~\ref{alg:MMR}. It consists of three main subroutines: coarsen, propagate, and refine, where the naming of these subroutines follows from \cite{gerber2017multiscale}. In our multiscale framework, we use $k=1$ to denote the coarsest level, which is also the first level on which we plan to solve an SDP. As $k$ increases, the partitions become finer. The algorithm works with the following objects:
\vspace{2mm}
\begin{enumerate}
    \item Quadrature points $X = \{q_1,\ldots,q_M\}\subset \mathbb{R}^{d}$. We use the quadrature points to discretize the state space $\mathcal{X}$, and we denote their associated quadrature weights $w = [w_{q_1},\ldots,w_{q_M}]$.
    \item Partitions $\{S^{(k)}_l\}_{l=1}^{M^{(k)}}$, $k=1,\ldots,K$. The partitions satisfy the following properties:
    \begin{itemize}
    \item Disjoint partition at each level $k$: $S^{(k)}_l\subset X$ for all $l$,  $\cup_{l=1}^{M^{(k)}} S^{(k)}_l = X$, and $S^{(k)}_{l_1} \cap S^{(k)}_{l_2} = \emptyset$ if $l_1\neq l_2$. 
    \item Nestedness: for any $m=1,\ldots, M^{(k+1)}$, the $m$-th subset $S^{(k+1)}_{m}$ at the $(k+1)$-th level is properly contained in one and only one $S^{(k)}_{l}$, where $l \in \{1,\ldots,M^{(k)}\}$. If $S^{(k+1)}_{m} \subset S^{(k)}_{l}$, we call $S^{(k+1)}_{m}$ the \emph{child} of $S^{(k)}_{l}$ (and $S^{(k)}_{l}$ the \emph{parent} of  $S^{(k+1)}_{m}$).
    \item At the finest level $k=K$, $M^{(K)} = M$, and the partition is given by
    \[
    S^{(K)}_l = \{\{q_l\} \,:\, l=1,\ldots,M \}.
    \]
    \end{itemize}
    \item Selected parts $\{p_i(k)\}_{i=1}^N$, $k=1,\ldots,K$. In the multiscale strategy, for each particle we only discretize the cost $H_{ij}$ around some locations where $\mu_i$ is non-zero.
    The choice of $p_i(k) \subset \{1,\ldots, M^{(k)} \}$ indicates that the $i$-th particle is believed to be on one of the grid points in $\cup_{l\in p_i(k)} S^{(k)}_l$. For convenience we denote $S^{(k)}_{p_i(k)}=\cup_{l\in p_i(k)} S^{(k)}_l$.
\end{enumerate}
\vspace{2mm} 

Besides the basic framework of Alg.~\ref{alg:MMR}, there are several parameters involved in processing the solutions of the inner 2-marginal relaxations, which are of the form \eqref{eq:sdp}. Here we briefly introduce these parameters and their functionalities.
\vspace{2mm}
\begin{enumerate}
    \item Thresholds $\{\eta^{(k)}\}_{k=1}^{K}$: For each particle $i$, the multiscale strategy relies on determining the support of the 1-marginal $\mu_i$ in order to identify parts $p_{i}(k)$ in which the true solution may be located. In practice we find the support by thresholding, i.e., at level $k$ the entries of $\mu_i$ are identified to be 0 if they are less than $\eta^{(k)}$.
    \item Upper bounds $\{u^{(k)}\}_{k=1}^{K}$: In our multiscale algorithm we generally expect the discrete 1-marginals recovered from our inner 2-marginal relaxations of the form~\eqref{eq:sdp} to be exactly or nearly delta measures. However, this useful property can have the following adverse effect. At the coarse levels of the multiscale algorithm, discretization error may yield discrete coarse problems for which the optimal solution does not correspond to the optimal solution of the underlying continuous problem, particularly if this underlying problem has many near-optimal local minima. In this case, the support of the solution of our coarse SDPs may fail to include the true optimizer. Hence we want to force $\mu_{i}$ to be less concentrated in order to `hedge our bets' until the best solution can be determined after sufficient refinement. Numerically this is implemented by including an additional entrywise upper bound constraint on the 2-marginal $\{\mu_{ij}\}_{i,j=1}^{N}$ in \eqref{eq:sdp}. Since $\mu_{ij}\in[0,1]$ already, our upper bound $u^{(k)}$ should lie in $(0,1]$, and we enforce the constraint $0\leq \mu_{ij}\leq u^{(k)}$.
\end{enumerate}
\vspace{2mm}

\begin{algorithm}[!htb]
\caption{Pseudocode for Multiscale Marginal Relaxation (MMR)} \label{alg:MMR}
\begin{algorithmic}[1]
\Procedure{MMR}{ $\{H_{ij}\}_{i,j=1}^N, \ w, \  \{ \{S^{(k)}_l\}_{l=1}^{M^{(k)}} \}_{k=1}^K $ }   
\For{$k=1$ \text{to} $K$}
    \If{$k>1$}
        \State $p_i(k)=\{l:S^{(k)}_{l}\subset S^{(k-1)}_{p_i(k-1)}\}$, $i=1,\ldots,N$
    \EndIf
    \State $\{H^{(k)}_{ij}\}_{i,j=1}^N$ $\gets$ Coarsen($\{H_{ij}\}_{i,j=1}^N, w, \{S^{(k)}_i\}_{i=1}^{M^{(k)}}, \{p_i(k)\}_{i=1}^{N}$)
    \State $\{\tilde p_i(k)\}_{i=1}^{N}$ $\gets$
    Propagate($\{H^{(k)}_{ij}\}_{i,j=1}^N$, $w$)
    \State $\{p_i(k)\}_{i=1}^{N}$ $\gets$ Refine($\{H_{ij}\}_{i,j=1}^N, \{\tilde p_i(k)\}_{i=1}^{N}$)
\EndFor
\EndProcedure
\end{algorithmic}
\end{algorithm}

In Algorithm~\ref{alg:MMR}, at the beginning of each level $k$, we first initialize the selected parts $\{p_i(k)\}_{i=1}^N$. Specifically, at the starting level $k=1$ the initial guess is the entire discretized state space, i.e., $p_i^{(1)}=\{1,2,\ldots,M^{(1)}\}$ for $i=1,2,\ldots,N$. For $k>1$, we initialize $\{p_i(k)\}_{i=1}^N$ based on the selected parts from the last level $\{p_i(k-1)\}_{i=1}^N$. Assume that we are given $p_i(k-1)$ from the last level, indicating that at level $k-1$, the $i$-th particle is near grid points contained in
$ S^{(k-1)}_{p_i(k-1)}$.
Using this information at coarse level, we initialize $p_i(k)$ to include all children of the part at level $k-1$, i.e.,
\begin{equation}
 p_i(k)=\{l:S^{(k)}_{l}\subset S^{(k-1)}_{p_i(k-1)}\},\  i=1,\ldots,N.
\end{equation}

We now describe in detail the `Coarsen,' `Propagate,' and `Refine' steps of Algorithm~\ref{alg:MMR} which complete the pseudocode.

\subsubsection{Coarsen}
At each level $k$, we perform a coarsening of the objective function $H_{ij}$ for each $i,j=1,\ldots,N$ in order to obtain an SDP of manageable size. Although we could discretize $H_{ij}$ into a size $M^{(k)}\times M^{(k)}$ matrix, this is not helpful as at the finest scale, $M^{(K)} = M$, and $M$ can be prohibitively large. However, the initial guesses $\{p_i(k)\}_{i=1}^{N}$ provide us a way to restrict the problem to the selected parts for each particle and ultimately reduce the problem size. Indeed, we can restrict $H_{ij}$ to the product set $S^{(k)}_{p_i(k)}\times S^{(k)}_{p_j(k)}$:
\begin{equation}
\tilde H_{ij}^{(k)} := H_{ij}\left(S^{(k)}_{p_i(k)}, S^{(k)}_{p_j(k)}\right).\label{eq:cost1},
\end{equation}
where we interpret the right-hand side as a matrix of size
\begin{equation}
\sum_{S\in S^{(k)}_{p_i(k)}} \vert S \vert \times \sum_{S'\in S^{(k)}_{p_j(k)}} \vert S' \vert.
\end{equation}
Since $\tilde H_{ij}^{(k)}$ has values defined on the finest grid points, we coarsen it to form a matrix
\begin{equation}
H^{(k)}_{ij} \in \mathbb{R}^{\left\vert p_i(k) \right \vert\times \left \vert p_j(k) \right \vert},
\end{equation}
by taking averages within each part according to our quadrature weights, i.e., we form:
\begin{multline}
H^{(k)}_{ij} = \cr \left[\begin{smallmatrix}  w\left(S^{(k)}_{p_i(k)}(1)\right)^\top & & \\ & \ddots &  \\ & &  w\left(S^{(k)}_{p_i(k)}(\vert  p_i(k) \vert)\right)^\top \end{smallmatrix}\right] \tilde H^{(k)}_{ij} \left[\begin{smallmatrix}  w\left(S^{(k)}_{p_j(k)}(1)\right) & & \\ & \ddots &  \\ & &  w\left(S^{(k)}_{p_j(k)}(\vert  p_j(k)\vert)\right) \end{smallmatrix}\right].\label{eq:cost2}
\end{multline}

\subsubsection{Propagate}
At this point we have already obtained the coarsened cost at level $k$, restricted to selected parts in $k$-level partition. More precisely, each particle $i$ now is restricted to
$m^{(k)}_i := \left\vert p_i(k) \right\vert $
In order to determine which parts of the level $k$-partition the particles live in, in principle we want to solve the coarsened version of problem \eqref{eq:lp}:
\begin{align}
\underset{\{\mu^{(k)}_{ij}\}, \mu^{(k)}}{\mbox{minimize}}\quad\quad\quad & \sum_{i<j} \Tr [ {H^{(k)}_{ij}}^\top \mu^{(k)}_{ij}] \\
\mathrm{\mbox{subject to }\quad\quad\quad} & \mu^{(k)}_{ij}\ \text{is 2-marginal of}\ \mu^{(k)},\ i\neq j\in [N].  \nonumber \\ 
 & \mu^{(k)}\in \mathcal{P}\left( [m^{(k)}_1] \times  \cdots\times  [m^{(k)}_N] \right)  \nonumber,
\end{align}
where we denote $[m]:=\{1,\ldots ,m \}$.
As mentioned previously, this problem is impossible to solve for $N$ of even moderate size, so we solve its surrogate 2-marginal relaxation \eqref{eq:sdp}. Also, assuming $u^{(k)}< 1$, we inclue an additional upper bound constraint on the 2-marginal variables, yielding the following SDP:
\begin{align}
\label{eq:sdp_red}
\underset{\substack{\{\mu^{(k)}_{ij}\in\mathcal{P}([m^{(k)}_i]\times  [m^{(k)}_j]\}_{i<j},\\ \{\mu^{(k)}_{i}\in\mathcal{P}([m^{(k)}_i])\}_i}}{\mbox{minimize}}\quad\quad\quad & \sum_{i<j} \Tr [ {H^{(k)}_{ij}}^\top \mu^{(k)}_{ij}] \\
\mathrm{\mbox{subject to }\quad\quad\quad} & \mu^{(k)}_{ij}\mathbf{1}_{m^{(k)}_j}=\mu^{(k)}_{i},\ \ {\mu^{(k)}_{ij}}^{\top}\mathbf{1}_{m^{(k)}_i}=\mu^{(k)}_{j},\quad i\neq j \nonumber \\
 & G(\{\mu^{(k)}_{ij}\})\succeq0,\ \ \mu_{ij}^{(k)}\leq u^{(k)}.\nonumber
\end{align}

Finally, we solve SDP \eqref{eq:sdp_red} and define the `intermediate parts' that each particle lives in at level $k$ (i.e., the $\tilde p_i(k),\ i=1,\ldots,N$) as the thresholded supports of the $\mu^{(k)}_i$, which are obtained from the solution of \eqref{eq:sdp_red}. Specifically,
\[
\tilde{p}_i(k) := 
\{ l \in [m_i^{(k)}] : \mu_{i}^{(k)}(l) \geq \eta^{(k)} \}.
\]
The interpretation is that for each $i=1,\ldots,N$, it is likely that  $x_i^\star$ is near one of the grid points in $S^{(k)}_{\tilde p_i(k)}$. 

\subsubsection{Refine}
At this point, we have solved the SDP problem on our level-$k$ grid and replaced the initial part $p_i(k)$ for $i$-th particle with the intermediate part $\tilde p_i(k)$.
Due to discretization error, it is often the case that the true location is in fact covered by grid points that neighbor those of $S^{(k)}_{\tilde p_i(k)}$.
We take these points into consideration to make sure that we do not miss the true solution. More precisely, we will add them to the collection of points in our intermediate parts and solve another SDP on this expanded collection, yielding final estimates for the parts via another thresholding step. We denote the neighborhood of grid points of $S^{(k)}_{\tilde p_i(k)}$ as
\begin{equation}
\mathcal{N}^{(k)}_i =  S^{(k)}_{\tilde p_i(k)} \cup \bigcup_{l\in[M^{(k)}]\,:\, l \sim \tilde p_i(k)} S_l^{(k)},
\end{equation}
where we  use `$\sim$' to denote adjacency in the level-$k$ grid.
We now rediscretize $H_{ij}$ on this neighborhood, letting
\begin{equation}
\tilde R^{(k)}_{ij} = H_{ij}(\mathcal{N}^{(k)}_i,\mathcal{N}^{(k)}_j)
\end{equation}
and taking the quadrature weightings:
\begin{multline}
R^{(k)}_{ij} = \cr \left[\begin{smallmatrix}  w\left(\mathcal{N}^{(k)}_i(1)\right)^\top & & \\ & \ddots &  \\ & &  w\left(\mathcal{N}^{(k)}_i(\vert \mathcal{N}^{(k)}_i \vert)\right)^\top \end{smallmatrix}\right] \tilde R^{(k)}_{ij} \left[\begin{smallmatrix}  w\left(\mathcal{N}^{(k)}_{j}(1)\right) & & \\ & \ddots &  \\ & &  w\left(\mathcal{N}^{(k)}_{j}(\vert \mathcal{N}^{(k)}_{j} \vert)\right) \end{smallmatrix}\right].
\end{multline}

In order to determine the part within $\mathcal{N}^{(k)}_i$ where each $i$-th particle lives, we again solve an SDP:
\begin{align}
\label{eq:sdp_red2}
\underset{\substack{\{\mu^{(k)}_{ij}\in\mathcal{P}([\vert\mathcal{N}^{(k)}_i\vert] \times [\vert \mathcal{N}^{(k)}_j\vert ]\}_{i<j},\\ \{\mu^{(k)}_{i}\in\mathcal{P}([\vert\mathcal{N}^{(k)}_{i}]\vert\}}}{\mbox{minimize}}\quad\quad\quad & \sum_{i<j} \Tr [ {R^{(k)}_{ij}}^\top \mu^{(k)}_{ij}] \\
\mathrm{\mbox{subject to }\quad\quad\quad} & \mu^{(k)}_{ij}\mathbf{1}_{\vert \mathcal{N}^{(k)}_{j} \vert}=\mu^{(k)}_{i},\ \ {\mu^{(k)}_{ij}}^{\top}\mathbf{1}_{\vert\mathcal{N}^{(k)}_{i}\vert}=\mu^{(k)}_{j},\quad i\neq j \nonumber \\
 & G(\{\mu^{(k)}_{ij}\})\succeq0,\ \ \mu_{ij}^{(k)}\leq u^{(k)}.\nonumber
\end{align}

Just as the propagating step, the final refined parts $\{p_i(k)\}_{i=1}^N$ of the particles at level $k$ are obtained by determining the thresholded supports of the $\mu^{(k)}_i$ obtained by solving the problem \eqref{eq:sdp_red2}. This means that for each $i$, our new guess is that $x_i^\star$ is contained within the grid points $S^{(k)}_{p_i(k)}$. This refinement procedure can in fact be iterated several times to improve the guess for the refined parts $\{p_i(k)\}_{i=1}^N$.

\subsection{Exploring near optimal solutions}\label{sec:multiple soln}
At the finest level $K$ the discretization error is low, so ideally we expect the $\mu^{(K)}_i$ to be delta-measures, with support corresponding to the support of the true $\mu_i$. If this is the case, $S_{p_i(K)}^{(K)}$ is a singleton, and we can simply take its unique element $x_i^{\star}$ as the final MMR output. (Recall that $S_{l}^{(K)} = \{q_l\}$ for all $l$.) However, in some practical problems there exist multiple globally-optimal solutions or nearly-optimal solutions. In this case, the recovered $\mu^{(K)}_i$ often have non-singleton support. A direct strategy for recovering an unambiguous solution is to pick the largest entry of $\mu^{(K)}_i$. This strategy generally succeeds in recovering a near-optimal solution but does not allow us to explore the space of near-optimal configurations. We propose a post-processing step to tackle this problem. 

Consider the case where each $\mu^{(K)}_i$ may consist of multiple atoms, i.e., possibly $m^{(K)}_i = \vert p_i(K) \vert>1$. We first obtain a cost matrix $H_{ij}^{(K)}$ by the restriction as $H_{ij}( S_{p_i(K)}^{(K)} , S_{p_j(K)}^{(K)} )$. However, we random noise to the cost matrix in order to generate additional biases, i.e., we define $\hat{H}_{ij}^{(K)}=H_{ij}^{(K)}+ \lambda R_{ij}$ where $R_{ij}$ is a random matrix with independent standard normal entries. We rely on this noise to sample from multiple near-optimal solutions. Note that the magnitude $\lambda$ of the random noise should be comparable to the range of costs among  the near-optimal solutions. Then we solve the 2-marginal relaxation:
\begin{align}\label{eq:sdp_rand}
\underset{\substack{\{\tilde \mu^{(K)}_{ij}\in\mathcal{P}([m^{(K)}_i]\times [m^{(K)}_j]\}_{i<j},\\ \{\tilde \mu^{(K)}_{i}\in\mathcal{P}([m^{(K)}_i])\}}}{\mbox{minimize}}\quad\quad\quad & \sum_{i<j} \Tr [ \hat{H}_{ij}^{(K)\top} \tilde \mu^{(K)}_{ij}] \\
\mathrm{\mbox{subject to }\quad\quad\quad} & \tilde \mu^{(K)}_{ij}\mathbf{1}_{m^{(K)}_j}=\tilde \mu^{(K)}_{i},\ \ {\tilde \mu^{(K)\top}_{ij}}\mathbf{1}_{m^{(K)}_i}=\tilde \mu^{(K)}_{j},\quad i\neq j \nonumber \\
 & G(\{\tilde \mu^{(K)}_{ij}\})\succeq0.\nonumber
\end{align}

Here we do not impose upper bounds on the 2-marginals since we want to recover a unique near-optimal solution. If the 1-marginals $\tilde \mu_i^{(K)}$ that we recover are delta-measures, then it is trivial to extract a solution. Otherwise, we compute the top eigenvector $v_1$ of the 2-marginal matrix $G(\{\tilde \mu^{(K)}_{ij}\})$. Note that this eigenvector has nonnegative entries by the Perron-Frobenius theorem. Assuming that the relaxation is reasonably tight, the top eigenvector should approximately satisfy
\begin{equation}
v_1 \approx [{\delta_{\tilde{x}_1^{\star}}^\top } ,\cdots, {\delta_{\tilde{x}_N^{\star}}^\top}]^\top,
\end{equation}
where the $\tilde{x}_i^\star$ correspond to optimal solution of the appropriate combinatorial optimization problem induced by $\hat{H}_{ij}$.
Then for each particle we pick the largest entry in the corresponding block of $v_1$ as our position for this particle. By solving different random instances of this problem we explore various near-optimal solutions.

\section{Numerical experiments}\label{sec:numerical}
Below we present numerical experiments for the MMR method. In Section \ref{sec:41} we apply the method to a sensor network localization (SNL) problem, and in Section \ref{sec:42} we apply the method to the exploration of the many near-optimal configurations of the Lennard-Jones (LJ) potential.

\subsection{Sensor network localization}\label{sec:41}
Here we consider $N=50$ particles in two spatial dimensions. Each particle has a state space restricted to the square region $\mathcal{X}=[0,10]^{2} \subset \mathbb{R}^2$. We sample $50$ points uniformly from the square as our ground truth for the particle positions and let $D_{0}$ denote the corresponding pairwise distance matrix. Then we let 
\[
D(i,j) = D_0 (i,j) + b_{ij} z_{ij},
\]
for $i\neq j$, where $b_{ij} \in \text{Bernoulli}(\sigma)$ and $z_{ij} \sim \text{Unif}[0,3]$ are all independently distributed. Here $\sigma \in (0,1)$ is a parameter that we will vary, corresponding to the expected proportion of entries of the distance matrix to be contaminated by noise. We view $D$ as the (corrupted) measurements specifying our SNL problem. We also introduce a `sensing radius' $D_{\max}$. If the distance between two particles is greater than this radius, then we do not have any observation (even corrupted) about their pairwise distance.
This problem specification is reflected in the cost function
\begin{equation}
\label{eq:SNLpairObj}
    H_{ij}(x_{i},x_{j})=\left\{
    \begin{aligned}
    &\sqrt{(||x_{i}-x_{j}||_{2}-D(i,j))},\quad &&i\neq j,\  D_0 (i,j) \leq D_{\max}, \\
    & 0, &&i=j.
    \end{aligned}
    \right.
\end{equation}

Now we describe the choice of method parameters for our numerical experiments. For our discretization points $q_l$ we consider a $2^{7}\times2^{7}$ regular grid, corresponding to the finest level $K=6$ of our multiscale discretization. In our experiments we always simply choose uniform quadrature weights. We coarsen the partition by successively merging $2\times 2$ blocks into one. Therefore $M^{(k)}=4^{(k+1)}$, and specifically $M^{(1)}=16$. We let $u^{(k)}=1$ for $k=1,\ldots,6$ (so the upper bound constraint is effectively inactive in \eqref{eq:sdp_red}). The thresholding parameters for the 1-marginals $\mu_i^{(k)}$ are chosen to be $\eta^{(k)}=5\times 10^{-2}$ for $k=1,\ldots,6$. In addition if there are fewer than $3$ entries larger than the threshold at each level, we force the algorithm to pick $3$ largest entries from the 1-marginal $\mu_i$. In the refinement step, we consider grid points to be `adjacent' if one is in the Moore neighborhood of the other. The Moore neighborhood of a grid point consists of the eight grid points surrounding this grid point (except on the boundaries), so the inequality $|\mathcal{N}_{i}^{(k)}|\leq 9\vert \tilde p_i(k)\vert$ controls the size of each optimization problem in the refinement step. We iterate the refinement step 3 times, or until there is no change to $\tilde p_i(k)$.

In order to fix the degeneracy of the SNL problem with respect to rigid motions, we fix the positions of the first 3 particles (which we call anchor particles) to their ground truth values. It is then equivalent either to view our optimization problem as an optimization problem over $N-3$ particles in the square with a suitably modified pairwise objective, or simply to constrain the discretization grids for the first three particles to be the appropriate singleton sets.

We solve the SDPs in MMR using the SDPNAL+ package~\cite{sdpnal+}. The SDP subproblems are solved to moderate accuracy (tolerance $=10^{-3}$). We will compare the results of our method with the results obtained both by directly minimizing the nonlinear cost function and by simulated annealing  (SA) via the MATLAB functions \texttt{fmincon} and \texttt{simulannealbnd}, respectively, initialized with particle locations chosen uniformly from the square $[0,10]^2$. Additionally, we will compare our approach with SNLSDP~\cite{snlsdp}, a semidefinite relaxation approach for sensor network localization.

\subsubsection{Comparison with other methods} We compare our approach with \texttt{fmincon}, SA, and SNLSDP. We also use the output at the finest level of MMR as the initial guess for \texttt{fmincon} in order to mitigate the effect of discretization error. Since SNLSDP does not directly optimize~\eqref{eq:SNLpairObj}, we also can view it as an initial guess for \texttt{fmincon} as we do below in our tests.

For each $\sigma, D_{\max}$ considered, we test on 100 independent instances of the problem. We let $\epsilon_{p}$  denote the 2-norm error of the recovered particle positions (averaged over the particles), and we let $\epsilon_{e}$ denote the value of the recovered cost or `energy' (i.e., the cost function~\eqref{eq:SNLpairObj} evaluated at the recovered particle positions). In this subsection, we fix $D_{\max}=6$, causing roughly half of distance measurements to be omitted. For $\sigma$, we consider the two cases $\sigma=0.1$ and $\sigma=0.2$. The results are summarized in Tables \ref{table:1} and \ref{table:ext1}.

\begin{table}[h]
\centering
 \caption{Mean position error $\epsilon_{p}$ and mean energy $\epsilon_{e}$, together with corresponding standard deviations. The sensing radius is $D_{\max} =6$, and the noisy proportion is $\sigma=0.1$. Results are obtained from 100 independent realizations.}\label{table:1}
\begin{tabular}{|c|c|c|c|c|c|c|}
\hline
 & \multirow{2}{*}{\texttt{fmincon}} & \multirow{2}{*}{SA} & \multirow{2}{*}{SNLSDP} & SNLSDP & \multirow{2}{*}{MMR} & MMR\\
   & & & & +\texttt{fmincon} & & +\texttt{fmincon}\\
\hline
\multirow{2}{*}{$\epsilon_{p}$}  & 3.6152 & 3.8465 & 0.5940 & 0.2585 & 0.0913 & 0.0439\\
 & $\pm$1.1835 & $\pm$0.5120 & $\pm$0.9724 & $\pm$0.8321 & $\pm$0.4486 & $\pm$0.4506\\
\hline
\multirow{2}{*}{$\epsilon_{e}$}  & 1058.70 & 1545.60 & 773.74 & 359.01 & 480.34 & 327.6 \\
 & $\pm$211.13 & $\pm$109.62 & $\pm$73.90 & $\pm$98.68 & $\pm$40.71 & $\pm$32.15\\
\hline
\end{tabular}
\end{table}
\begin{table}[h]
\centering
 \caption{Mean position error $\epsilon_{p}$ and mean energy $\epsilon_{e}$, together with corresponding standard deviations. The sensing radius is $D_{\max} =6$, and the noisy proportion is $\sigma=0.2$. Results are obtained from 100 independent realizations.}\label{table:ext1}
\begin{tabular}{|c|c|c|c|c|c|c|}
\hline
 & \multirow{2}{*}{\texttt{fmincon}} & \multirow{2}{*}{SA} & \multirow{2}{*}{SNLSDP} & SNLSDP & \multirow{2}{*}{MMR} & MMR\\
   & & & & +\texttt{fmincon} & & +\texttt{fmincon}\\
\hline
\multirow{2}{*}{$\epsilon_{p}$}  & 3.7020 & 3.8914 & 1.0030 & 0.4456 & 0.1388 & 0.0926\\
 & $\pm$0.8574 & $\pm$0.5358 & $\pm$1.1545 & $\pm$1.0540 & $\pm$0.4941 & $\pm$0.5030\\
\hline
\multirow{2}{*}{$\epsilon_{e}$}  & 1000.01 & 1227.14 & 772.08 & 355.36 & 436.18 & 313.03 \\
 & $\pm$121.53 & $\pm$87.68 & $\pm$78.63 & $\pm$95.54 & $\pm$40.77 & $\pm$34.95\\
\hline
\end{tabular}
\end{table}

\par
\begin{figure}[!htb]
\centering
\subfloat[$D_{\max}=6$, $\sigma=0.1$]{{\includegraphics[trim=0.5cm 0.1cm 0.5cm 0.5cm, clip=true,width=0.5\textwidth]{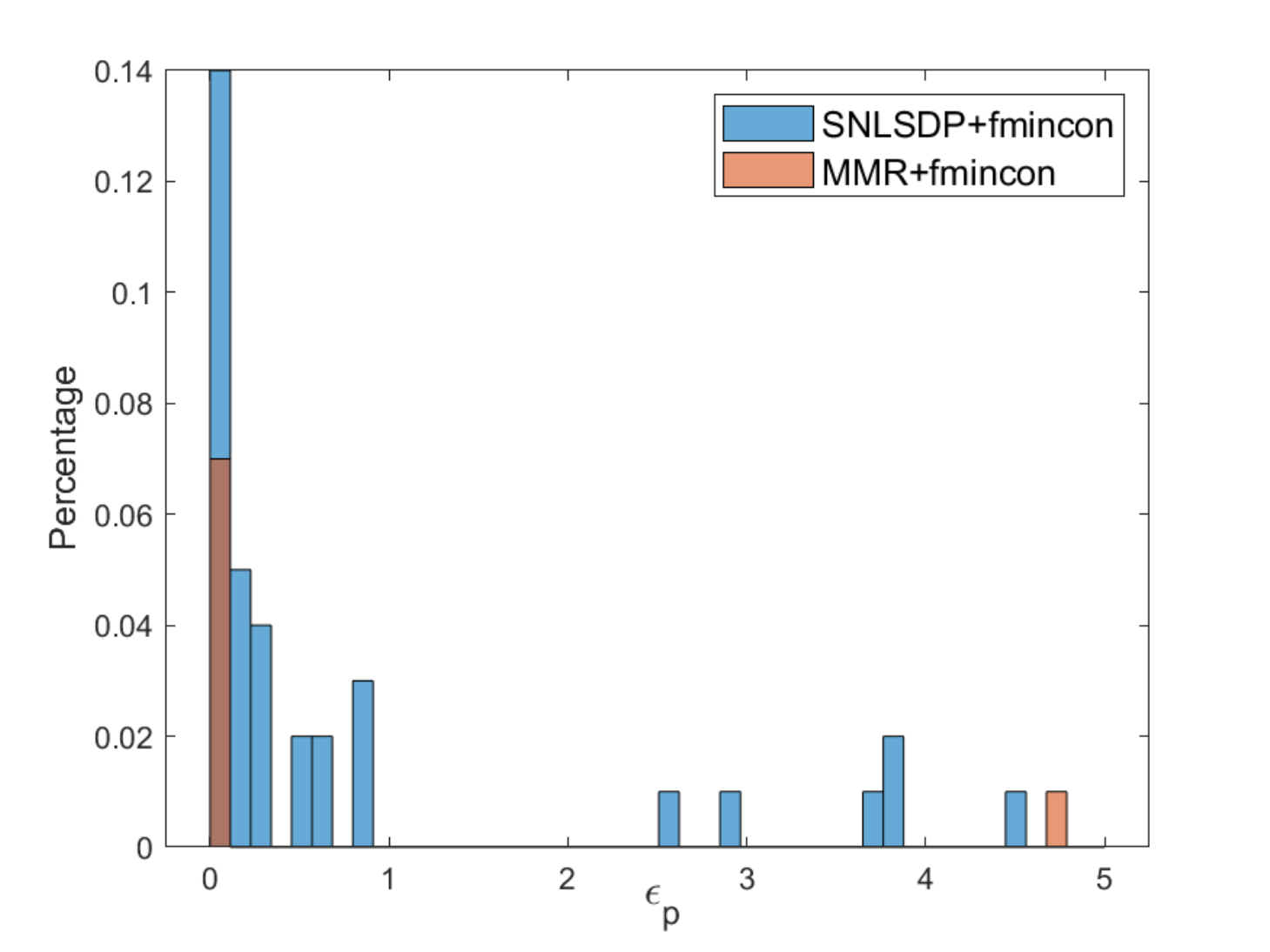}}}
\subfloat[$D_{\max}=6$, $\sigma=0.2$]{{\includegraphics[trim=0.5cm 0.1cm 0.5cm 0.5cm, clip=true,width=0.5\textwidth]{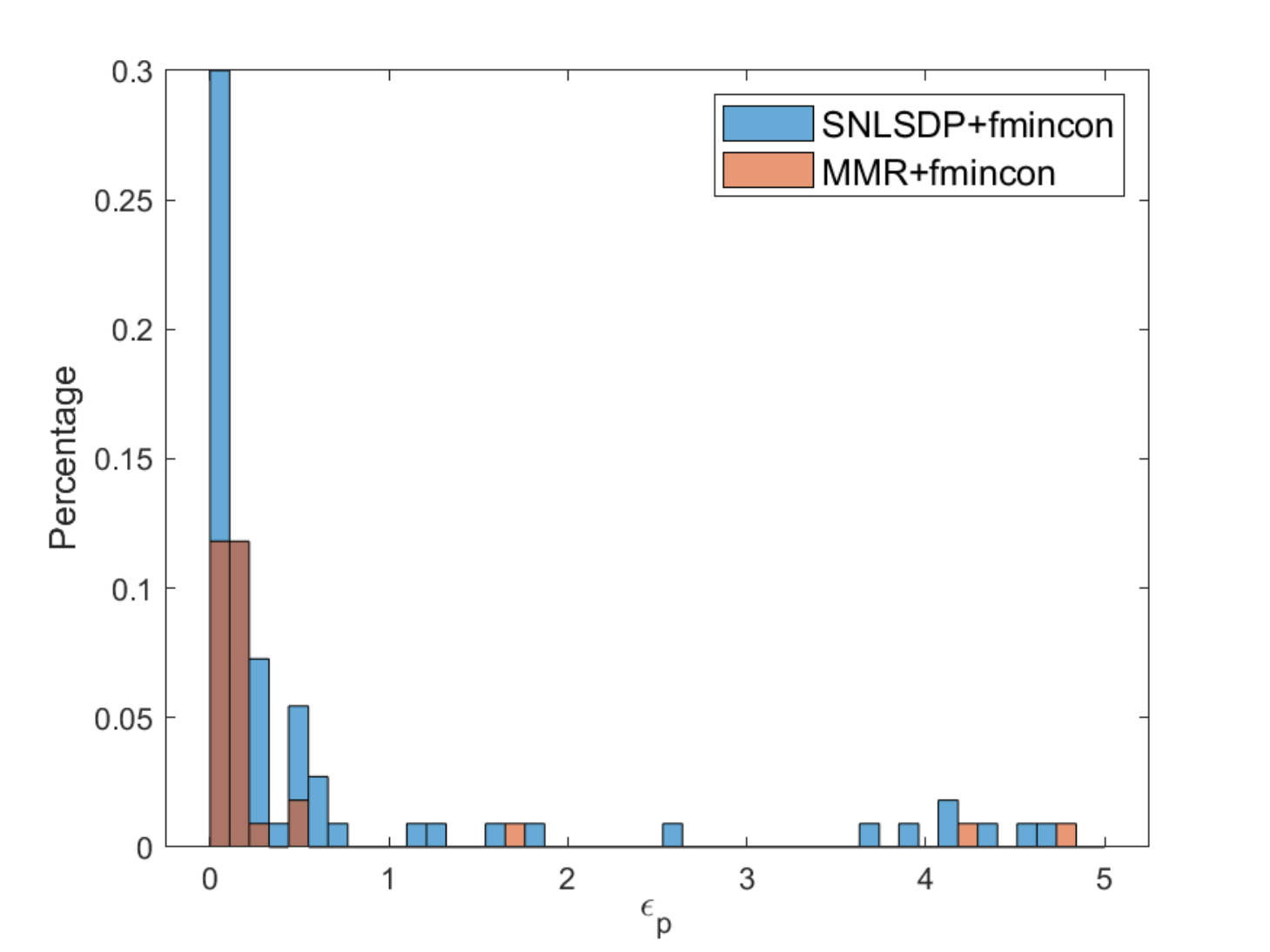}}}
  \caption{$\epsilon_p$ histograms of MMR$+\texttt{fmincon}$ and SNLSDP$+\texttt{fmincon}$ in the `failed' cases (i.e., cases where $\epsilon_p \geq 10^{-5}$). (a) $D_{\max}=6$, $\sigma=0.1$ and (b) $D_{\max}=6$, $\sigma=0.2$. Note that the error range is $[10^{-5},5]$. In these two cases, MMR$+\texttt{fmincon}$ achieves an exact recovery rate of (a) $92\%$ and (b) $73\%$, while SNLSDP$+\texttt{fmincon}$ achieves (a) $64\%$ and (b) $28\%$.}\label{fig:ext1}
\end{figure}

 We conclude that both $\texttt{fmincon}$ and SA rely heavily on the initialization. MMR outperforms SNLSDP in terms of both mean and standard deviation. However, for MMR and SNLSDP, the standard deviation of the position error is much larger than the mean $\epsilon_p$ since both methods can fail for a few outlier examples. To take a closer look at the error distribution, we say that the sensors are `exactly recovered' if $\epsilon_p<10^{-5}$. Then MMR$+\texttt{fmincon}$ achieves an exact recovery rate of $92\%$, $73\%$ under the two settings respectively, compared to $64\%$, $28\%$ for SNLSDP$+\texttt{fmincon}$. Moreover we plot the histograms the empirical position error distributions for both MMR$+\texttt{fmincon}$ and SNLSDP$+\texttt{fmincon}$ in Fig.~\ref{fig:ext1}. Note that the histograms only include `failed' cases, i.e., cases where $\epsilon_p>10^{-5}$, for ease of readability. We observe that even when MMR$+\texttt{fmincon}$ does not yield exact recovery, the error is much smaller and more concentrated near zero than that of SNLSDP$+\texttt{fmincon}$.

\subsubsection{Varying the MMR depth} 
\begin{figure}[h]
\subfloat{{\includegraphics[width=0.5\textwidth]{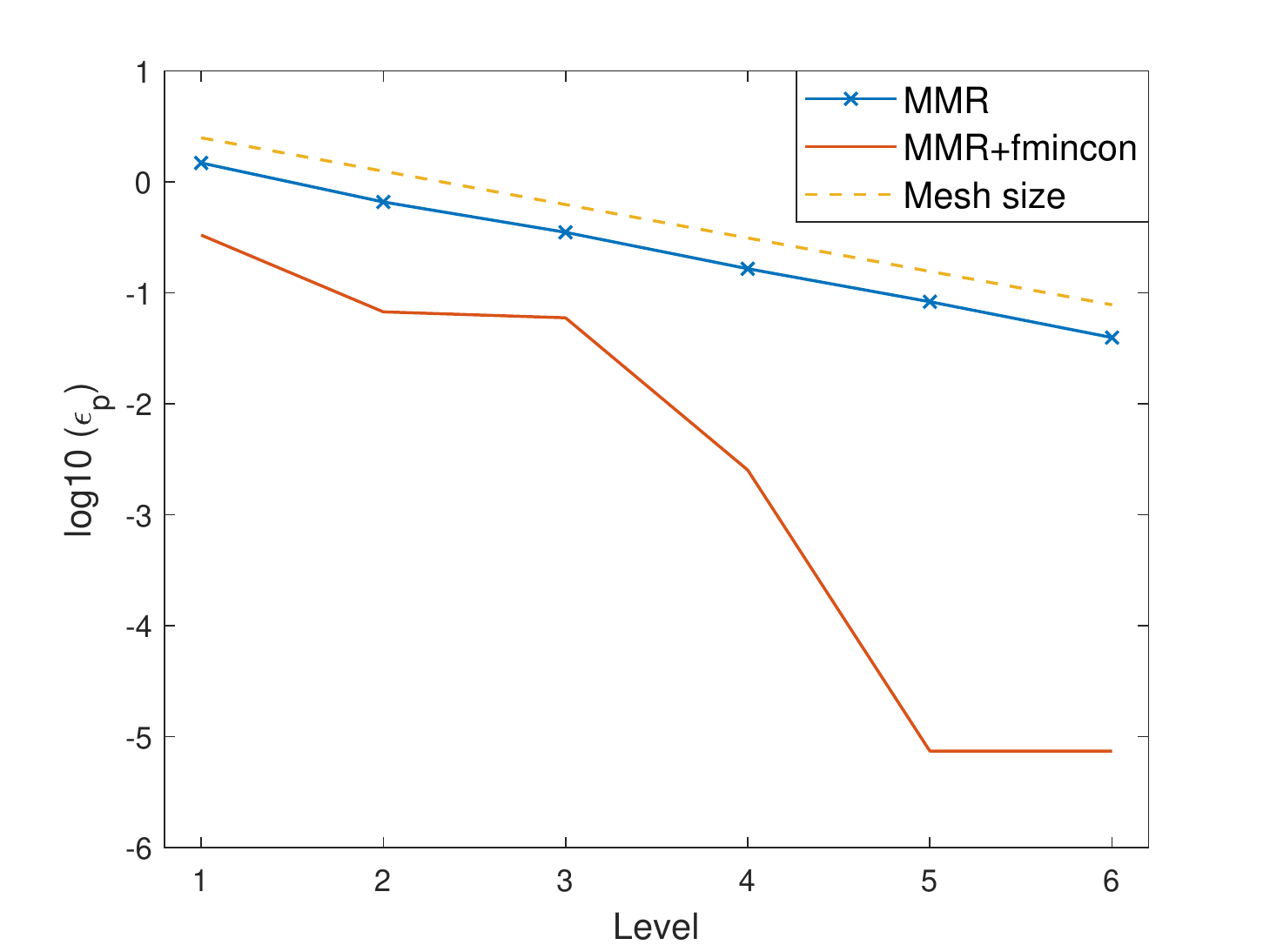}}}
\subfloat{{\includegraphics[width=0.5\textwidth]{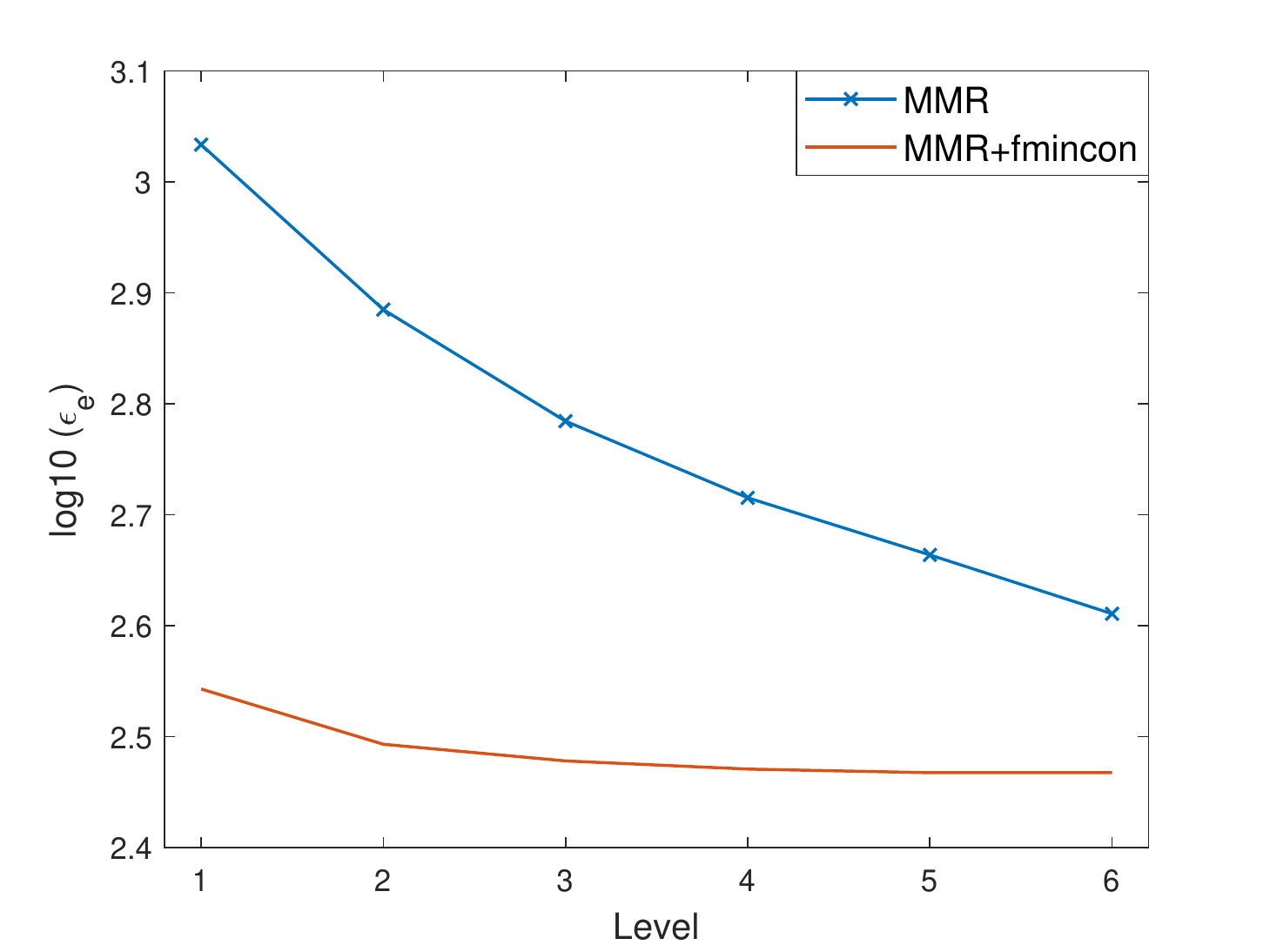}}}
  \caption{(a) Position error $\epsilon_{p}$ and (b) energy $\epsilon_{e}$ as functions of $K$. In (a), we include an additional reference line indicating the grid spacing at each level.}\label{fig:1}
\end{figure}

In this subsection we demonstrate how the performance of MMR depends on the number of levels $K$. We fix the parameters $\sigma=0.1$ and $D_{\max}=5$. We pick one example for which our algorithm succeeds.
In Fig.~\ref{fig:1} we plot the position error and energy against $K$. We observe that the decay of position error follows that of the grid spacing. We also observe that MMR+\texttt{fmincon} achieves perfect reconstruction (up to numerical error) by level 5. This observation validates the practical point that we can terminate MMR once we get good enough initial guess for \texttt{fmincon}.

\subsubsection{Varying $\sigma$ and $D_{\max}$}
\par In this subsection we vary the two model parameters $\sigma$ and $D_{\max}$, which determine, respectively, the proportion of contaminated entries in the distance matrix and how many measurements (corrupted or otherwise) are available.

We compare the output of the finest level $K=6$ of MMR to the ground truth particle positions. At the finest level, the grid spacing is $10/(2^{7})=0.078125$. For each particle, if the errors of both coordinates are less than the grid spacing we say that MMR was successful for this particle. The success rate is defined to be the proportion of successful particles out of the total 50 particles. In Tables \ref{table:2} and \ref{table:3} we present the success rate and the runtime of the algorithm (averaged over 30 independent realizations) across different settings.

\begin{table}[h]
\centering
 \caption{Success rate for different values of $\sigma$ and $D_{\max}$. Results are averaged over 30 independent realizations. }\label{table:2}
\begin{tabular}{|c|m{1cm}|m{1cm}|m{1cm}|m{1cm}|m{1cm}|}
\hline
 \diagbox{$\sigma$}{$D_{\max}$} &5& 6 & 7 & 8 & 9\\
\hline
0.1  & 90.00\% & 96.67\% & 93.33\% & 100\% & 100\%\\
\hline
0.2  & 50.00\% & 70.00\% & 86.67\% & 93.33\% & 96.67\%\\
\hline
0.3  & 3.33\% & 60.00\% & 73.33\% & 96.67\% & 96.67\% \\
\hline
\end{tabular}
\end{table}

\begin{table}[h]
\centering
 \caption{Runtime (hours) for different values of $\sigma$ and $D_{\max}$. Results are averaged over 30 independent realizations.}\label{table:3}
\begin{tabular}{|c|m{1cm}|m{1cm}|m{1cm}|m{1cm}|m{1cm}|}
\hline
 \diagbox{$\sigma$}{$D_{\max}$} & 5 & 6 & 7 & 8 & 9\\
\hline
0.1  & 2.3882 & 2.1316 & 2.0241  & 2.0187 & 2.0708   \\
\hline
0.2  & 2.3923 & 2.5364 & 2.5277 & 2.2704 & 2.2771 \\
\hline
0.3  & 1.8739 & 1.9579 & 2.0375 & 2.1147 & 2.1741\\
\hline
\end{tabular}
\end{table}

As we can observe from the tables, in harder cases (larger $\sigma$, smaller $D_{\max}$) the success rate usually drops. However the total runtime remains roughly the same. We observed empirically that the runtime mainly depends on the threshold value $\eta^{(k)}$ and the required accuracy of the SDP solver. We can potentially improve the success rate by propagating larger parts in the $k$-th level partition to the $(k+1)$-th level. 

\subsection{Lennard-Jones clusters}\label{sec:42}
Here we consider particles in the plane interacting via the Lennard-Jones potential and solve for the globally optimal configuration. The pairwise cost, i.e., the Lennard-Jones (LJ) potential, is given by 
\begin{align}
    H_{ij}(x_{i},x_{j})=\epsilon\left[\left(\frac{r_{ij}}{||x_{i}-x_{j}||_{2}}\right)^{12}-2 \left(\frac{r_{ij}}{||x_{i}-x_{j}||_{2}}\right)^{6}\right].
\end{align}
Here $x_i,x_j \in \mathbb{R}^2$, and $r_{ij}$ is the distance at which the LJ potential reaches its minimum value $-\epsilon$ for the pair $i,j$. In all the tests below, we set $\epsilon=1$. In the classic `LJ cluster,' the optimal distances $r_{ij}$ are the same for all pairs of particles, i.e., $r_{ij}=r$. In our numerical tests, we test two cases: (1) the classic LJ cluster setting $r_{ij}=r$, which we call the `symmetric case,' and (2) the asymmetric setting where the $r_{ij}$ are in general different, which we call the `asymmetric case.' In each case, we detail the application of our method and compare performance with that of the successful basin-hopping approach~\cite{basinhopping}.

Note that in principle the domain of the LJ potential is noncompact. However, the LJ potential approaches 0 when the distance between two particles tends to infinity, so solutions with particles very far from one another cannot be locally or globally optimal. We can then restrict our search domain based on the optimal distance $r_{ij}$ and the number of particles $N$. We will restrict the domain to square regions as in the SNL experiments and pursue the same discretization and coarsening strategies.

In this problem we pick a minimal upper bound value for the coarsest level $u^{(1)}=u_{\mathrm{min}}$ and gradually increase the upper bound to $1$ as we refine, following the formula
\begin{align}
    u^{(k+1)}=u^{(k)}+\alpha(1-u^{(k)}),\ k=1,\ldots,K-1
\end{align}
where $\alpha\in[0,1]$ is the increasing rate. We couple the upper bound to the threshold parameter by setting
\begin{align}
    \eta^{(k)}=\beta u^{(k)},
\end{align}
and we will control the new ratio parameter $\beta$.

In our refinement steps, we use a Neumann neighborhood to define adjacency (so each non-boundary grid point has 4 adjacent sites). We perform 3 iterations of refinement at each level, terminating prematurely if the selected neighborhoods in two consecutive iterations remain the same.

\subsubsection{Symmetric Case}\label{sec:421}
\par In the symmetric case, we set $r_{ij}=1$ for all $i,j$. This puts us in the setting of Section~\ref{sec:symmetricGeneral}, where $H_{ij} = H$, and $H$ is symmetric as a matrix. Accordingly, we consider the 2-marginal relaxation \eqref{eq:sdp_symm}. As mentioned in Section~\ref{sec:symmetricGeneral}, since the diagonal of $H$ is $+\infty$, we replace it with zeros and include in our SDP the constraint that $\mathrm{diag}(\gamma) = 0$. For simplicity we omit most the details of our multiscale strategy for solving this SDP, which follow similarly the details outlined above for asymmetric problems. In the symmetric case, we need only consider a single 1-marginal (as opposed to $N$ separately), though the support of this 1-marginal should be $N$ times as expansive as in the asymmetric case.

At the coarsest level of our multiscale algorithm we consider a $16 \times 16$ regular grid with upper-bounding parameter $u^{(1)} = 0.1$. Our threshold at the first level is $\eta^{(1)} = 0.002$. For $k>1$, we take $u^{(k)} = 1$ and $\eta^{(k)} = 0.02$, and we consider $K=6$ total levels. The above choices of $\eta^{(k)}$ correspond to a choice of $\beta=0.02$ for all levels in our earlier notation. As in the above experiments, each grid point in a given level is the parent of a $2\times 2$ block of equispaced grid points in the next level. Hence the finest level is a grid of size $2^9 \times 2^9$. Throughout the following, the tolerance of the SDP solver is set to be $10^{-6}$.

For this problem, we have a strong intuition that the support of the optimal 1-marginal is nearly a subset of a hexagonal lattice~\cite{blanc2015crystallization}. Therefore we remove the degeneracy with respect to rigid motions by fixing the first three `anchor' particles as vertices of an equilateral triangle in the center of the search domain. Concretely, this corresponds to fixing three entries of $\rho$ in \eqref{eq:sdp_symm} to be $1/N$.

First we shall focus on the $N=13$ case as a detailed illustration of our method. Then we compare the performance more systematically to the basin-hopping (BH) algorithm for larger particle numbers.

\begin{figure}[!htb]
\subfloat[$k=3$]{{\includegraphics[trim=0.7cm 0.7cm 0.7cm 0.7cm, clip=true,width=0.25\textwidth]{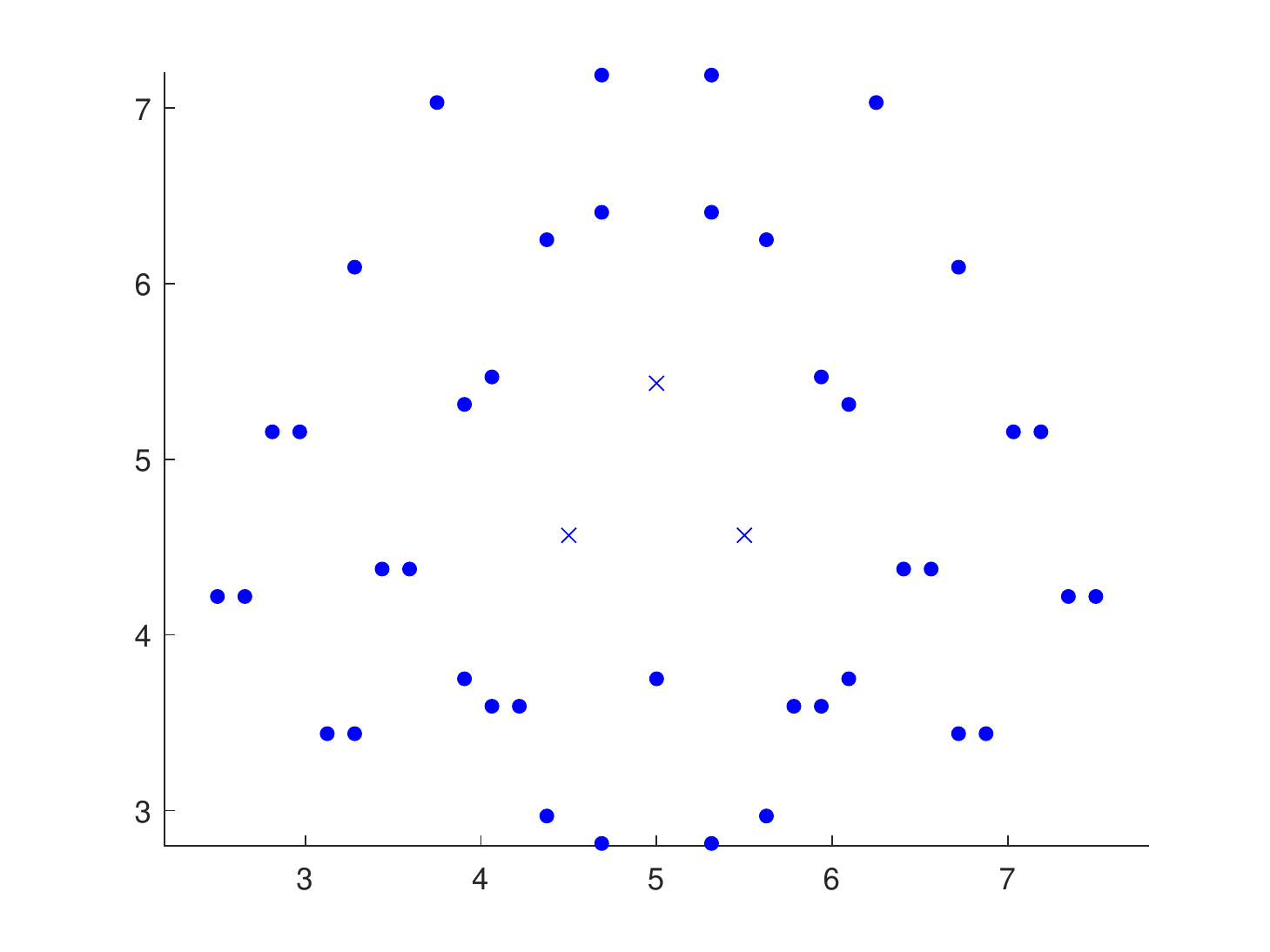}}}
\subfloat[$k=4$]{{\includegraphics[trim=0.7cm 0.7cm 0.7cm 0.7cm, clip=true,width=0.25\textwidth]{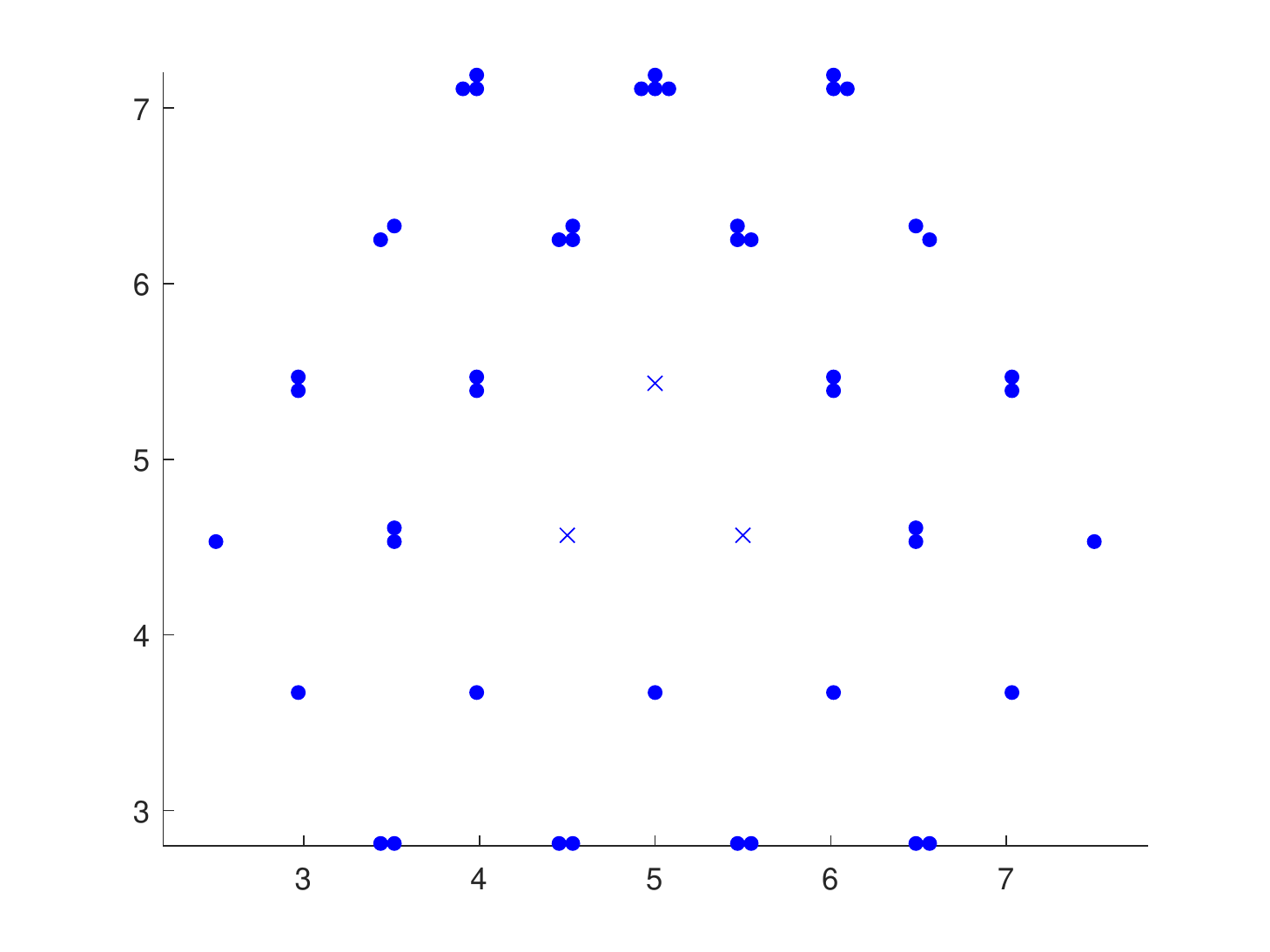}}}
\subfloat[$k=5$]{{\includegraphics[trim=0.7cm 0.7cm 0.7cm 0.7cm, clip=true,width=0.25\textwidth]{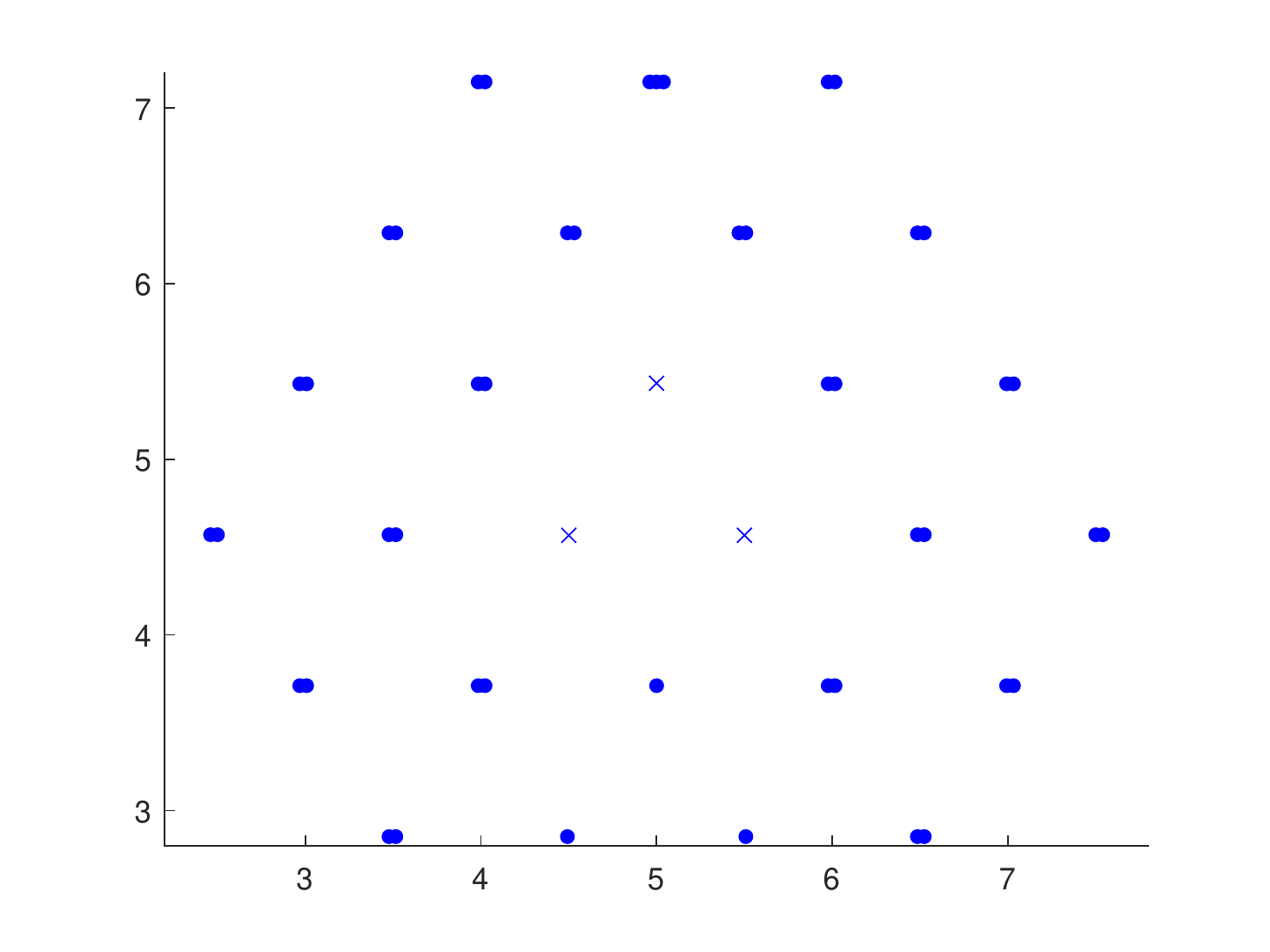}}}
\subfloat[$k=6$]{{\includegraphics[trim=0.7cm 0.7cm 0.7cm 0.7cm, clip=true,width=0.25\textwidth]{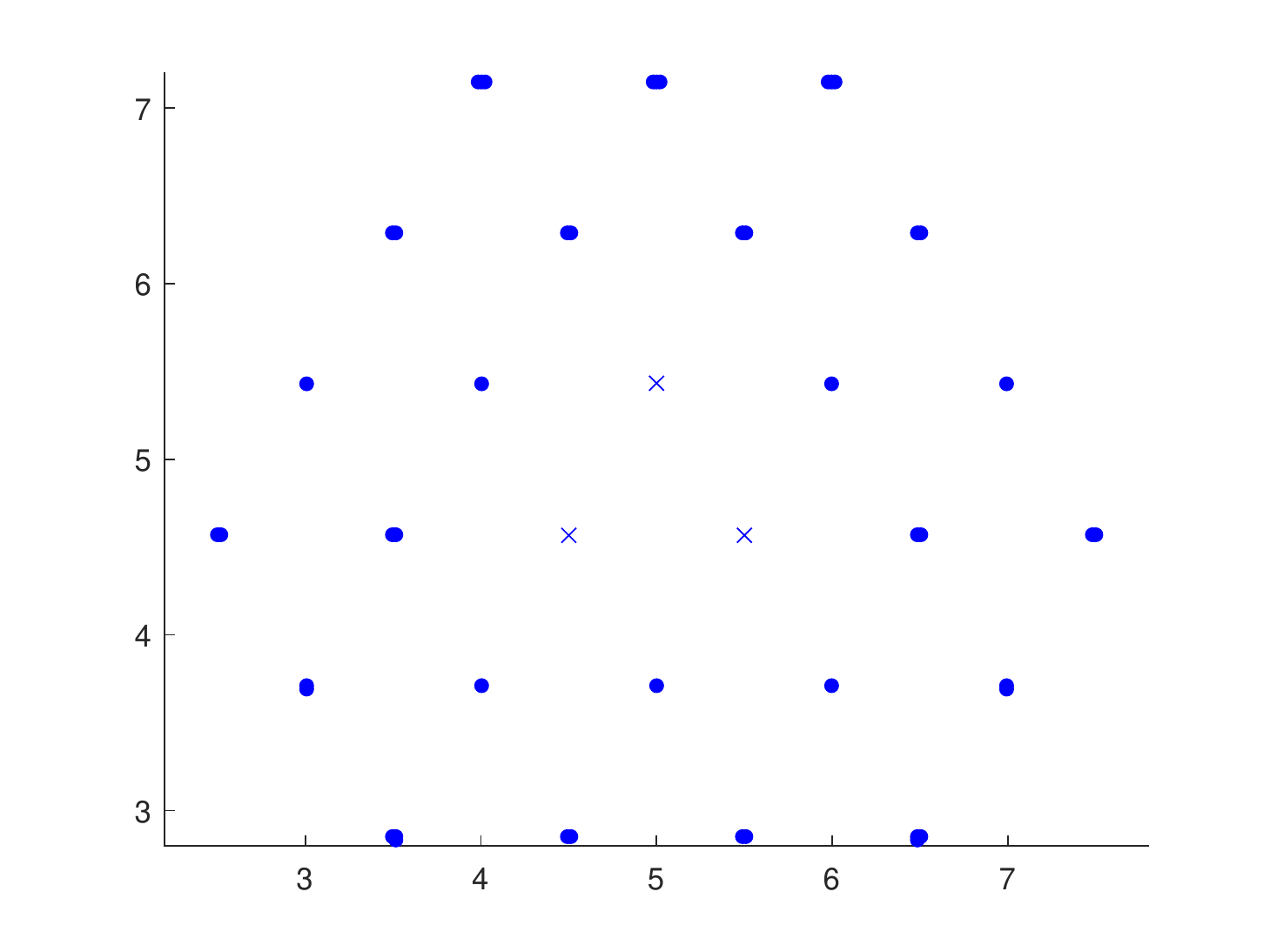}}}
  \caption{$N=13$. MMR output at levels (a) $k=3$, (b) $k=4$, (c) $k=5$, (d) $k=6$. Blue points indicate possible particle positions selected by MMR, and blue crosses indicate the anchor particles.}\label{fig:2}
\end{figure}

In the $N=13$ case, the search domain is set to be $[0,10]^{2}$, so the anchors are fixed at $x_{1}=(4.5,5-\frac{\sqrt{3}}{4})$, $x_{2}=(5.5,5-\frac{\sqrt{3}}{4})$, $x_{3}=(5,5+\frac{\sqrt{3}}{4})$.
We visualize the output of last four levels of MMR in Fig.~\ref{fig:2}, illustrating how MMR generates guesses and propagates them to the finer levels. In each plot, the blue points indicate the possible particle positions selected by MMR, and the three blue crosses indicate the anchor particles.

We see that the $k=3$ level is still too coarse for MMR to retrieve an accurate prediction of the support of optimal solutions, but at $k=4$ the support starts to separate into distinct clusters around possible particle locations. At the finest level $k=6$, we can see 27 clusters of support, and all of which are well-separated from each other. The fact that we obtain more clusters than $N=13$ particles is due to the existence of multiple near-optimal configurations of the LJ cluster with $N=13$ particles. Visually they correspond to different subsets of the hexagonal lattice recovered in Fig.~\ref{fig:2}(d) by MMR.

\begin{figure}[!htb]
\subfloat[Energy: -55.5767 (\textbf{-55.5889}) ]{{\includegraphics[trim=0.7cm 0.7cm 0.7cm 0.7cm, clip=true,width=0.33\textwidth]{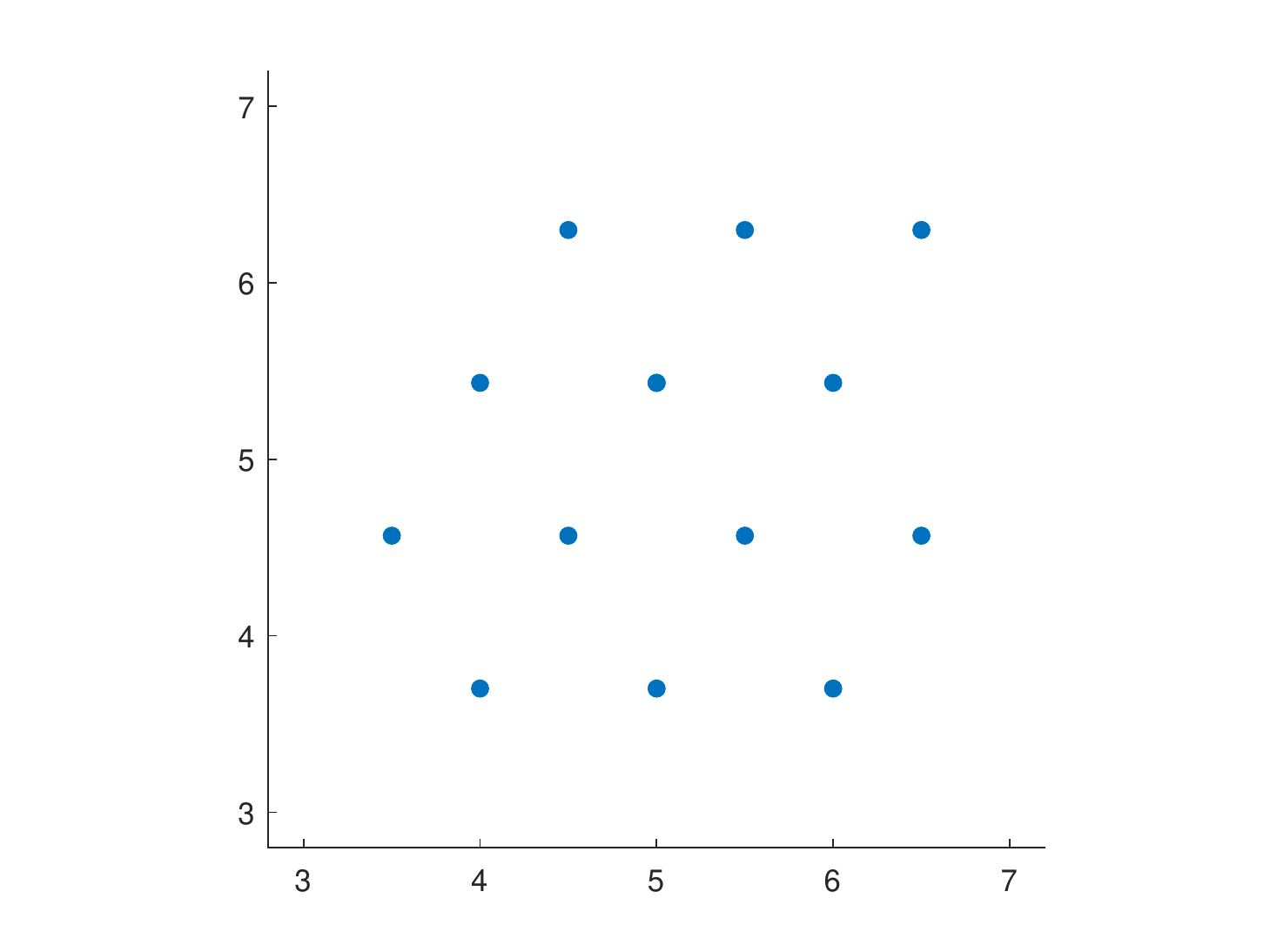}}}
\subfloat[Energy: -55.5735 (\textbf{-55.5889}) ]{{\includegraphics[trim=0.7cm 0.7cm 0.7cm 0.7cm, clip=true,width=0.33\textwidth]{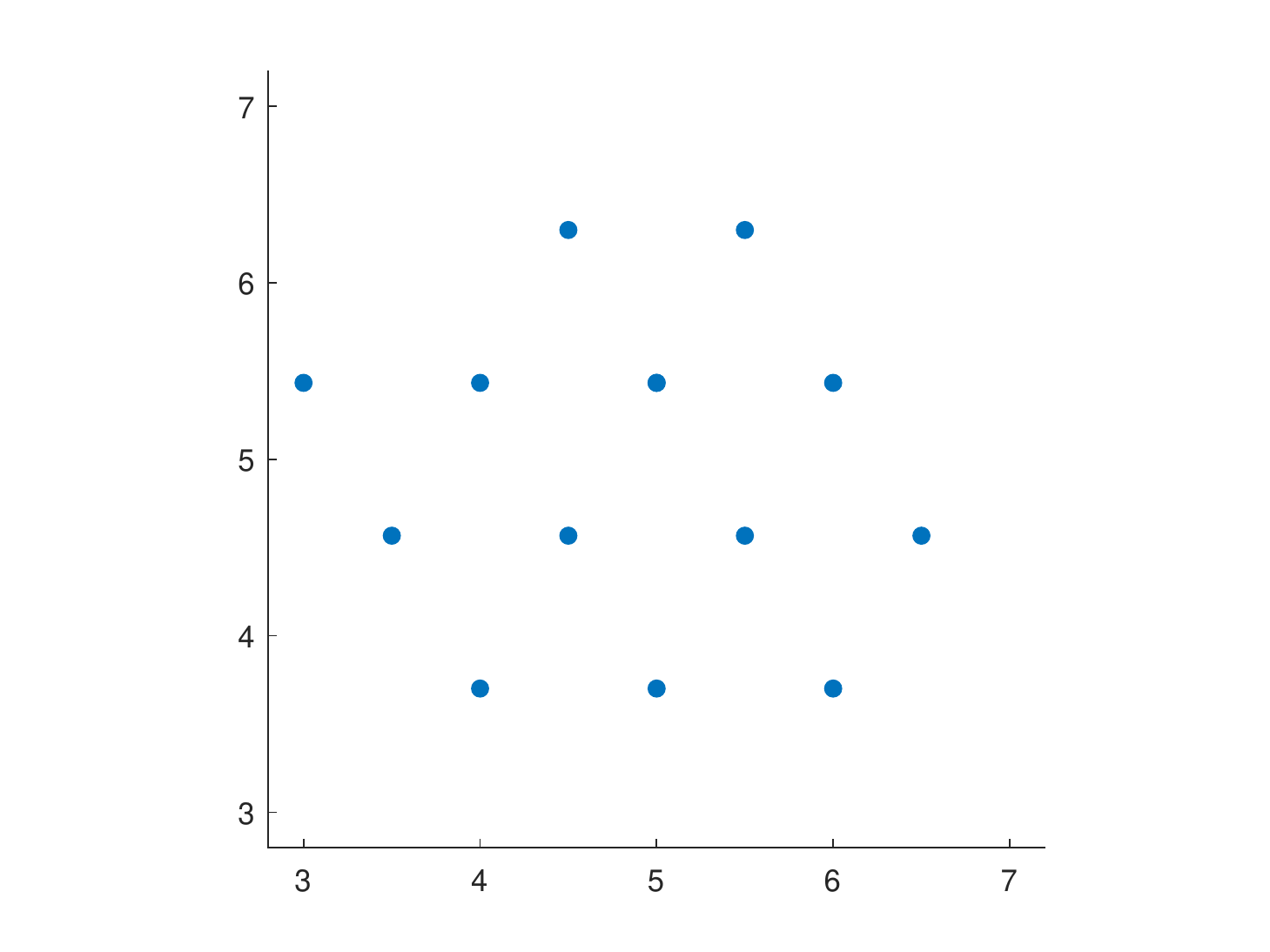}}}
\subfloat[Energy: -55.5735 (\textbf{-55.5889}) ]{{\includegraphics[trim=0.7cm 0.7cm 0.7cm 0.7cm, clip=true,width=0.33\textwidth]{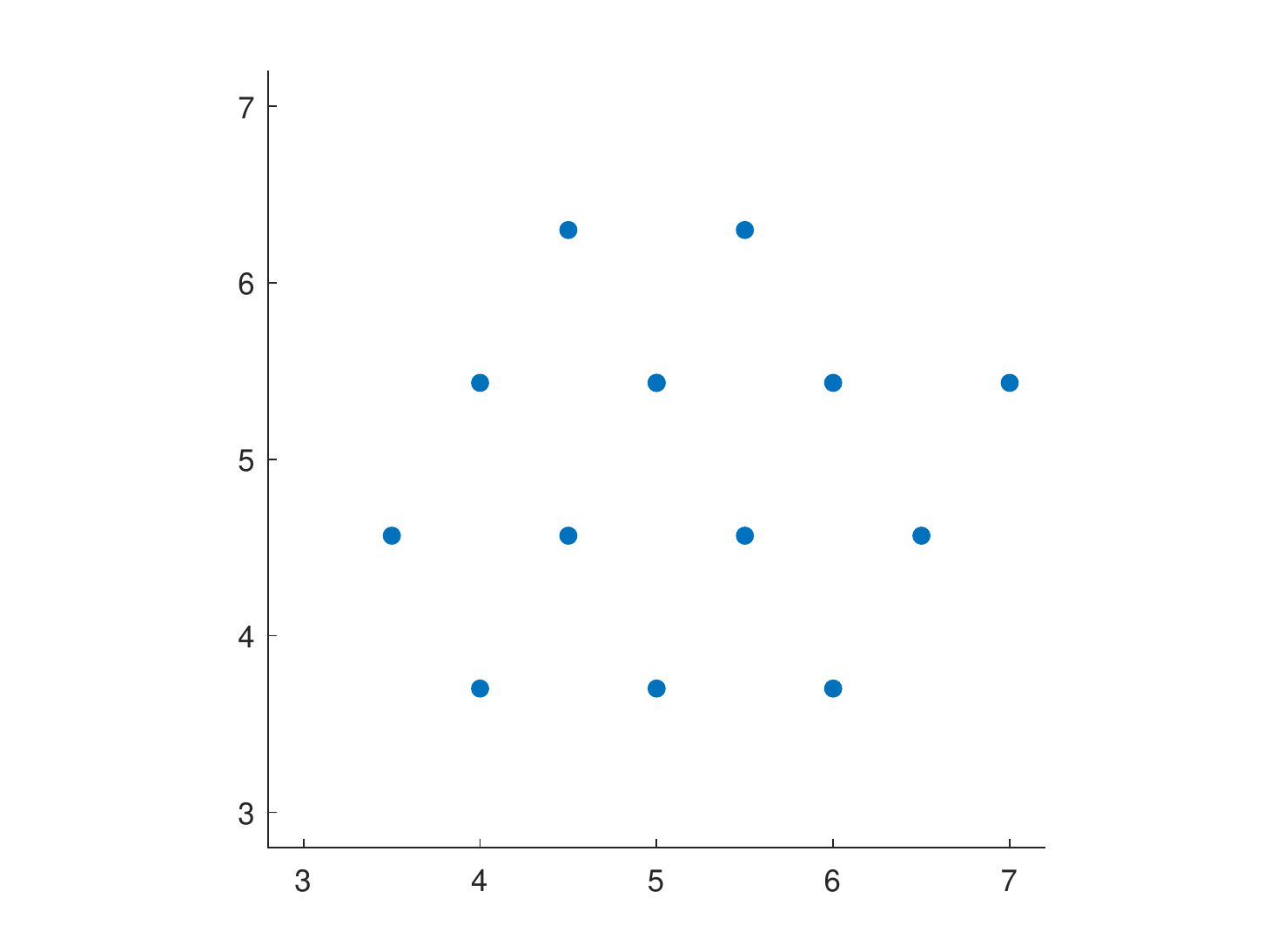}}}
  \caption{Three configurations generated by the randomized method with small noise, using the same MMR output but with different instantiations of noise in the cost matrix. We also show the energy of the configuration and the energy after post-processing by \texttt{fmincon} (\textbf{bold face}) to remove discretization error.}\label{fig:4}
\end{figure}

\begin{figure}[!htb]
\centering
\subfloat[][Energy:-55.5187

(\textbf{-55.5314})]{{\includegraphics[trim=0.7cm 0.7cm 0.7cm 0.7cm, clip=true,width=0.25\textwidth]{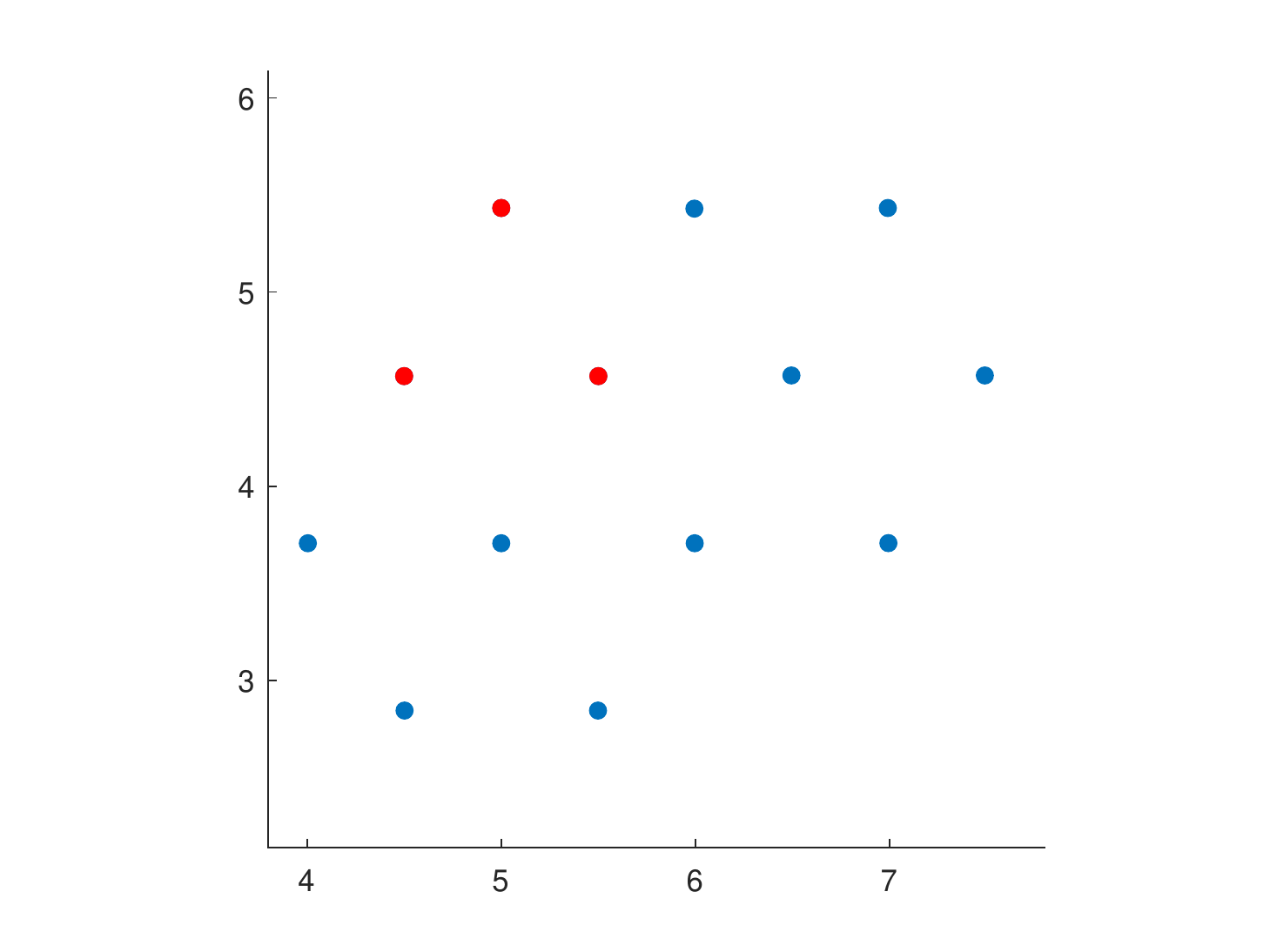}}}
\subfloat[][Energy:-51.2323

(\textbf{-51.2454})]{{\includegraphics[trim=0.7cm 0.7cm 0.7cm 0.7cm, clip=true,width=0.25\textwidth]{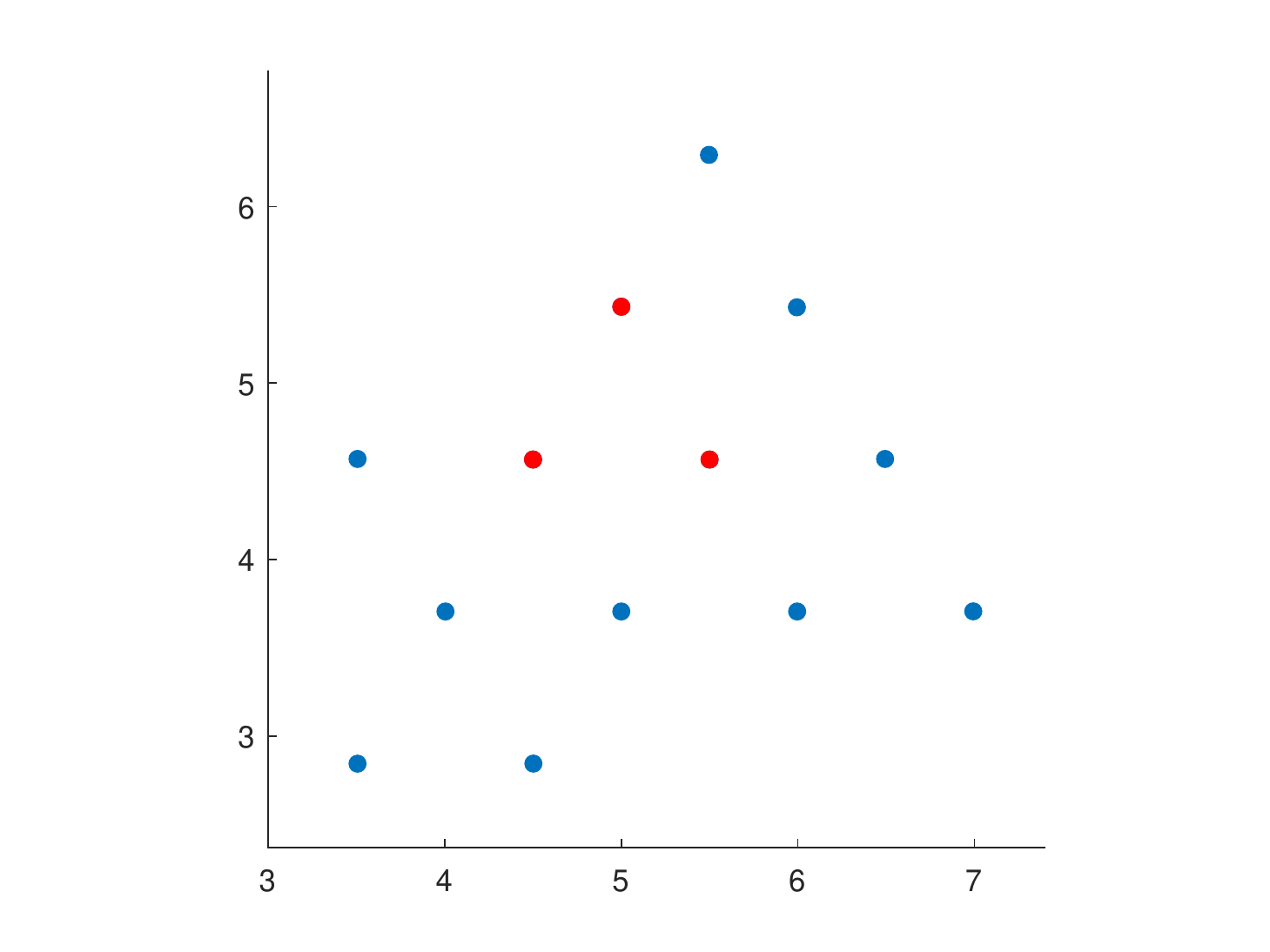}}}
\subfloat[][Energy:-53.5283

(\textbf{-55.5965})]{{\includegraphics[trim=0.7cm 0.7cm 0.7cm 0.7cm, clip=true,width=0.25\textwidth]{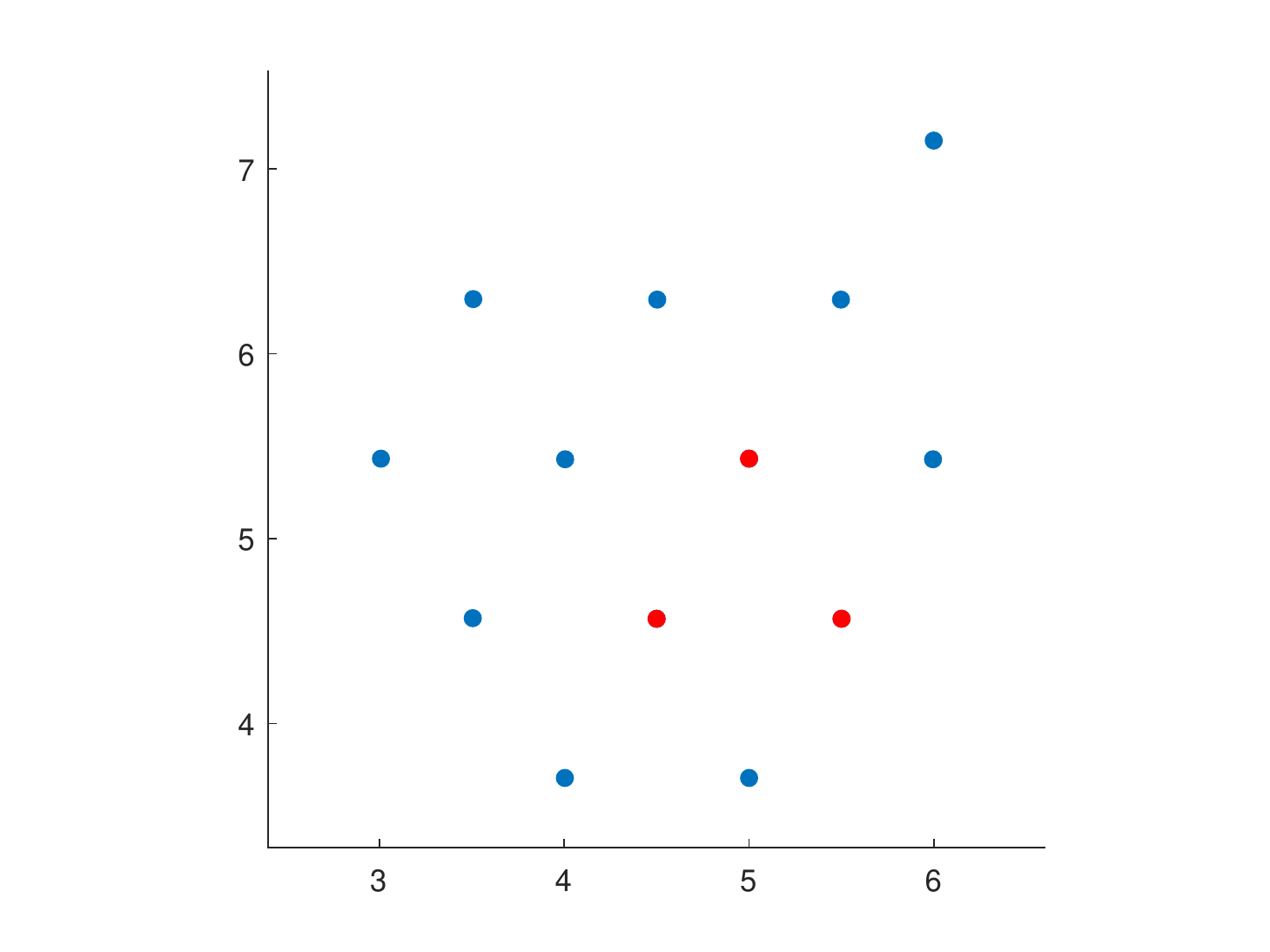}}}
\subfloat[][Energy:-53.5805

(\textbf{-53.5924})]{{\includegraphics[trim=0.7cm 0.7cm 0.7cm 0.7cm, clip=true,width=0.25\textwidth]{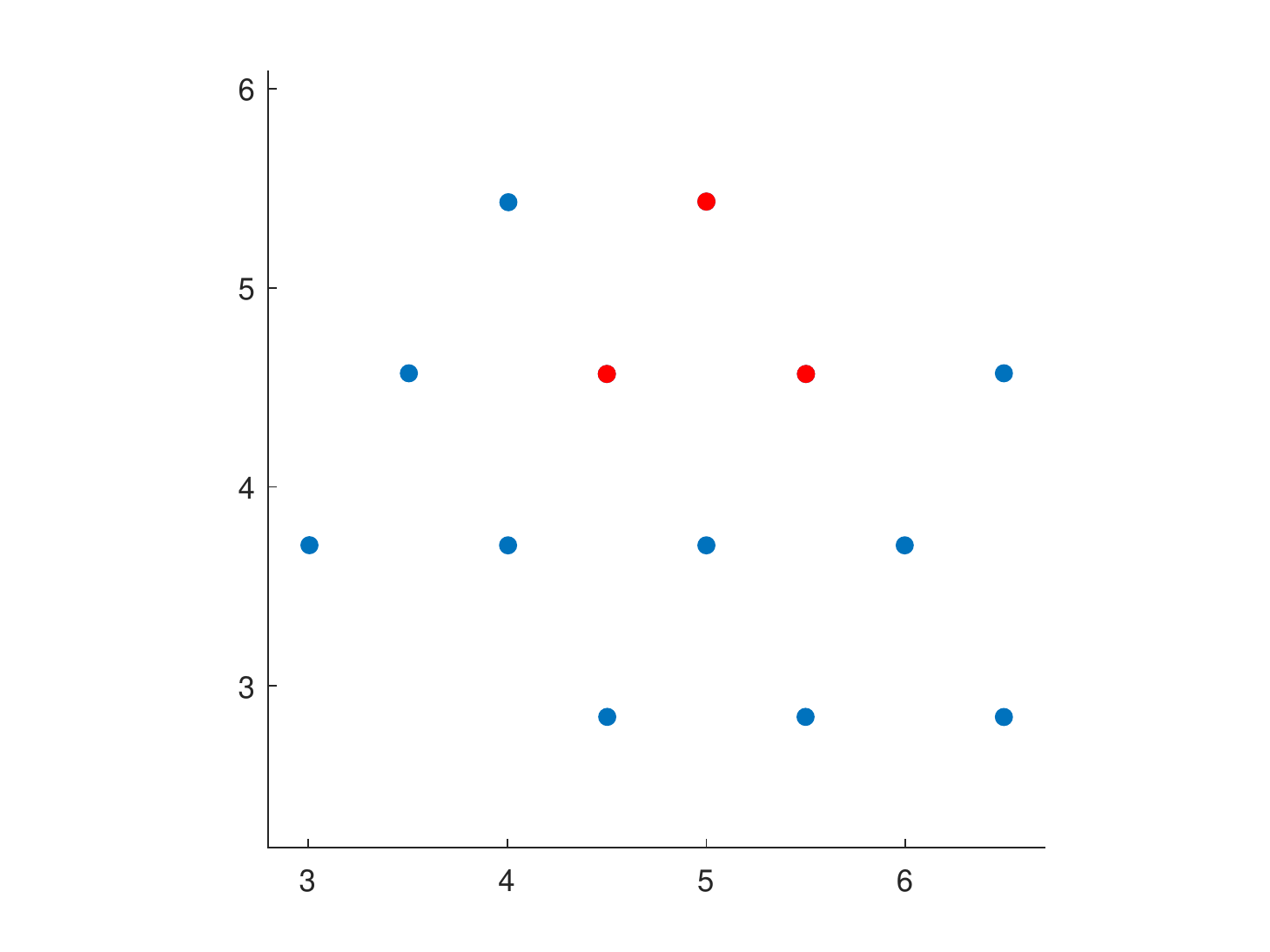}}}
  \caption{Four configurations generated by the randomized method with moderate noise, using the same MMR output but with different instantiations of noise in the cost matrix. We also show the energy of the configuration and the energy after post-processing by \texttt{fmincon} (\textbf{bold face}) to remove discretization error. The red points indicate the three fixed anchor particles.}\label{fig:5}
\end{figure}

Accordingly, we use the randomized method proposed in Section~\ref{sec:multiple soln} to explore the near-optimal configurations, i.e., we add a noise term to the cost in  \eqref{eq:sdp_symm}, restricted to grid points yielded by MMR at the finest level. Then we extract a solution via the top eigenvector of the entrywise nonnegative, semidefinite matrix $\Lambda$~\eqref{eq:lambda}, which should correspond to the outer product $\rho \rho^\top$ in the case of exact recovery of the 1-marginal. The eigenvector is nonnegative by Perron-Frobenius, and we generate a candidate solution by selecting its top $N$ entries. The results following 3 independent noise samples are reported in Fig.~\ref{fig:4}. There is almost no difference between the energies of different solutions. Next we show that the sampling procedure can also explore near global optima. By adding a larger noise term, we obtain configurations that have slightly higher energies, as shown in Fig.~\ref{fig:5}. 

Next we compare MMR with BH for particle numbers $N=7,13,20,30$. We take $[0,10]^2$ as the particle domain for $N=7,13,20$, and for $N=30$ we enlarge the the domain to $[0,14]^2$ since more particles are involved, though all other parameters remain the same. After obtaining the MMR output, we then used the randomized method of Section~\ref{sec:multiple soln} to generate 10 sample solutions. These solutions are further refined via local optimization of the LJ potential, mitigating the effect of discretization error.

Meanwhile, for BH, we choose the temperature parameter to be $T=0.2$~\cite{basinhopping}. We also choose \texttt{L-BFGS-B} as the local optimizer called by BH, and we initialize by selecting particles independently and uniformly from the domain $[0,5]^2$. We tuned these user choices for best performance on this example. Here we execute $10$ independent BH runs. We ran BH for as many iterations as it takes to match the runtime of MMR. For both algorithms, we pick the best energies obtained from the 10 samples. The minimal energies are reported in Table \ref{table:7}.
\begin{table}[!htb]
\centering
 \caption{Minimal energies of MMR+\texttt{fmincon} and BH for $N=7,13,20,30$.
 }\label{table:7}
\begin{tabular}{|c|c|c|c|c|c|}
\hline
 $N$ & $7$ & $13$ & $20$ & $30$\\
\hline
MMR & -25.0666 & -55.5889 & -95.1684 & -154.7580\\
\hline
BH  & -25.0666 & -55.6052 & -95.1847 & -154.7772\\
\hline
\end{tabular}
\end{table}

We see that the energy difference between MMR and BH is within $0.05\%$, and the optimal solutions furnished by the two algorithms are visually congruent, as shown in Fig.~\ref{fig:14} for the $N=13$ case and in Appendix~\ref{app:1} for the cases $N=20,30$. Hence the small gaps in energy appear to be attributable to numerical error in the local optimizer. (Note that highly accurate local optimization is nontrivial due to the singularities of the LJ potential.) In Appendix~\ref{app:1} we also visualize several MMR sample configurations for the cases $N=20,30$ to further demonstrate the exploration of nearly-optimal configurations.
\begin{figure}[!htb]
\centering
\subfloat[MMR optimal solution configuration with energy $-55.5889$]{{\includegraphics[trim=0.7cm 0.7cm 0.7cm 0.7cm, clip=true,width=0.3\textwidth]{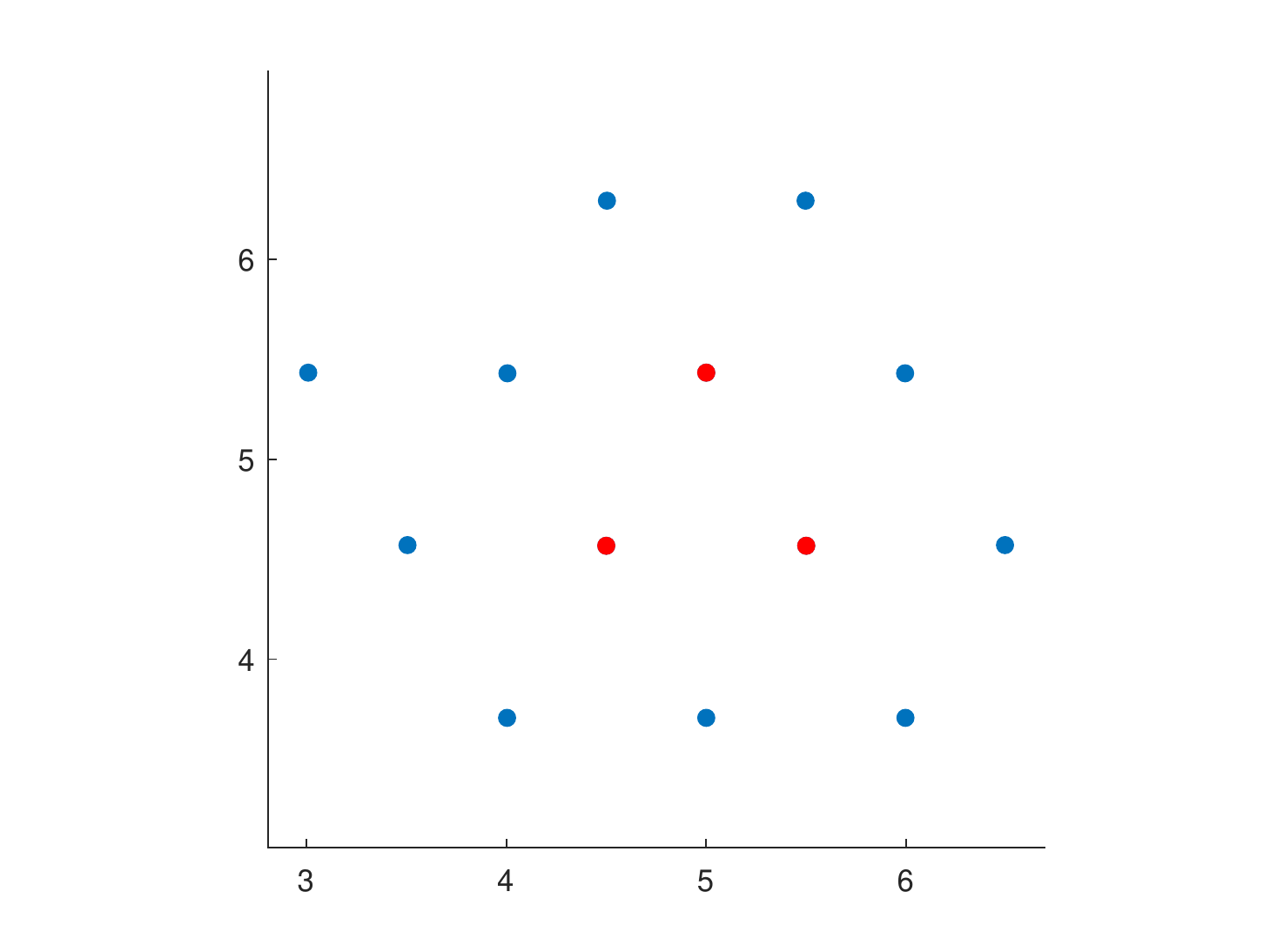}}}\ \ 
\subfloat[BH optimal solution configuration with energy $-55.6052$]{{\includegraphics[trim=0.7cm 0.7cm 0.7cm 0.7cm, clip=true,width=0.3\textwidth]{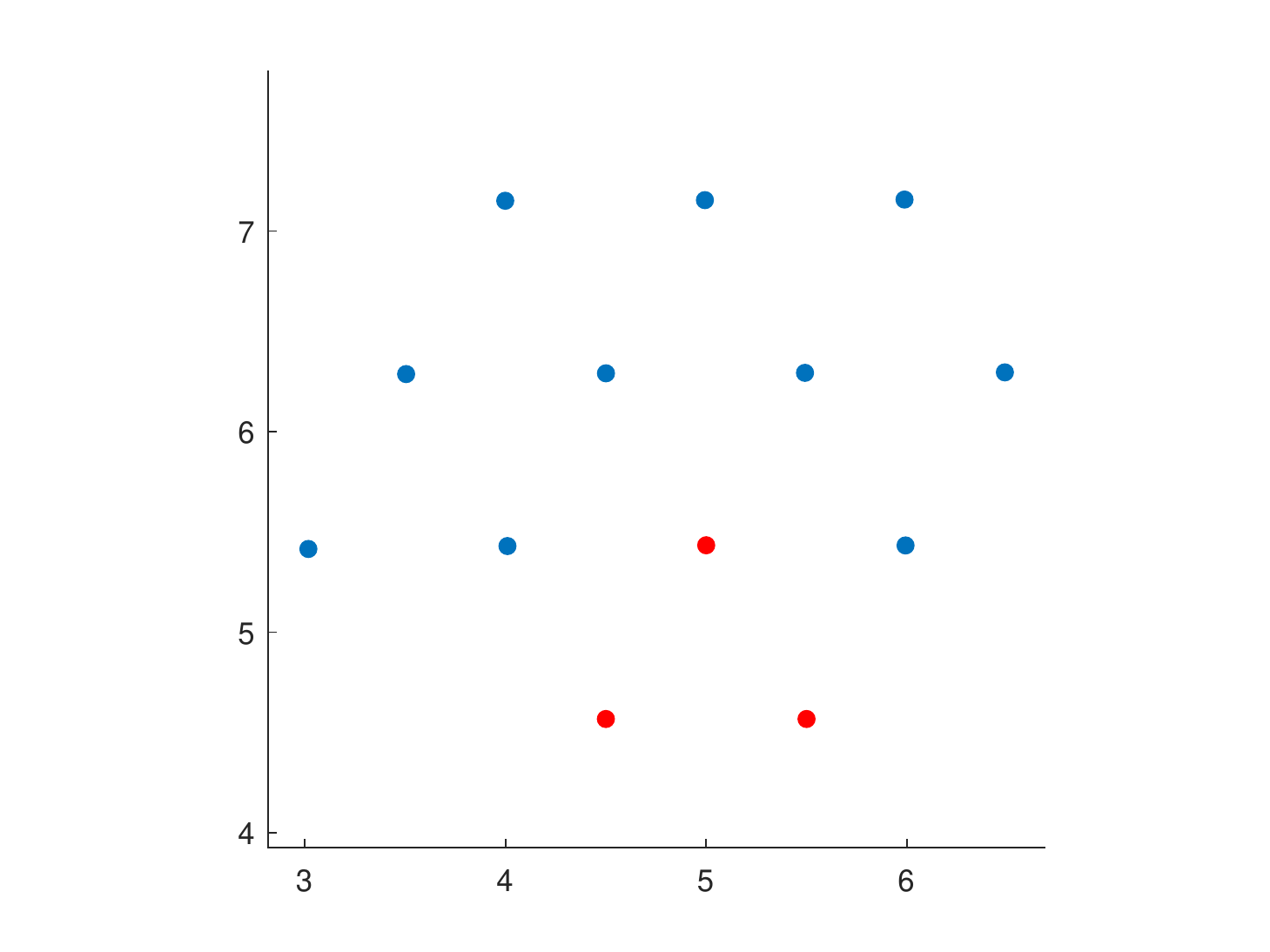}}}
  \caption{Illustration of MMR and BH optimal solution configuration. The red points indicate three fixed anchor particles. The two solution configurations are nearly congruent under rotation and reflection.}\label{fig:14}
\end{figure}

\subsubsection{Asymmetric Case}
For the asymmetric case, we independently sample $r_{ij} \sim \text{Unif}(0.5,1.5)$.  We examine the case with particle number $N=13$. For the particle domain we choose $\mathcal{X}=[0,10]^{2}$. In this example, it is difficult to fix the positions of any anchor particles, so we resort to the following completely general method for removing the degeneracy of the cost with respect to rigid motions in the plane. We fix the first particle to be the center of the domain, i.e.,  $x_{1}=(5,5)$. Then we constrain the second particle to have the same $y$-coordinate as the first particle but to have larger $x$-coordinate, i.e., $[x_{2}]_1 \geq 5$, $[x_{2}]_2 = 5$. Finally we constrain the third particle to have $y$-coordinate larger than the first two, i.e., $[x_{3}]_2 \geq 5$.

We again use regular grids at each level, in which each point is a parent of an equispaced $2\times 2$ block. At the coarsest level $k=1$, our grid is $8 \times 8$. We set the initial upper bound to be $u_{\mathrm{min}}=0.2$ and increase it with rate $\alpha=0.8$. The threshold $\eta^{(k)} = \beta u^{(k)}$ is fixed by the choice $\beta = 0.01$. 

\begin{figure}[!htb]
\subfloat[$k=4$]{{\includegraphics[trim=0.7cm 0.7cm 0.7cm 0.7cm, clip=true,width=0.250\textwidth]{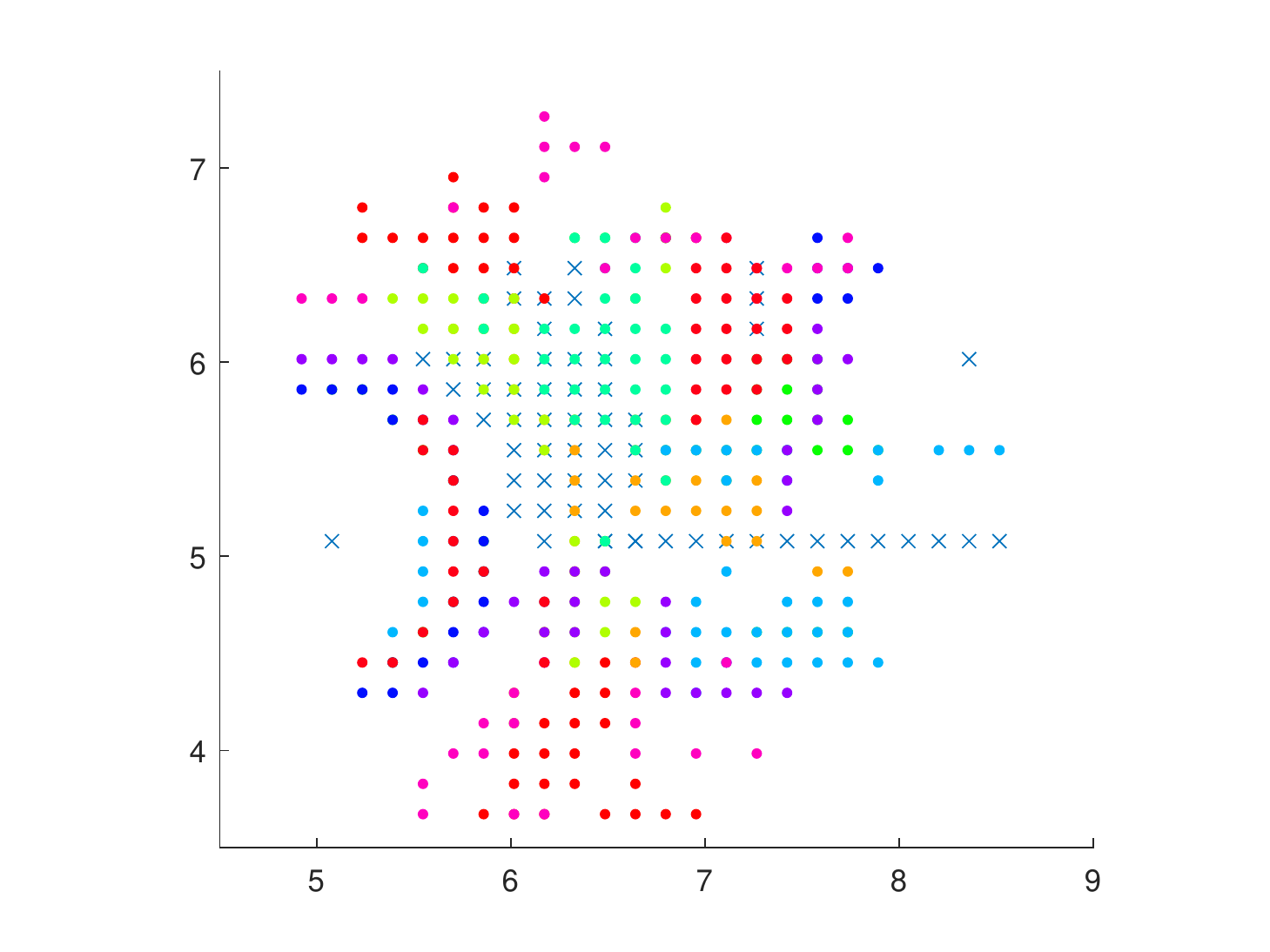}}}
\subfloat[$k=5$]{{\includegraphics[trim=0.7cm 0.7cm 0.7cm 0.7cm, clip=true,width=0.250\textwidth]{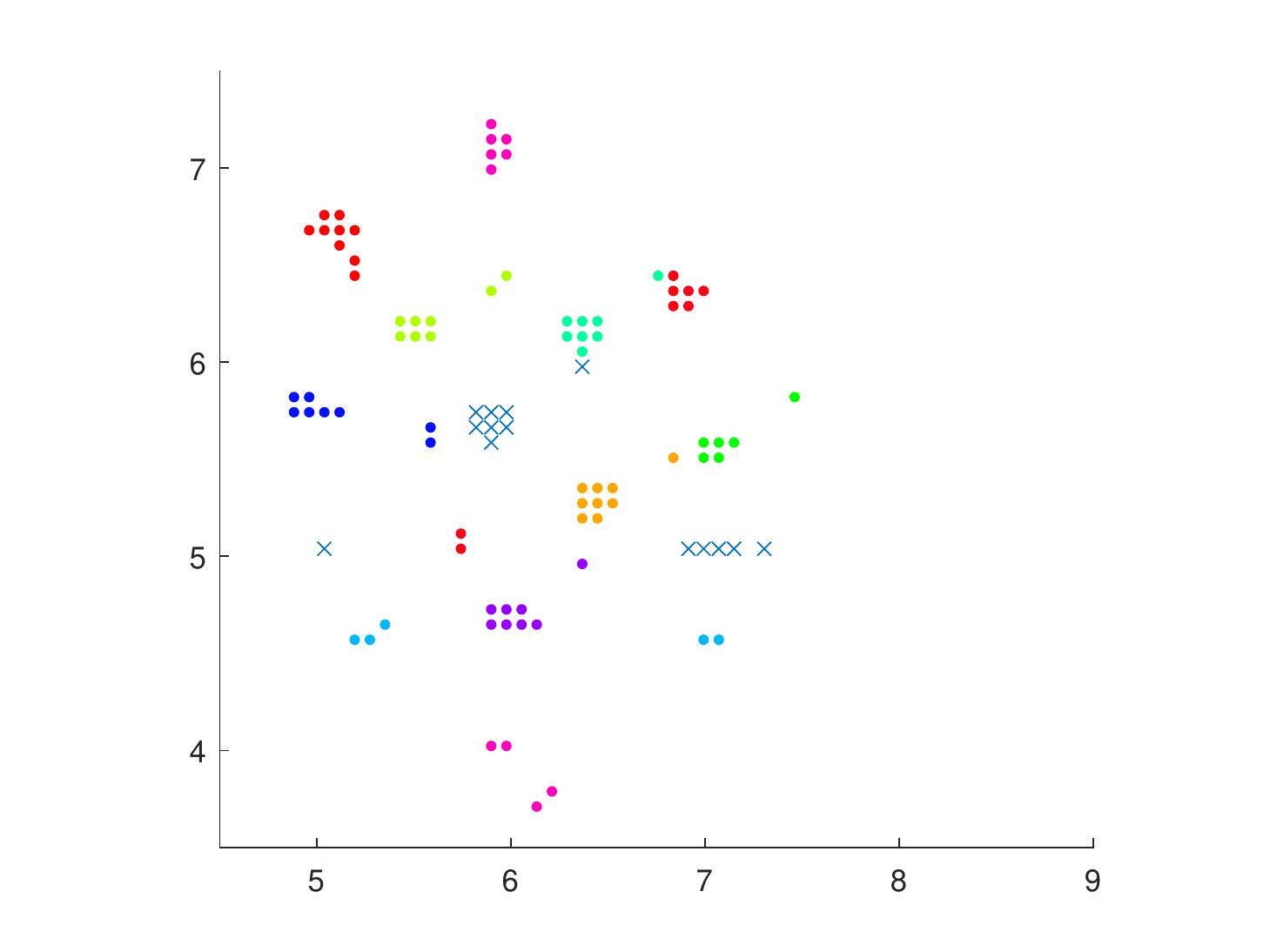}}}
\subfloat[$k=6$]{{\includegraphics[trim=0.7cm 0.7cm 0.7cm 0.7cm, clip=true,width=0.250\textwidth]{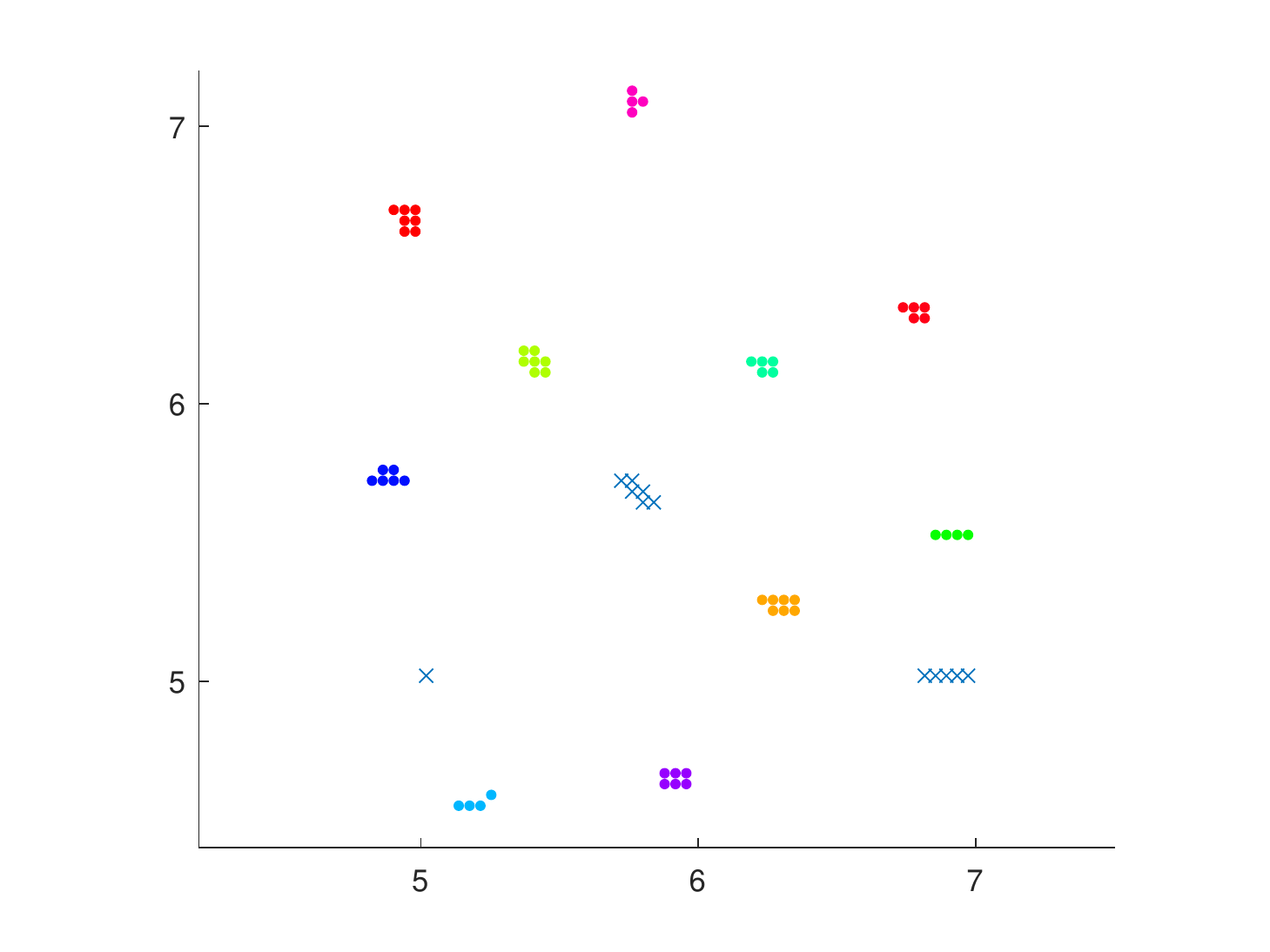}}}
\subfloat[$k=7$]{{\includegraphics[trim=0.7cm 0.7cm 0.7cm 0.7cm, clip=true,width=0.250\textwidth]{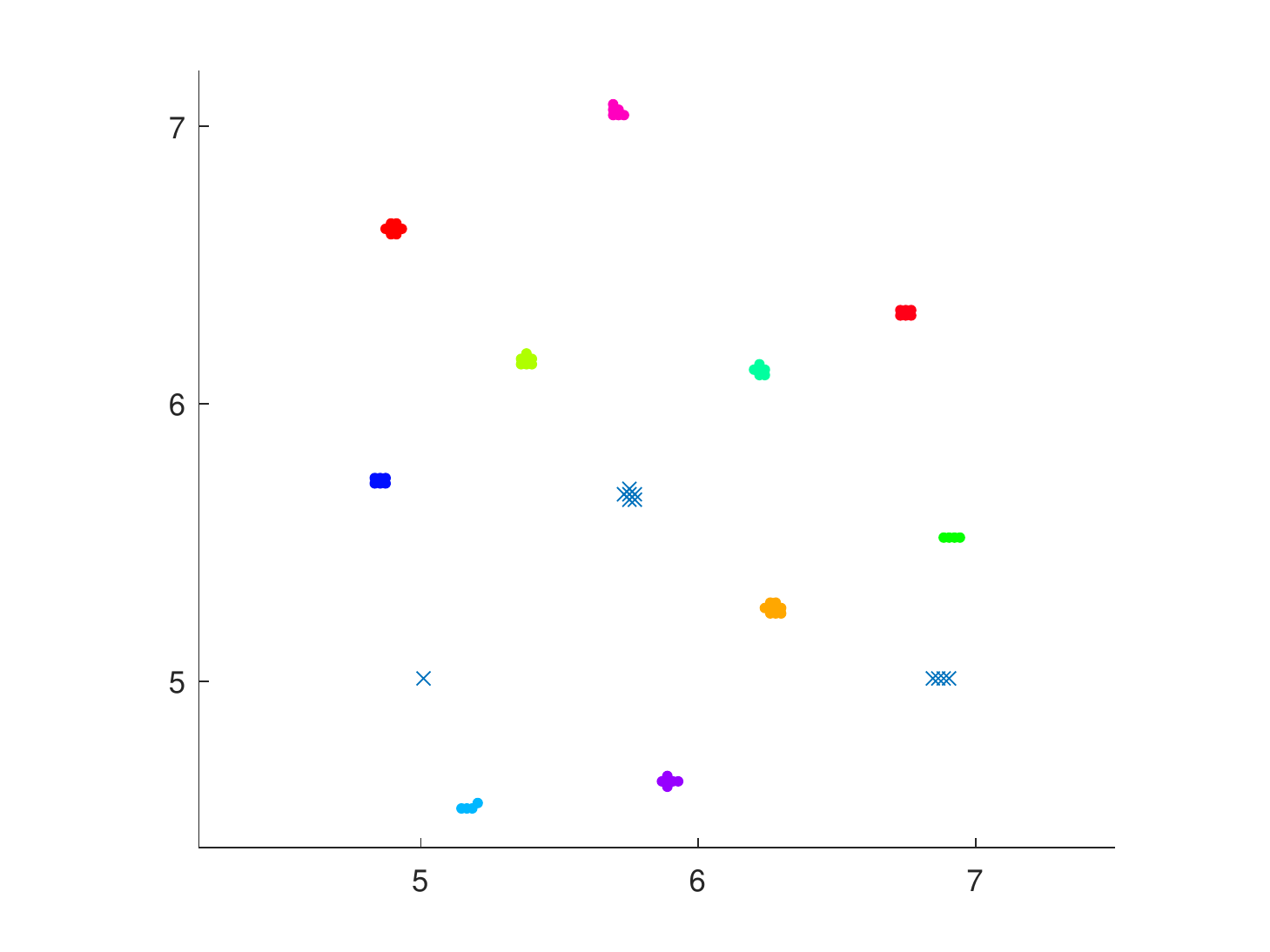}}}
  \caption{Last four layers of MMR output for sample instance of the asymmetric LJ potential, $N=13$. Different colors indicate different particles, and the blue crosses indicate the first three particles, which are specially constrained (though not completely fixed) to remove the degeneracy with respect to rigid motions.}\label{fig:8}
\end{figure}

In Fig. \ref{fig:8} we visualize the output of each of the last four levels of MMR as a scatter plot of all the grid points with 1-marginal larger than the threshold based on the solution of the SDP. Different colors are used to indicate different particles. At each level MMR identifies the points that are likely to contain the global optimal solution and propagates them to the next level.

Note that in principle, no two particles in any optimal configuration cannot be very close since the LJ potential tends to infinity when the pairwise distance approaches 0. However we can observe that at the coarse level $k=4$, the supports of different particles overlap significantly with each other. This observation highlights the advantage of including the upper bound $u^{(k)}$ in the SDP. The upper bound forces us to make more conservative (possibly overlapping) guesses for the particle positions so that we do not rule out the global optimizer before we reach the finer levels.

\par As MMR propagates to finer levels ($k=5$ for example), different particles start to separate from each other, and the solution support of each particle starts to take shape. We can observe the existence of multiple nearly optimal solutions in Fig.~\ref{fig:8}. For example, the dark blue particle has two possible locations. As we move to even finer levels ($k=6$ and $k=7$), though some of the possibilities are eliminated, we still have more support clusters than particles. Again, we can apply the methodology of Section~\ref{sec:multiple soln} to sample candidate solutions. For this particular example, the discrete solution attains an energy of \textbf{-90.2943}. By applying \texttt{fmincon}, we reach a final energy of \textbf{-90.4882}.

Next we compare the performance of MMR and BH on the asymmetric LJ potential. In MMR, we extract one solution using the method of Section~\ref{sec:multiple soln} and refine using \texttt{fmincon}. For BH we set the temperature to $T=1$ (tuned for best performance) and the maximum number of iterations to be $2000$, chosen such that BH runs for the same amount of time as MMR. Here we perform $10$ independent BH runs, and for each run the particles are independently and uniformly randomly initialized over $[0,5]^{2}$ (also tuned for best performance). This procedure is repeated over four independent realizations of the $r_{ij}$ defining the asymmmetric LJ potential. Results are summarized in Fig.\ref{fig:10}. A similar comparison for the $N=20$ case is reported in Appendix \ref{app:2}.

\begin{figure}[!htb]
\centering
\subfloat[Example 1]{{\includegraphics[width=0.20\textwidth]{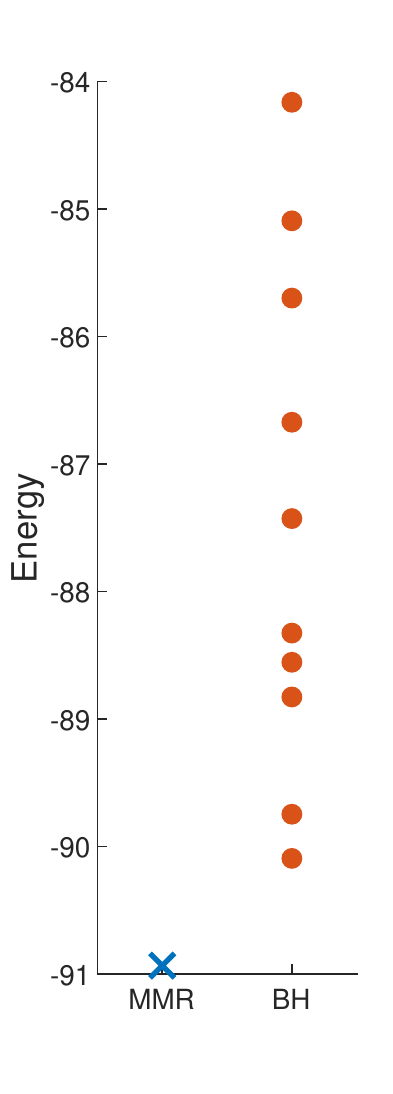}}}
\subfloat[Example 2]{{\includegraphics[width=0.20\textwidth]{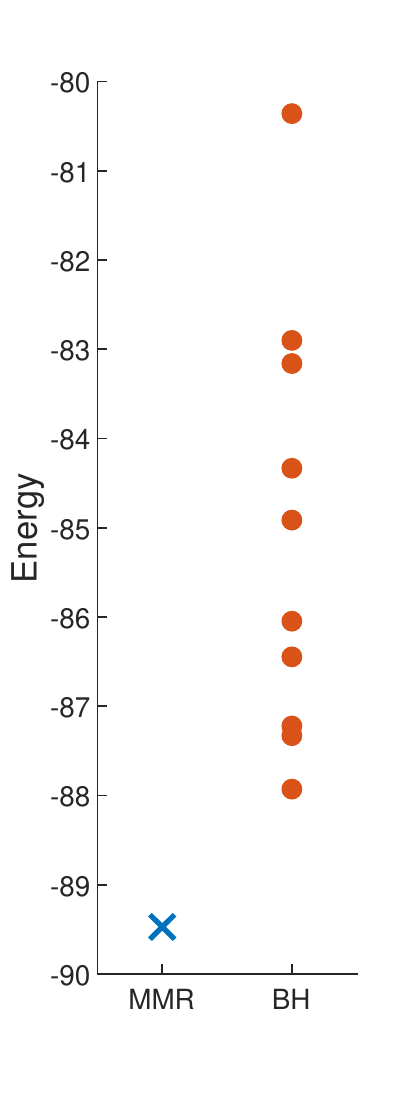}}}
\subfloat[Example 3]{{\includegraphics[width=0.20\textwidth]{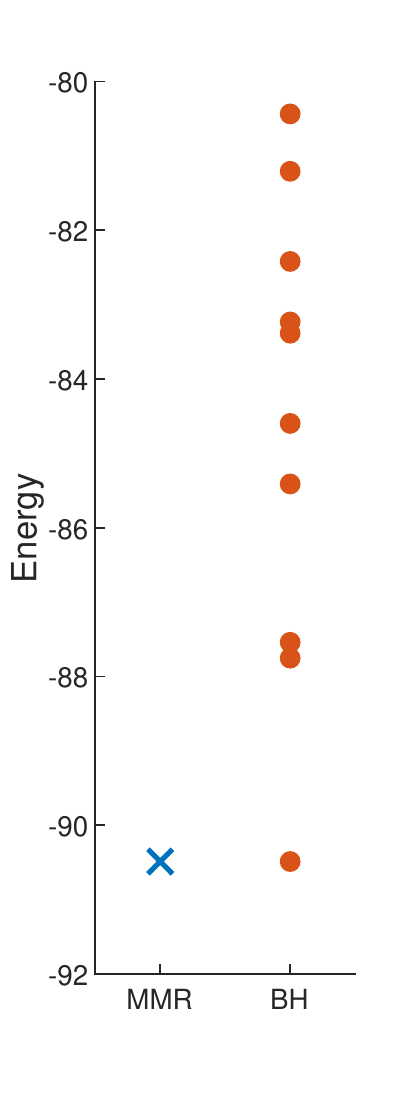}}}
\subfloat[Example 4]{{\includegraphics[width=0.20\textwidth]{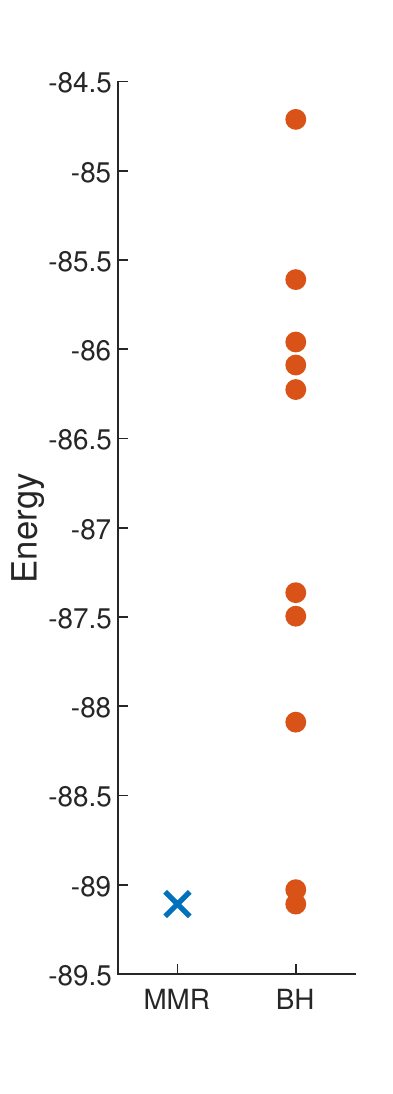}}}
  \caption{Energies of MMR and BH for $N=13$. The blue crosses represent the final MMR energies (after local optimization), and the red circles represent the energies of 10 independent BH tests with different initializations.}\label{fig:10}
\end{figure}

\section{Conclusion}

In this paper, we consider the global optimization of pairwise objective functions from the point of view of probability measure optimization. To avoid the exponential growth of the optimization space, we propose a relaxation and enforce semidefinite constraints on the 2-marginals. Then the relaxed problem can be solved by semidefinite programming. We embed a multiscale scheme into the SDP optimization framework in order to make the approach feasible for continuous state spaces. 

The key feature of the proposed algorithm is that it does not directly solve the global optimization problem, which includes a formidable semidefinite constraint of size $\sum_{i=1}^N \cdot |\mathcal{X_i}|$. Instead, it successively reduces the allowed support within the state space as it descends to finer discretizations. Compared with previous work on 2-marginal relaxations of continuum problems \cite{khoo2019convex}, our approach can achieve better time complexity via the multiscale framework. Furthermore, the algorithm is fundamentally global in that it starts by directly approaching a coarsened global optimization problem. In particular it does not rely on randomness to escape local optima like simulated annealing and basin-hopping~\cite{basinhopping}, though randomness can be used in a different way to explore candidate solutions as a post-processing step. Our numerical results for Lennard-Jones clusters indicate that our approach has more consistent results and can achieve better global energy, compared to the widely-used basin-hopping method, within the same budget of time.

When the pairwise objective function is measures the deviation between pairwise particle distances and some possibly noisy pairwise distance observations, the problem becomes a sensor network localization (SNL) problem. Compared with the SNLSDP algorithm~\cite{snlsdp}, our 2-marginal relaxation SDP is of larger size since it also relies on the discretization of the state space. However, we demonstrate our algorithm's ability to reconstruct certain globally rigid graphs instead of lifted higher-dimensional solutions provided by SNLSDP, for which dimension is not directly built into the algorithm. Our numerical results also indicate that MMR achieves much better accuracy and higher reconstruction probability when the measurements are noisy and incomplete.

The success of MMR on these two different types of problems is noteworthy, since existing state-of-the-art approaches are very different in nature. It is natural to then pursue the extension of MMR to other problems with similar structure.

However, note that the size of the SDP can still grow significantly when the size of the sensor network grows or when multiple near-optimal solutions exist (so that many points must be retained in the support at each level). The asymptotic complexity of the method may therefore depend on the problem. Moreover, several parameters are tuned carefully in this study, such as the upper bound, the threshold parameter, and the accuracy of the SDP solver. Improper parameter values may cause the computational complexity to grow dramatically on the one hand or may fail to capture all near-optimal solutions on the other. Though we present intuitions about their relationship, it is hard to devise a systematic way to choose these parameters.

Finally, in the case where multiple near-optimal solutions exist, MMR yields a 2-marginal resembling a convex combination of these solutions. In Section~\ref{sec:multiple soln} we propose an efficient spectral method for extracting solutions, but we remark that future work may uncover better ways to identify each convex component.

\newpage
\appendix
\section{Additional figures for symmetric LJ potential}\label{app:1}
\begin{figure}[!htb]
\subfloat[MMR+\texttt{fmincon} solution for $N=20$. Energy: $-95.1684$.]{{\includegraphics[trim=0.7cm 0.7cm 0.7cm 0.7cm, clip=true,width=0.23\textwidth]{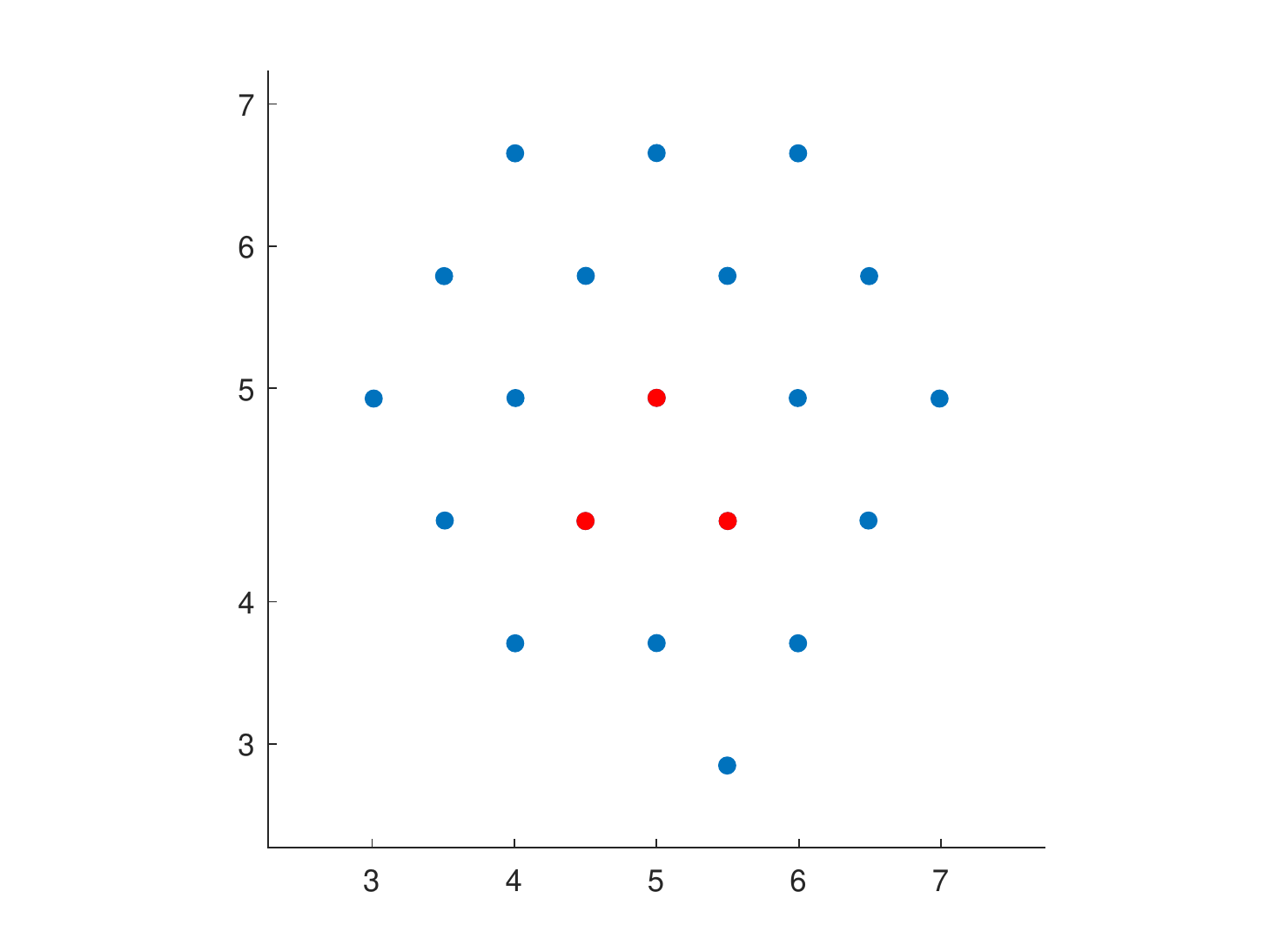}}}\ \ 
\subfloat[BH solution for $N=20$. Energy: $-95.1847$.]{{\includegraphics[trim=0.7cm 0.7cm 0.7cm 0.7cm, clip=true,width=0.23\textwidth]{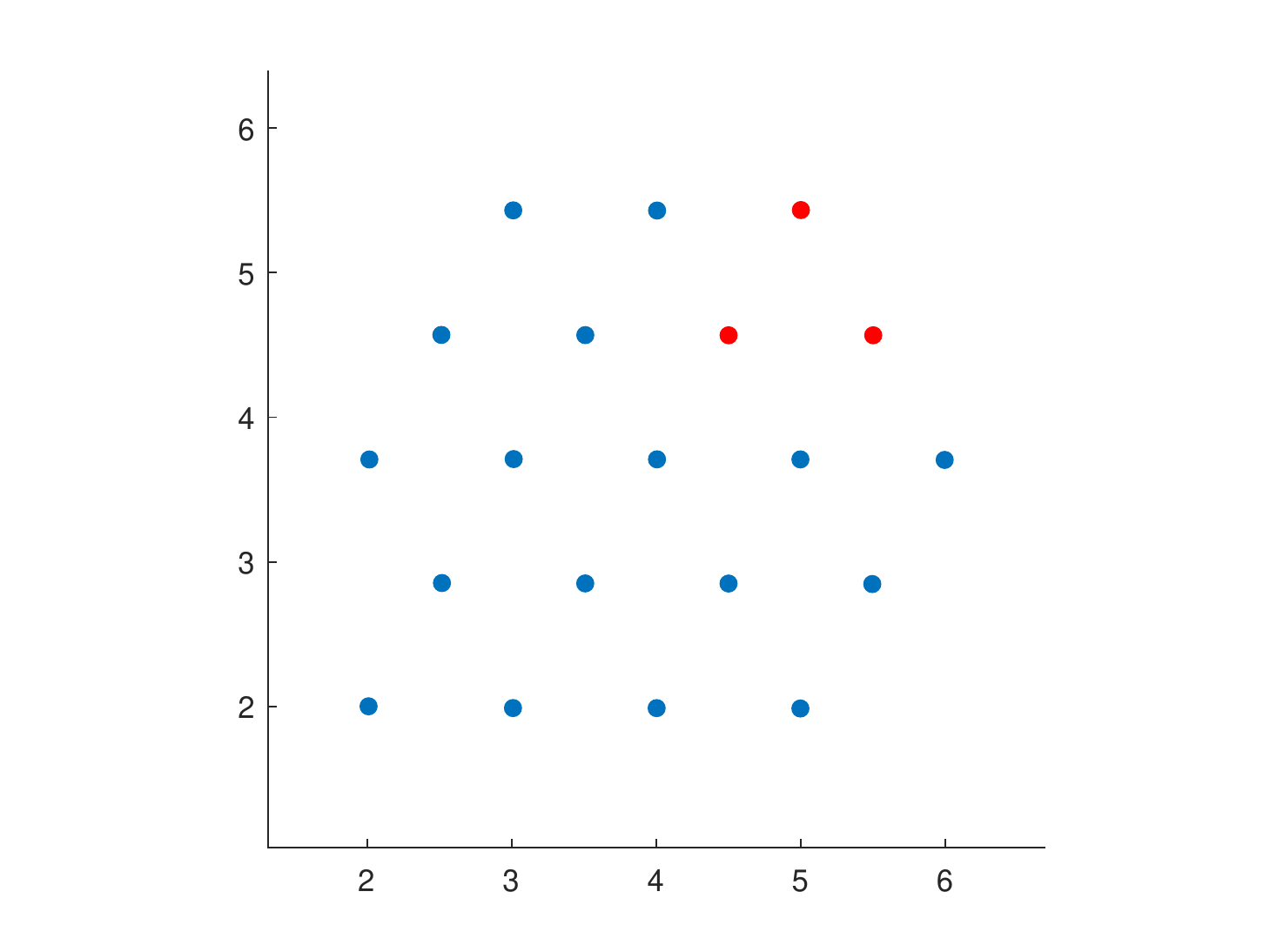}}}\ \ 
\subfloat[MMR+\texttt{fmincon} solution for $N=30$. Energy: $-154.7580$.]{{\includegraphics[trim=0.7cm 0.7cm 0.7cm 0.7cm, clip=true,width=0.23\textwidth]{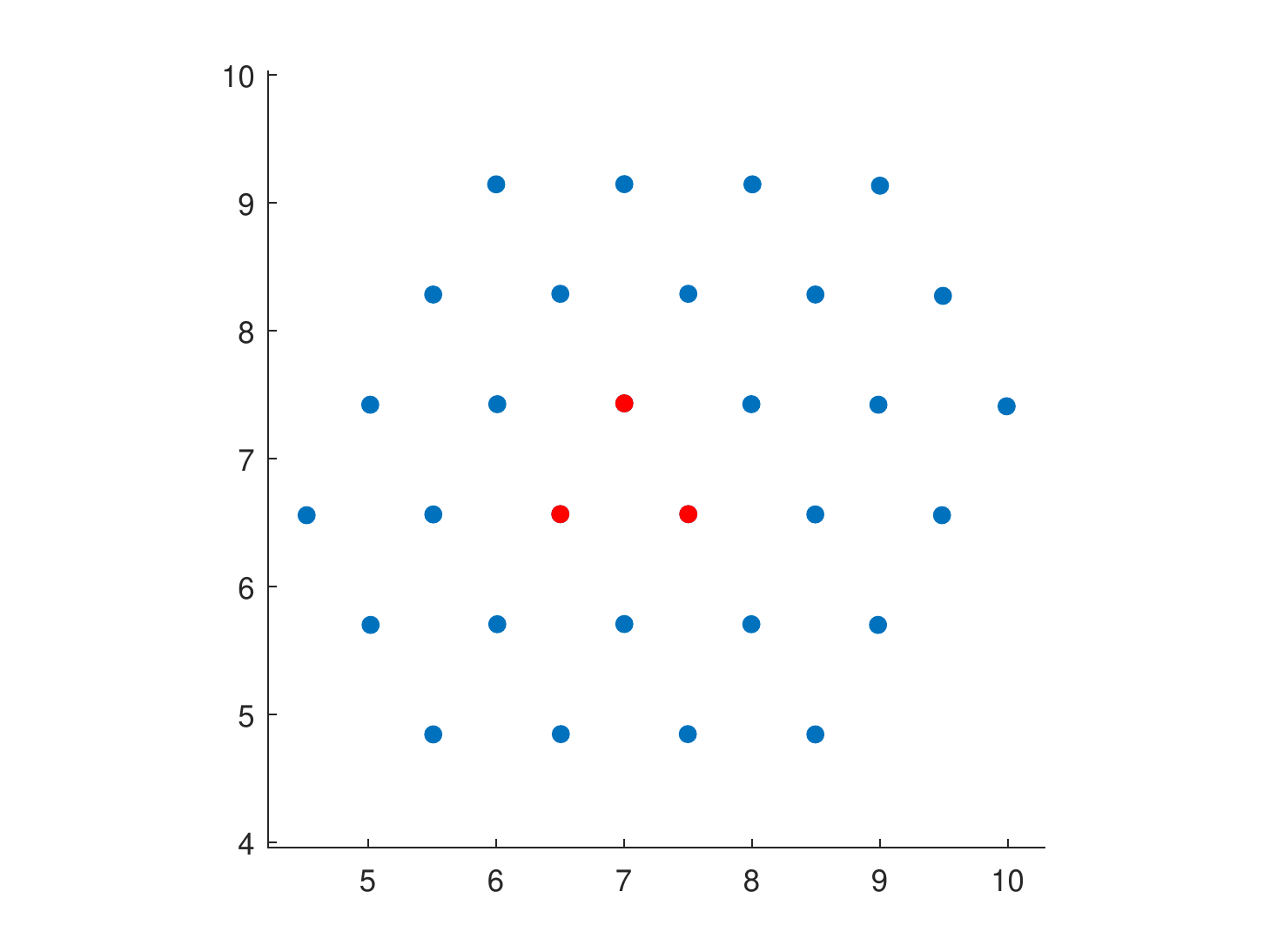}}}\ \ 
\subfloat[BH solution for $N=30$. Energy: $-154.7772$.]{{\includegraphics[trim=0.7cm 0.7cm 0.7cm 0.7cm, clip=true,width=0.23\textwidth]{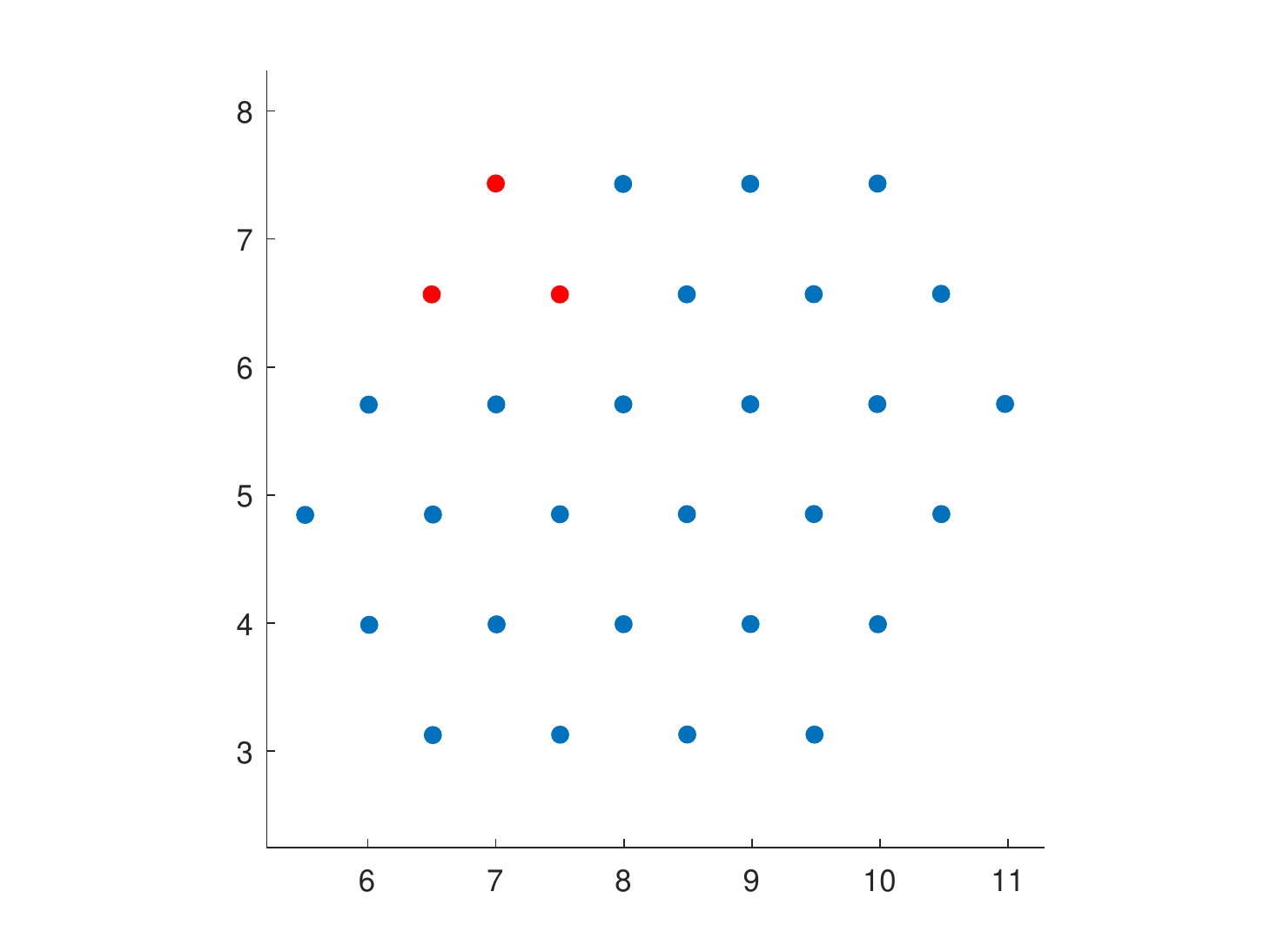}}}
  \caption{Optimal MMR and BH configurations. The red points indicate three fixed anchor particles. The configurations in (a) and (b) are congruent up to numerical error, as are the configurations in (c) and (d).}\label{fig:15}
\end{figure}

In Fig.~\ref{fig:15}, we plot the optimal MMR and BH configurations for $N=20,30$ in order to demonstrate that they are qualitatively the same, i.e., congruent.

\begin{figure}[!htb]
\subfloat[Energy: -95.0763 (\textbf{-95.0972}) ]{{\includegraphics[trim=0.7cm 0.7cm 0.7cm 0.7cm, clip=true,width=0.25\textwidth]{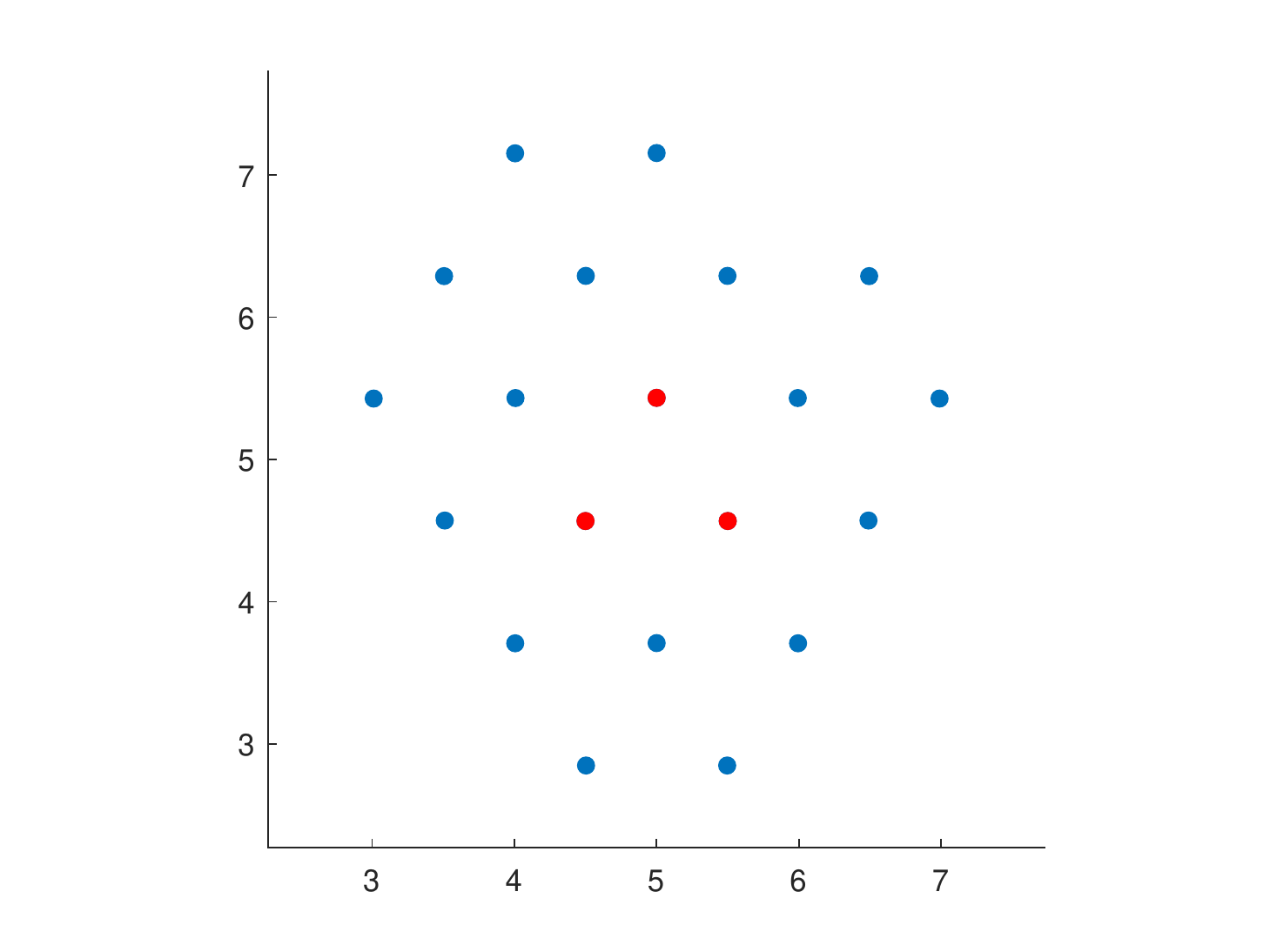}}}
\subfloat[Energy: -93.2332 (\textbf{-93.2541}) ]{{\includegraphics[trim=0.7cm 0.7cm 0.7cm 0.7cm, clip=true,width=0.25\textwidth]{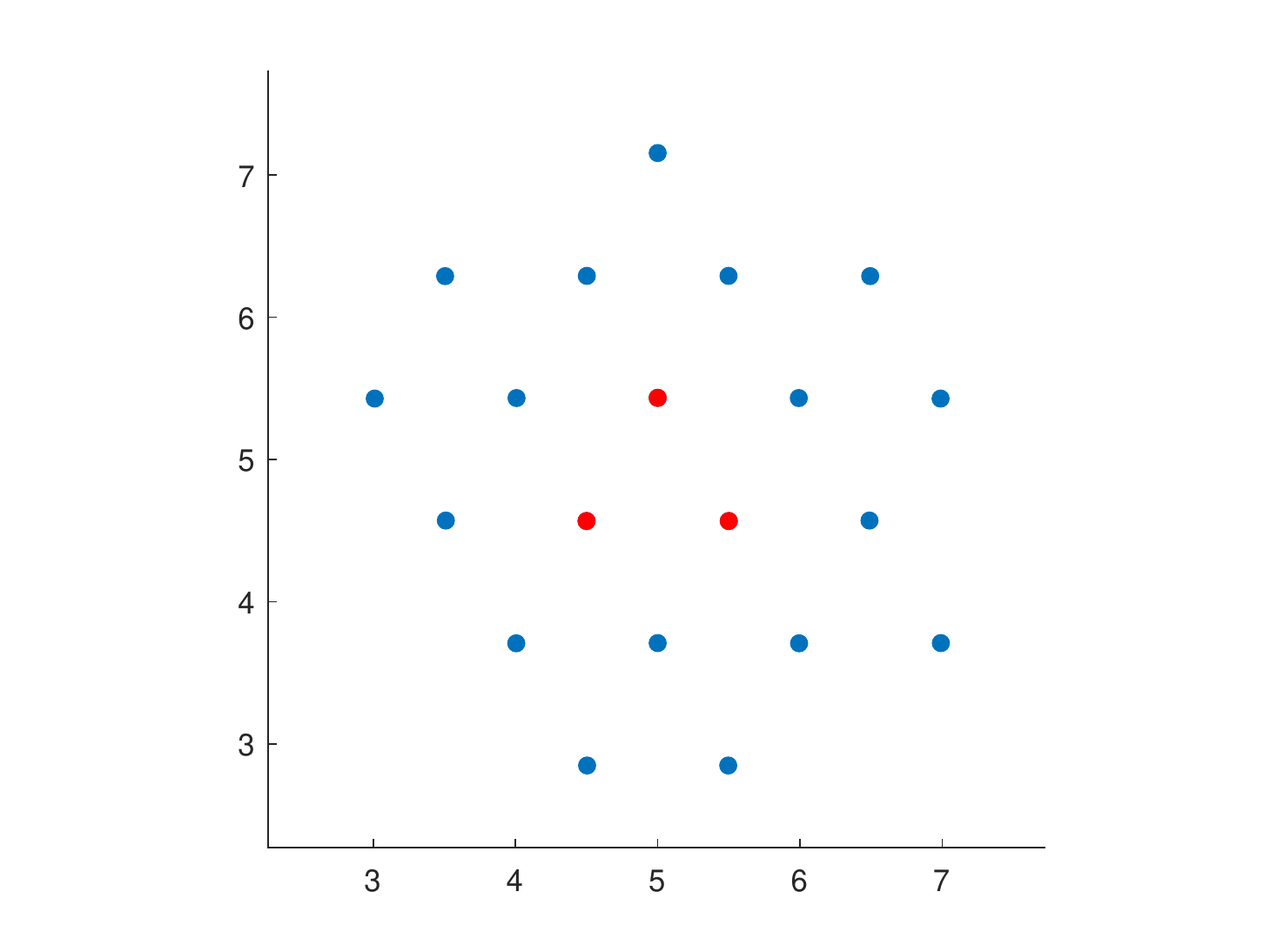}}}
\subfloat[Energy: -93.2348 (\textbf{-93.2493}) ]{{\includegraphics[trim=0.7cm 0.7cm 0.7cm 0.7cm, clip=true,width=0.25\textwidth]{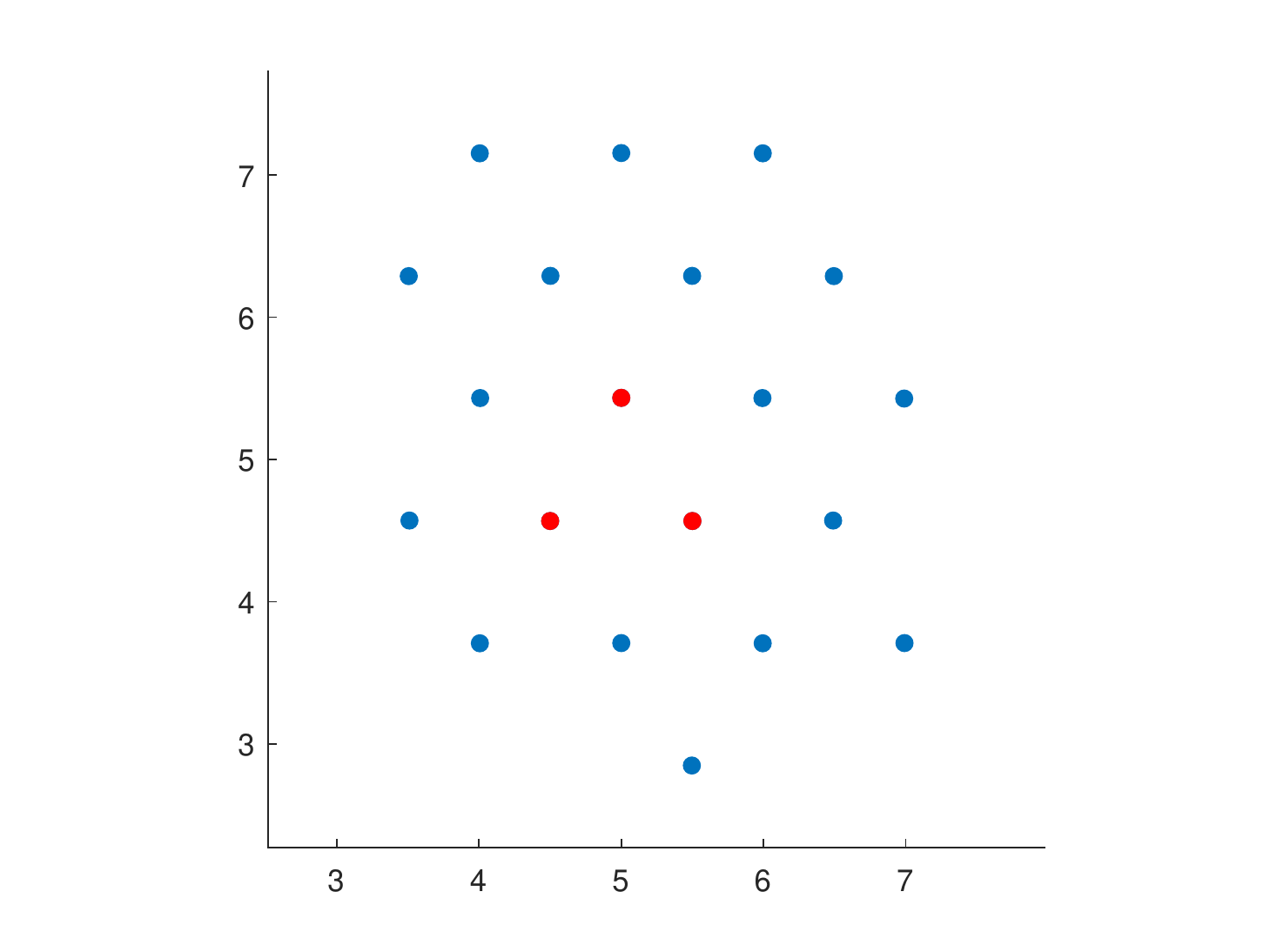}}}
\subfloat[Energy: -95.0763 (\textbf{-95.0972}) ]{{\includegraphics[trim=0.7cm 0.7cm 0.7cm 0.7cm, clip=true,width=0.25\textwidth]{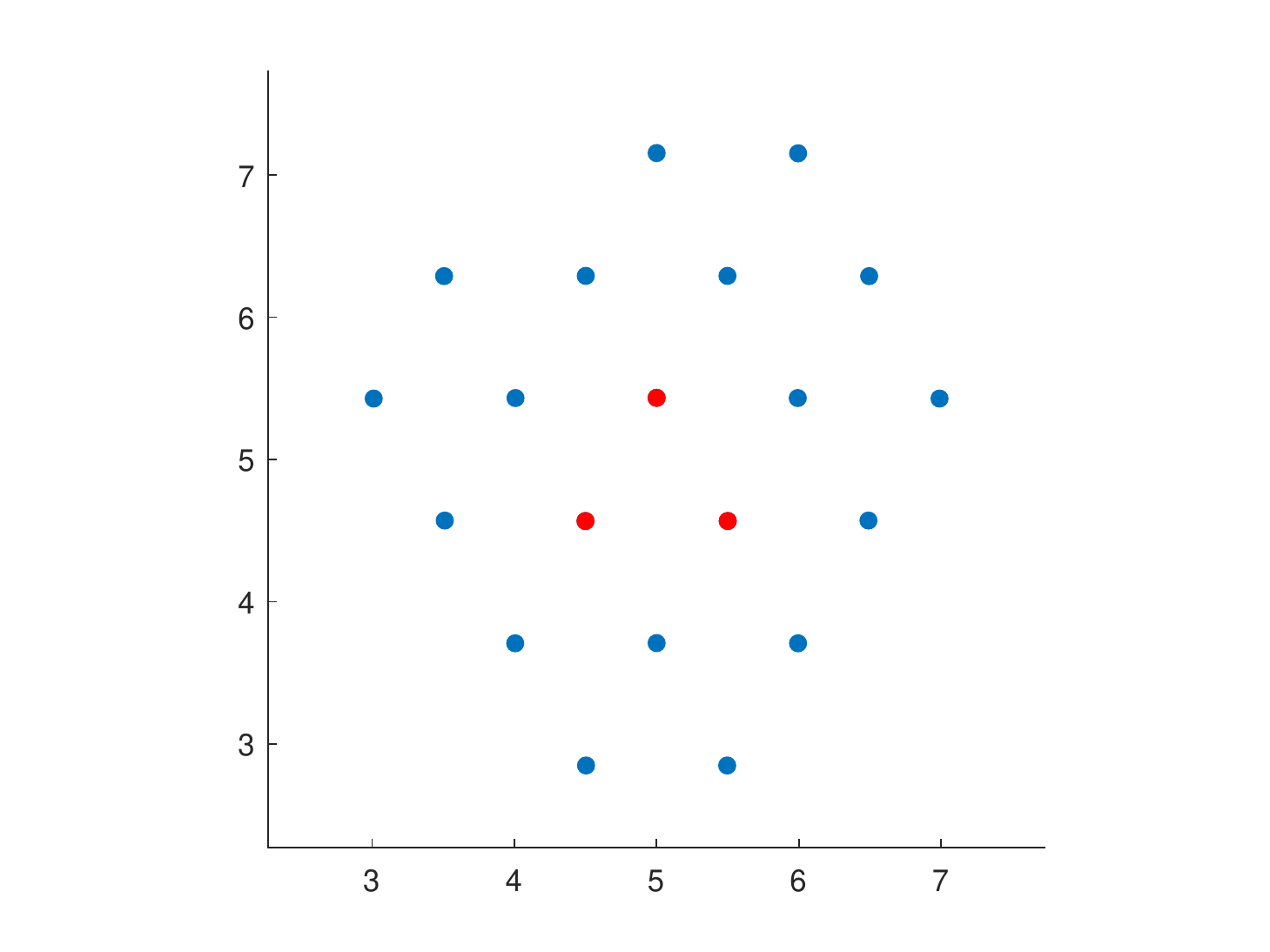}}}
  \caption{Four sample configurations via MMR for $N=20$. The red points indicate three fixed anchor particles. We also show the energy as well as the energy after refinement by \texttt{fmincon} (\textbf{bold face}) to remove discretization error.}\label{fig:12}
\end{figure}

\begin{figure}[!htb]
\subfloat[Energy: -150.8375 (\textbf{-150.9012}) ]{{\includegraphics[trim=0.7cm 0.7cm 0.7cm 0.7cm, clip=true,width=0.25\textwidth]{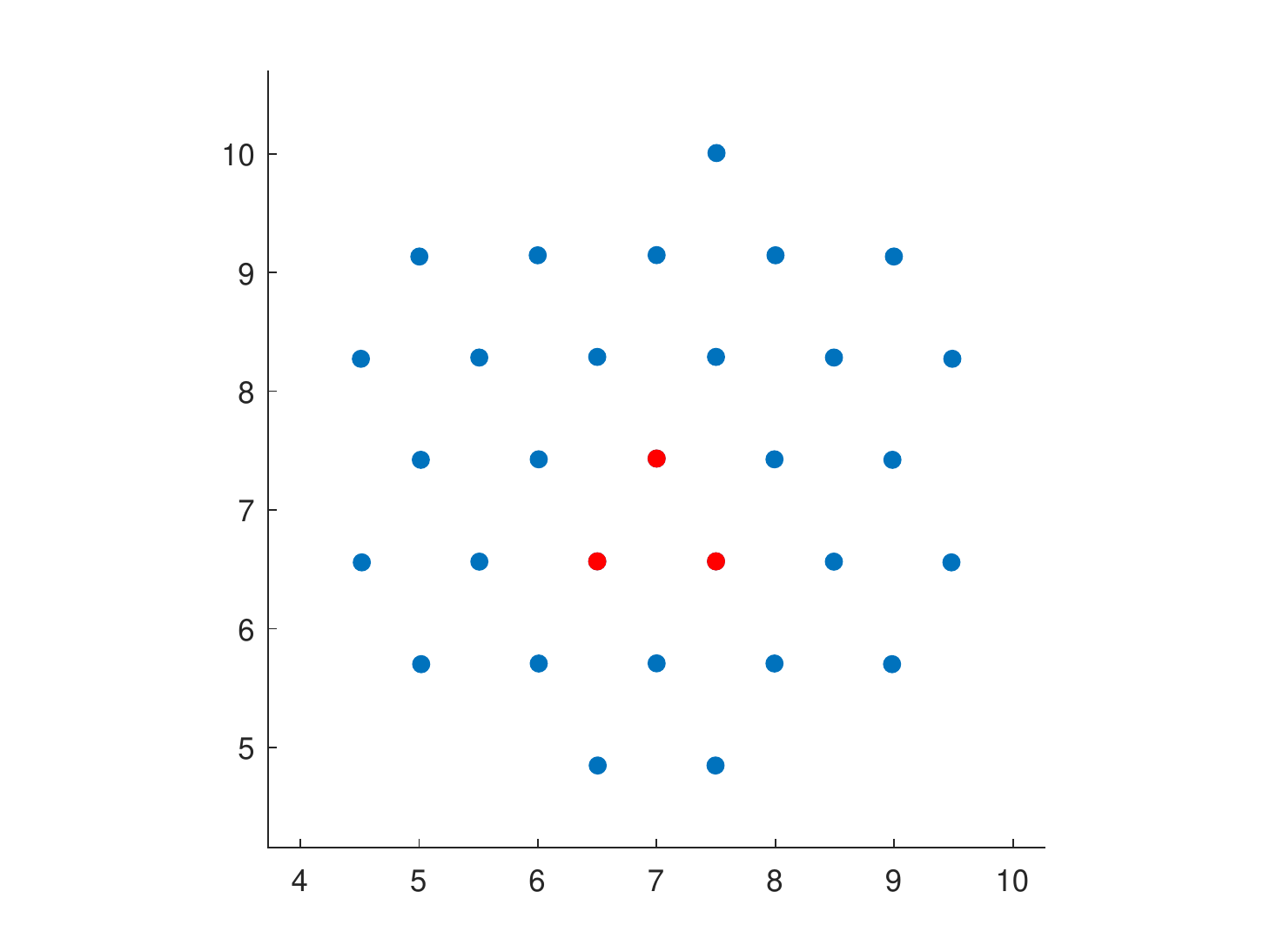}}}
\subfloat[Energy: -150.8707 (\textbf{-150.9460}) ]{{\includegraphics[trim=0.7cm 0.7cm 0.7cm 0.7cm, clip=true,width=0.25\textwidth]{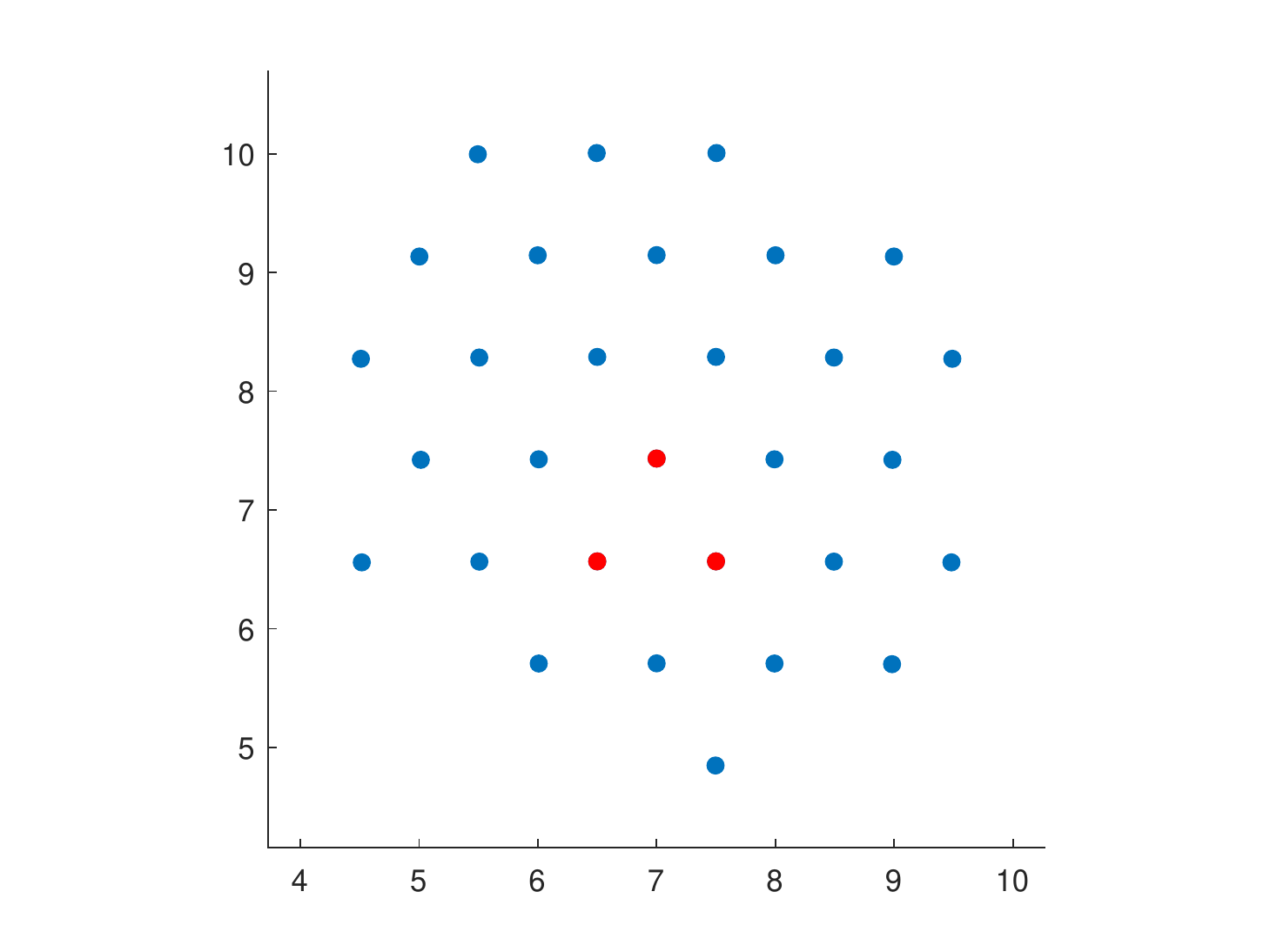}}}
\subfloat[Energy: -150.8274 (\textbf{-150.9010}) ]{{\includegraphics[trim=0.7cm 0.7cm 0.7cm 0.7cm, clip=true,width=0.25\textwidth]{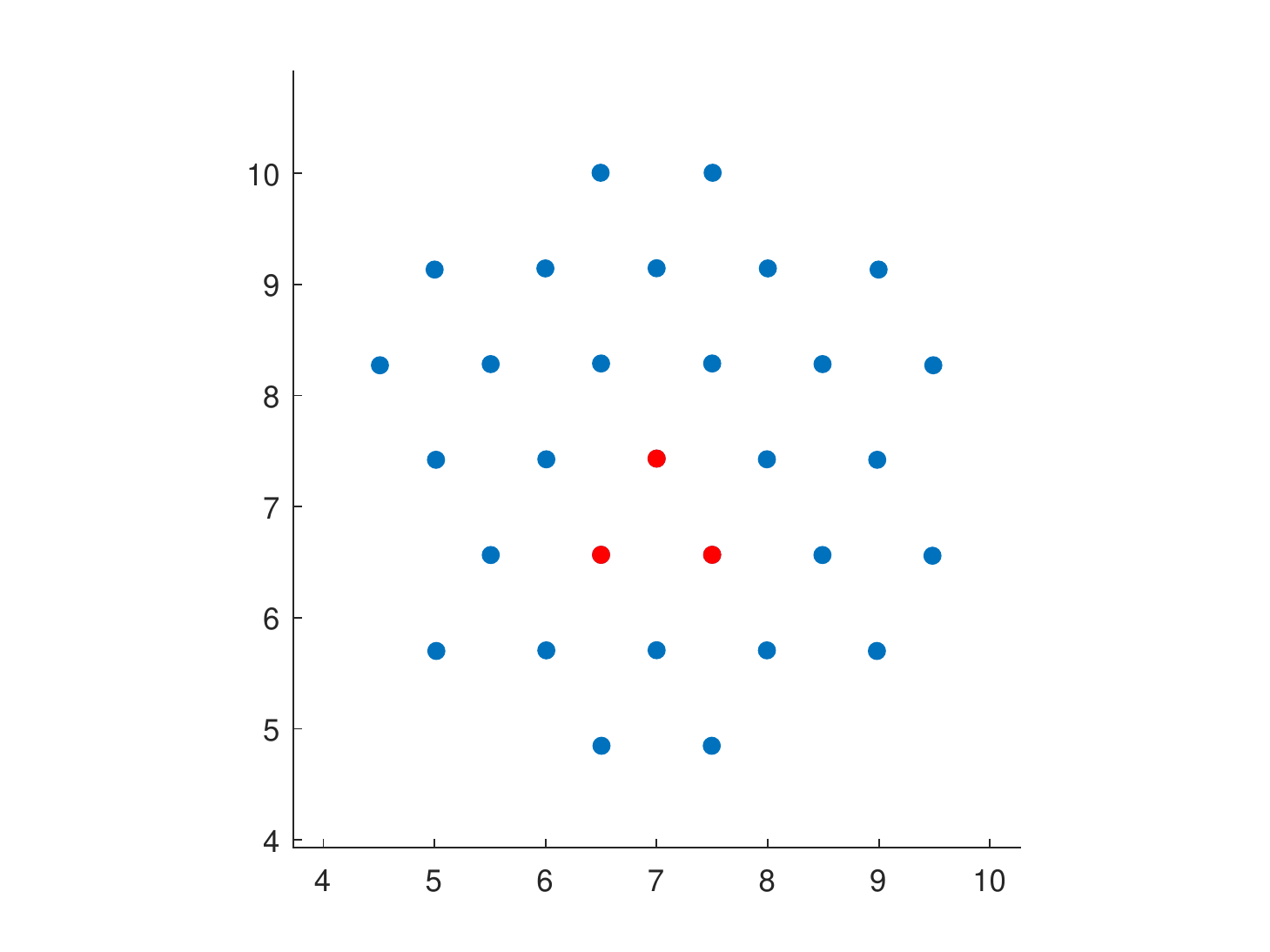}}}
\subfloat[Energy: -152.7218 (\textbf{-152.7930}) ]{{\includegraphics[trim=0.7cm 0.7cm 0.7cm 0.7cm, clip=true,width=0.25\textwidth]{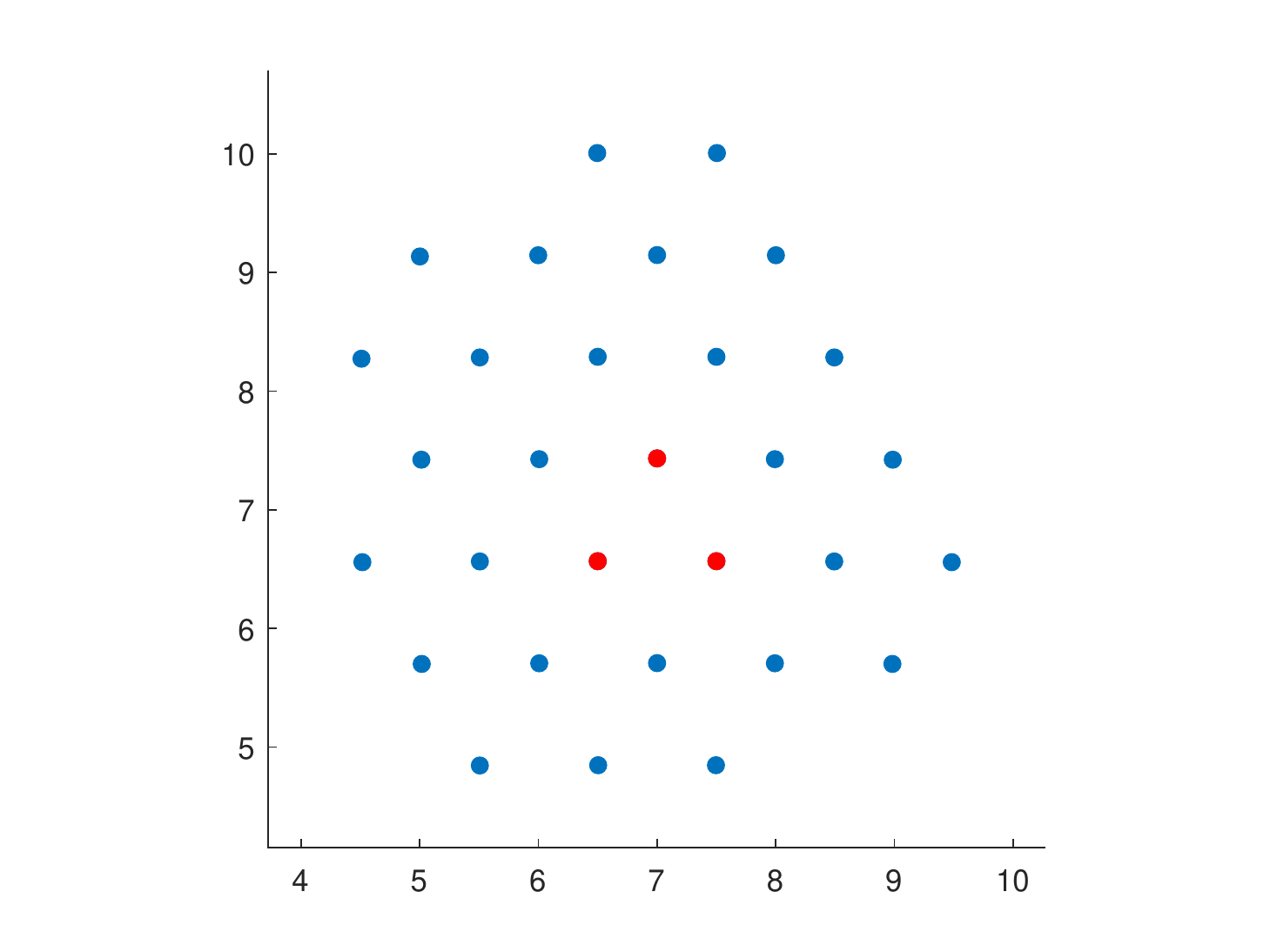}}}
  \caption{Four sample configurations via MMR for $N=30$. The red points indicate three fixed anchor particles. We also show the energy as well as the energy after refinement by \texttt{fmincon} (\textbf{bold face}) to remove discretization error.}\label{fig:13}
\end{figure}

Moreover, we use the randomized method of Section~\ref{sec:multiple soln} to explore various near-optimal configurations. We show four sample configurations each for the $N=20$ and $N=30$ cases in Figs.~\ref{fig:12} and~\ref{fig:13}, respectively. 

\section{Comparison of MMR and BH for asymmetric LJ potential, $N=20$}\label{app:2}
\begin{figure}[!htb]
\centering
\subfloat[Example 1]{{\includegraphics[width=0.20\textwidth]{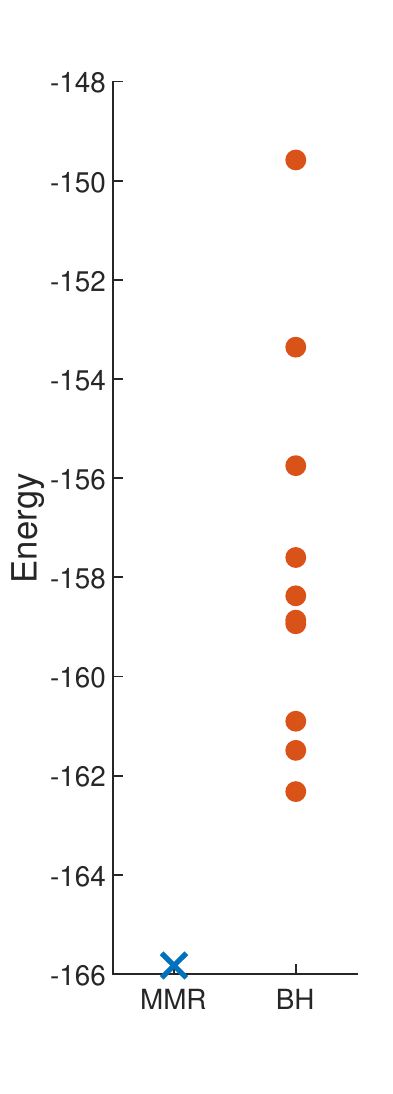}}}
\subfloat[Example 2]{{\includegraphics[width=0.20\textwidth]{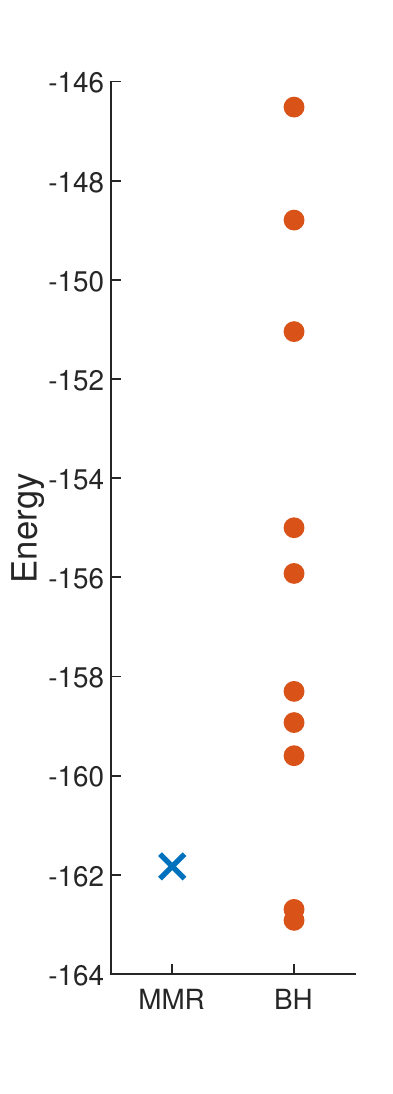}}}
\subfloat[Example 3]{{\includegraphics[width=0.20\textwidth]{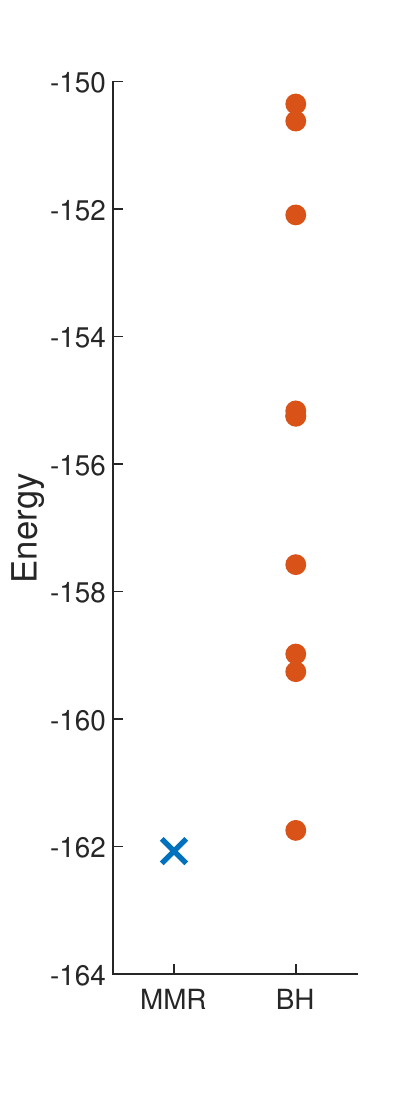}}}
\subfloat[Example 4]{{\includegraphics[width=0.20\textwidth]{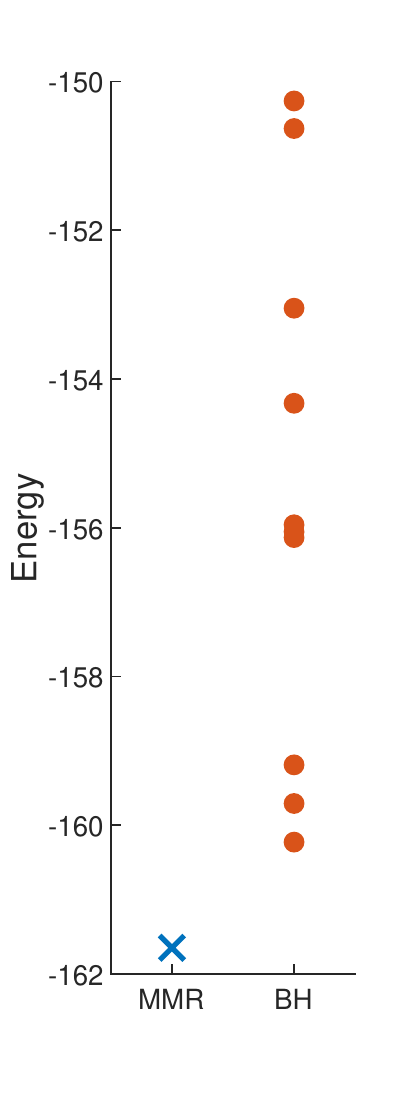}}}
  \caption{Energies of MMR and BH for $N=20$. The blue crosses represent the final MMR energies (after local optimization), and the red circles represent the energies of 10 independent BH tests with different initializations.}\label{fig:11}
\end{figure}

In Fig.~\ref{fig:11} we compare the performance of MMR and BH on the asymmetric LJ potential for $N=20$. For BH we set the temperature to $T=1$ and the maximum number of iterations to $5000$, chosen so that BH runs for the same length of time as the MMR algorithm. Here we perform $10$ independent BH runs, and for each run the particles are independently and uniformly randomly initialized over  $[0,5]^{2}$. This procedure is repeated over four independent realizations of the $r_{ij}$ defining the asymmmetric LJ potential.

\newpage
\bibliographystyle{siam}
\bibliography{mmr}

\end{document}